%% file: Voronoi.tex
\newcommand\webcite[1]{\texttt{\def~{\~{}}#1}}

\def\tZ{Z^*}
\def\tcol{\col^*}
\def\pmax{p_{\max}}
\def\short{{L}}
\def\shortm{{\widetilde L}}
\def\ymax{{y_+}}
\def\ymin{{y_-}}
\def\ymaxmin{{y_\pm}}
\def\etatwo{\eta}

\documentclass{article}
\usepackage{amsmath,amsfonts,amsthm,epsfig,color}
\newtheorem{theorem}{Theorem}
\newtheorem{lemma}[theorem]{Lemma}

\newcounter{claim_base}
\newtheorem{claim}{Claim}[claim_base]

\newtheorem{corollary}[theorem]{Corollary}

\def\fun{t}
\def\Prone{{\widetilde\Pr}}
\def\Prtwo{{\widetilde\Pr}}

\def\chiN{{\chi}}
\def\chiA{{\chi^{\mathrm{area}}}}
\def\dist#1#2{{\overline{#1#2}}}
\renewcommand\epsilon{\varepsilon}
\def\whp{{\bf whp}}
\newcommand\E{{\mathop{\mathbb E{}}\nolimits}}
\renewcommand\Pr{{\mathop{\mathbb P{}}\nolimits}}
\newcommand\Bi{{\mathop{\mathrm{Bi}}\nolimits}}

\newcommand\area{{\mathop{\mathrm{area}}\nolimits}}
\newcommand\diam{{\mathop{\mathrm{diam}}\nolimits}}

\def\Es{{E_{\rm dense}}}
\def\pbad{p_{\mathrm{bad}}}
\def\pneut{p_{\mathrm{neutral}}}
\def\pgood{p_{\mathrm{good}}}
\def\pbadn#1{p_{{\mathrm{bad}},#1}}
\def\pgoodn#1{p_{{\mathrm{good}},#1}}
\def\col{{\rm col}}

\def\Egl{{E_{\rm good}}}
\def\Etdc{{E_3^{\rm crude}}}
\def\Etd{{E_3^{\rm disc}}}

\def\CL{{\mathrm{C}}}
\def\cA{{\mathcal A}}
\def\cB{{\mathcal B}}

\def\Pow{{\mathcal P}}
\def\Z{{\mathbb Z}}
\def\RR{{\mathbb R}}
\def\TT{{\mathbb T}}

\begin{document}
\title{The critical probability for random Voronoi percolation in the plane is $1/2$}
\date{September 27, 2005}

\author{B\'ela Bollob\'as\thanks{Department of Mathematical Sciences,
University of Memphis, Memphis TN 38152, USA}
\thanks{Trinity College, Cambridge CB2 1TQ, UK}
\thanks{Research supported in part by NSF grant ITR 0225610 and DARPA grant
F33615-01-C-1900}
\thanks{Research partially undertaken during a visit to the Forschungsinstitut f\"ur Mathematik,
ETH Z\"urich}
\and Oliver Riordan$^{\dag\S}$%
\thanks{Royal Society Research Fellow, Department of Pure Mathematics
and Mathematical Statistics, University of Cambridge, UK}}
\maketitle

\begin{abstract}
We study percolation in the following random environment: let $Z$ be a Poisson process
of constant intensity on $\RR^2$, and form the Voronoi tessellation of $\RR^2$ with respect
to $Z$. Colour each Voronoi cell black with probability $p$, independently of the other cells.
We show that the critical probability is $1/2$. More precisely,
if $p>1/2$ then the union of the black cells contains an infinite component with probability $1$,
while if $p<1/2$ then the distribution of the size of the component of black cells containing a given point decays
exponentially. These results are analogous to Kesten's results for bond percolation in $\Z^2$.

The result corresponding to Harris' Theorem for bond percolation in $\Z^2$ is known:
Zvavitch noted that one of the many proofs of this result can easily be adapted
to the random Voronoi setting. For Kesten's results, none of the existing proofs seems to adapt.
The methods used here also give a new and very simple
proof of Kesten's Theorem for $\Z^2$; we hope they will
be applicable in other contexts as well.
\end{abstract}

\setcounter{tocdepth}{1}
\tableofcontents

\begin{figure}[htb]
 \[\epsfig{file=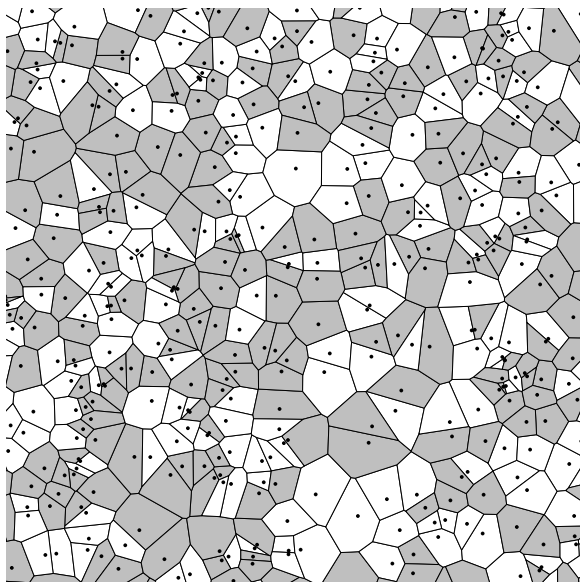,height=3in}\]
\caption{Critical random Voronoi percolation.}
\label{fig_CVP}
\end{figure}

\section{Introduction and results}\label{sec_thms}

Percolation theory is a very active area, with many difficult basic questions
still unanswered. Perhaps the most basic
question is to establish the critical probability above which
percolation occurs; this question was posed by Broadbent and
Hammersley~\cite{BH} in 1957 in a wide variety of contexts. In
general, it seems impossible to answer: the exact value is known
only in a small number of cases in which something
very special happens, involving duality. 
The best known case is Kesten's celebrated result
that for bond percolation in $\Z^2$, which is self-dual, the critical
probability is $1/2$ (more about this later). For other exact results,
see Grimmett~\cite{Grimmett}. Another very natural percolation process
for which
self-duality occurs is {\em random Voronoi percolation}, introduced in
the context of first-passage percolation by Vahidi-Asl and Wierman~\cite{VW}; we come to
the formal description in a moment. Random Voronoi percolation has
been studied by many people; see, for example \cite{VW,VW2,VW3,Zv,FreedmanProj,Aizenman,BS,BBQ}.
Unlike the classical lattice examples, here the environment in which percolation occurs is
itself random; as we shall see, this means that the techniques needed to
establish the critical probability are rather different.

The formal set-up is as follows.
A certain parameter $p$ is given. We construct a Poisson process
$Z$ with intensity 1 on $\RR^2$, and, given $Z$, assign to each point $z$ of $Z$ independently a colour $\col(z)
\in\{\hbox{black, white}\}$,
with $\Pr(\col(z)=\hbox{black})=p$. Equivalently, we may generate $Z$ as the union
of two Poisson processes $Z_b$, $Z_w$ with intensities $p$, $1-p$, corresponding to the black and white points,
i.e., to $\{z\in Z:\col(z)=\hbox{black}\}$ and $\{z\in Z:\col(z)=\hbox{white}\}$.
Throughout, we shall write $\Pr_p^{\RR^2}$, or simply $\Pr_p$, for the associated probability
measure. 
We use colouring terminology, as in~\cite{Smirnov}, for example,
rather than the more common `open or closed' terminology;
this emphasizes the symmetry and is more natural for figures.

Given $Z$, we construct the associated Voronoi tiling: the {\em Voronoi cell} of $z\in Z$ with respect to $Z$
is the set
\[
 V(z) = V_Z(z)= \{x\in \RR^2:\: d(x,z)=\inf_{z'\in Z}d(x,z')\},
\]
where $d(x,y)$ is the Euclidean metric.
\begin{figure}[htb]
 \[\epsfig{file=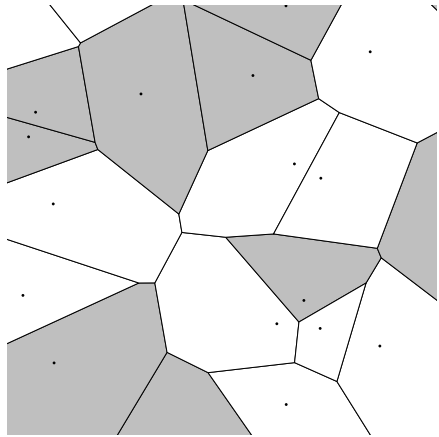,height=2.25in}\]
\caption{Part of the black (grey) and white-coloured random Voronoi tiling in $\RR^2$. The dots
are the points of $Z$.}
\label{fig_Vp}
\end{figure}

As we shall note below, with probability 1, every cell is a closed convex $k$-gon for
some $k$, every vertex of every cell lies in exactly three cells,
and if two cells meet, they share an entire edge. Note, however, that cells
$V(z_1)$ and $V(z_2)$ may meet without their common edge containing the midpoint
of $z_1z_2$. (There are several examples in Figure~\ref{fig_Vp}.)
The {\em graph} of the Voronoi tiling is the countably infinite graph $G$ with $Z$
as vertex set, in which two points $z_1$, $z_2$ are adjacent if their Voronoi cells meet.
The notions of percolation on the Voronoi tiling will be the usual notions of site percolation
on the (random) graph $G$.

We take the colour of a Voronoi cell $V_Z(z)$ to be the colour of the point $z$, and say that
each $x\in V_Z(z)$ also has this colour. (Thus, a point $x\in \RR^2$ may be both black and white
if it is in the boundary of two (or three) Voronoi cells.)
A {\em path of cells} corresponds to a path in $G$.
Let $z_0$ be the (with probability 1 unique) point $z_0\in Z$ such that the origin lies in $V_Z(z_0)$.
In the graph $G$, the {\em black vertex cluster containing $v_0$} is the set $C_0^G$ of all black vertices of $G$
connected to $v_0$ by a path consisting of black vertices.
Correspondingly, the {\em  black component of the origin} $C_0$
is the maximal connected set of black points
in $\RR^2$ containing the origin, i.e., the union of the cells $V_Z(z)$, $z\in C_0^G$.

Writing $\area(\cdot)$ and $\diam(\cdot)$ for the usual geometric area and diameter of a set in $\RR^2$,
consider the following three conditions: $|C_0^G|=\infty$, $\area(C_0)=\infty$ and $\diam(C_0)=\infty$.
It is easy to see that, except on a set of measure zero, either all three of these conditions hold or none.
Let 
\[
 \theta(p)= \Pr_p(|C_0^G|=\infty).
\]
Equivalently, $\theta(p)$ may be defined in terms of $C_0$. We say that
the random Voronoi percolation with parameter $p$ {\em percolates} if $\theta(p)>0$.
As translationally invariant events have probability $0$ or $1$, 
it is easy to show that if $\theta(p)>0$
then with probability $1$ there is an infinite black component somewhere in $\RR^2$,
while if $\theta(p)=0$ then with probability $1$ there is no such component.

As $\theta(p)$ is increasing, there is a critical probability
\[
 p_H = \inf\{p:\theta(p)>0\}.
\]
Here, following Welsh (see~\cite{SW}), the $H$ is in honour of Hammersley.

There is another natural definition of critical probability.
Writing $\E_p$ for the expectation associated to $\Pr_p$, let
\[
 \chiN(p) = \E_p |C_0^G|,
\]
and let
\[
 \chiA(p) = \E_p \area(C_0).
\]
Of course, we cannot determine the exact relationship between $\chiN(p)$ and $\chiA(p)$.
However, it will follow from our results that, unsurprisingly,
one is finite if and only if the other is. As $\chiN(p)$ is increasing in $p$,
there is another critical probability,
\[
 p_T = \inf\{p : \chiN(p)=\infty\},
\]
where the $T$ is in honour of Temperley.

The quintessential percolation problem is the study of {\em bond percolation in $\Z^2$}.
Here $\Z^2$ is the planar square lattice, i.e., the graph with vertex set $\Z^2$ in which
two vertices are adjacent if they are at Euclidean distance $1$. Each edge of this graph is taken
to be {\em open} with probability $p$ independently of all other edges. Edges
that are not open are {\em closed}. The {\em open cluster} containing a vertex $v$
is the maximal connected subgraph of $\Z^2$ that contains $v$, all of whose edges are open.
In this context, one can define $p_H(\Z^2)$ and $p_T(\Z^2)$ as above.

In 1957, Broadbent and Hammersley~\cite{BH} posed the problem of determining $p_H$
in a wide variety of contexts, including bond percolation in $\Z^2$. Hammersley~\cite{H2,H4,H5}
proved general results implying in particular that $0.35<p_H(\Z^2)<0.65$.

As the square lattice $\Z^2$ is self-dual, if $p=1/2$ and
$R$ is an $n+1$ by $n$ rectangle then the probability that there
is a path of open edges joining some vertex on the left-hand side of $R$ to some vertex on
the right-hand side is {\em exactly} $1/2$. This observation suggests very strongly that
$p_H(\Z^2)=1/2$. Nevertheless, proving this turned out to be highly non-trivial. The first major
progress was made by Harris~\cite{Harris} in 1960, when he proved that $p_H(\Z^2)\ge 1/2$.
Harris' proof used self-duality as a starting point, but this trivial observation is {\em just} a starting point;
several substantial ideas were needed to deduce the result.
From this point it became
a well-known problem to show that $p_H(\Z^2)=p_T(\Z^2)=1/2$. Indeed, by the time Kesten~\cite{Kesten1/2}
finally proved this result 20 years later, the problem was sufficiently well known that Kesten
quite rightly felt it needed no introduction: the abstract of his paper, whose title is `The critical probability
of bond percolation on the square lattice equals $1/2$', reads `We prove the statement in the title of the paper'!

The random Voronoi percolation we consider here also has a `self-duality' property which
implies that for $p=1/2$ the appropriately defined crossing probability for a square
is exactly $1/2$. In the
light of this property and the Harris-Kesten results, it is extremely natural to conjecture
that for the Voronoi set-up, $p_T=p_H=1/2$ holds also.
By adapting one of the many known proofs of Harris' Theorem,
Zvavitch~\cite{Zv} has shown that $\theta(1/2)=0$, so $p_H\ge 1/2$. The proof Zvavitch adapted is due
to Burton and Keane~\cite{BurtonKeane} and Zhang; see~\cite[pages 198 and 289]{Grimmett}.
Other known proofs of Harris' Theorem do not seem to adapt in this way.

Here we shall prove the analogue of Kesten's results for bond percolation in $\Z^2$, showing that
$p_T=p_H=1/2$ holds for random Voronoi percolation. Although
no known proof of Kesten's results seems to translate to the random
Voronoi context, the methods used in this paper do give a new 
and very simple proof of Kesten's Theorem, which we shall describe in~\cite{ourKesten}.

\begin{theorem}\label{th_doesperc}
For random Voronoi percolation in the plane, $p_H\le 1/2$.
In other words, if $p>1/2$ then $\theta(p)>0$.
\end{theorem}

In the subcritical case, we establish exponential decay of the `volume' of the black component
containing a given point. Here we may take $|C_0|$ to be $\area(C_0)$, $\diam(C_0)$, or
the number $|C_0^G|$ of Voronoi cells making up $C_0$.

\begin{theorem}\label{th_decay}
For any $p<1/2$ there is a constant $c(p)>0$ such that
\[
 \Pr_p (|C_0|\ge n) \le \exp(-c(p)n)
\]
for every $n\ge 1$. In particular, $p_T\ge 1/2$.
\end{theorem}

As $p_T\le p_H$, it follows that $p_T=p_H=1/2$.

One of our intermediate results also gives an alternative proof of Zvavitch's result that
$\theta(1/2)=0$; we state this as a theorem for ease of reference. In principle, our method
gives a bound on the rate of decay of $\Pr_{1/2}(|C_0^G|\ge n)$ as $n\to\infty$, but this
bound would be rather weak and we have not calculated it explicitly.

\begin{theorem}[Zvavitch \cite{Zv}]\label{th_doesnt}
For random Voronoi percolation in the plane, $\theta(1/2)=0$.
\end{theorem}

The rest of the paper is organized as follows. In Section~\ref{sec_ingred} we collect
the results of general probabilistic combinatorics we shall use, and in Section~\ref{sec_prelim} we present
some basic facts about random Voronoi percolation.
Our first intermediate result, Theorem~\ref{th_RSW},
which will play the role of the Russo-Seymour-Welsh (RSW) Theorem of~\cite{Russo}
and~\cite{SW}, is stated and proved in Section~\ref{sec_RSW}.

In the remaining sections, our aim is to apply a sharp-threshold result of Friedgut and Kalai~\cite{FK}; see
Section~\ref{sec_ingred}. In order to do this, we work with 
random Voronoi percolation not in $\RR^2$, but in a torus,
introduced in Section~\ref{sec_torus}. We need to approximate the relevant Poisson process by
a discrete process: preliminary lemmas for doing this are given in Section~\ref{sec_approx}. The actual
application of the sharp-threshold result is given in Section~\ref{sec_above}; the proofs
of the main results follow in Section~\ref{sec_proofs}. In the course of the paper, for
technical reasons we need several deterministic lemmas that are not terribly appealing. The proofs of these
are postponed to the appendix.

There are two parts of the overall proof that are longer than one might expect. The first
is our proof of an analogue of the RSW Theorem, given in Section~\ref{sec_RSW}.
For bond percolation in $\Z^2$, several proofs of the RSW Theorem
are now known, and it is possible to give a very short proof. However, none of the existing proofs carries
over to the random Voronoi setting, which lacks an important independence property needed. This is discussed
further in Section~\ref{sec_RSW}.

Secondly, to apply the Friedgut-Kalai sharp-threshold result, we need to approximate a Poisson process
with a discrete probability space in some way. We work with a torus $T$, which has finite area, rather than $\RR^2$.
Of course, one can `discretize' a Poisson process on $T$, essentially by rounding the coordinates of all points
to multiples of a small quantity $\delta$. However, there will be a limit as to how small we can take $\delta$
in terms of the area of $T$, so this approximation will introduce `defects',
i.e., places where the discretized Poisson process does not tell us which of the cells in the Voronoi
tiling meet. There are not many defects, and, intuitively, it is clear that these defects cannot
affect the critical probability. However, proving this turns out to be rather taxing.

Let us remark that this difficultly with defects is not unique to the present paper. Benjamini and Schramm~\cite{BS}
prove a certain `conformal invariance' property of random
Voronoi percolation. This is {\em not} conformal invariance in the
sense of the celebrated
conjecture of Aizenman, Langlands, Pouliot and Saint-Aubin~\cite{Langlands_confinvar},
which is believed to hold in a very wide variety of contexts, and has been
proved for site percolation in the triangular lattice by Smirnov~\cite{Smirnov}.
Rather, Benjamini and Schramm prove (essentially) the following statement specific
to Voronoi percolation: fix a region 
$R\subset\RR^2$, and two segments $S_1$ and $S_2$ of its boundary.
Consider the Voronoi percolation associated with a Poisson process of intensity $\lambda$ on $R$,
using a certain metric ${\rm d}s$ to form the Voronoi cells, rather than the usual
Euclidean metric. Then, as $\lambda\to \infty$, a fixed conformal change in the metric ${\rm d}s$ does not affect
the probability that there is a black path from $S_1$ to $S_2$ by more than $o(1)$. Benjamini and Schramm
also prove a corresponding statement in 3 dimensions.
This statement is extremely unsurprising, at least in $\RR^2$:
because the change in the metric is conformal and, as $\lambda\to\infty$,
the density of cells is very high, the graphs associated to the Voronoi tilings for the two metrics coincide
almost everywhere. In other words, there are only a few `defects', where changing the metric causes
different cells to meet. In fact, in $2$ dimensions, Benjamini and Schramm note that there
are (in expectation) only a {\em bounded} number of defects. In $3$ dimensions
there are more, but still very few. Despite this fact, the result of Benjamini and Schramm is not at all easy:
even dealing with these very few defects requires a lot of work.

The discretization 
problems described above would not arise in the related setting of {\em random discrete Voronoi percolation},
in which, instead of considering
a Poisson process $Z$ on $\RR^2$, we take a random subset $L_\pi$ of a lattice $L$ (for example, $\Z^2$), 
where $L_\pi$ is formed by selecting points of $L$ independently with a certain probability $\pi$.
Letting $\pi\to 0$ and rescaling suitably, $L_\pi$
of course converges to $Z$ in a natural sense. If we define the Voronoi tiling with respect to $L_\pi$, $\pi>0$,
and then colour the cells black with probability $p$ independently of one another, we can show that percolation does
occur for any $p>1/2$ and $\pi>0$, with the following proviso. In this setting, more than three Voronoi cells may meet
at a vertex (a decreasing proportion of cells do so as $\pi\to 0$),
and we must consider all cells meeting at a vertex to be connected. In this setting, while 
we still need an equivalent of Theorem~\ref{th_RSW} (which follows from the proof
in Section~\ref{sec_RSW}), most of the remaining complications in this paper can
be avoided. We shall return to this in future work~\cite{ourKesten2}.

\section{External ingredients}\label{sec_ingred}

The proofs presented here will be mostly self-contained; we shall make use of two 
results from probabilistic combinatorics and two observations concerning percolation.

Let $X$ be a fixed ground set with $n$ elements, let $\Pow(X)$ denote
its power-set, and let ${\bf p}=(p_x)_{x\in X}$ be a vector
of probabilities. Let $X_{\bf p}$ be a random subset
of $X$ obtained by selecting each $x\in X$ independently with probability $p_x$.
For $\cA\subset \Pow(X)$, let $\Pr_{\bf p}(\cA)$ be the probability that $X_{\bf p}\in\cA$.
As usual, we say that $\cA$ is {\em increasing} if $A\in \cA$ and $A\subset B\subset X$
imply $B\in \cA$.

The first result we shall need is Harris' lemma, from his 1960 paper~\cite{Harris} in which he proved that $p_H\ge 1/2$
for bond percolation in $\Z^2$.
\begin{lemma}\label{l_harris}
 If $\cA$ and $\cB$ are increasing, then $\Pr_{\bf p}(\cA\cap \cB) \ge \Pr_{\bf p}(\cA)\Pr_{\bf p}(\cB)$.
\end{lemma}
In other words, increasing events in the weighted cube (the product of two-element probability spaces)
are positively correlated.
Harris' original statement is for the case $\bf p$ constant, but the extension to general ${\bf p}$ is
essentially equivalent. In fact, once one thinks of the statement, the proof of either form is
extremely simple using induction on $n$.
In~\cite{Harris}, the set $X$ was a finite set of edges of $\Z^2$, and $X_{\bf p}$ was the subset
of $X$ consisting of the open edges.
The extension to infinite $X$ is trivial.
Harris' Lemma, which was rediscovered by Kleitman~\cite{Kleitman} in a different context,
led to a series of generalizations culminating in the 
`Four-functions Theorem' of Ahlswede and Daykin~\cite{AD-FFT}.
However, in the context of percolation, it is often
exactly Harris' original lemma that is needed.

The second result we shall need is a modified form of
a sharp-threshold result of Friedgut and Kalai~\cite{FK},
which is itself a consequence of a result of Kahn, Kalai and Linial~\cite{KKL}
(see also~\cite{BKKKL}) concerning the influences of coordinates in a product space.
Let $W_{p_-,p_+}$ be the (weighted) three-element probability space $\{-1,0,1\}$
where the elements have respective probabilities $p_-$, $1-p_--p_+$ and $p_+$. We shall work
in the $n$th power $W_{p_-,p_+}^n$ of this space, writing $\Pr_{p_-,p_+}^n$ for the
corresponding probability measure on $\{-1,0,1\}^n$.

An event $E\subset\{-1,0,1\}^n$
is {\em increasing} if, whenever $x\in E$ and $x\le x'$ holds pointwise, we have $x'\in E$.
This is a natural generalization of the notion of increasing for subsets of $\Pow(X)$.
The event $E$ is {\em symmetric} if there is a permutation group acting transitively
on $[n]=\{1,2,\ldots,n\}$ whose induced action on $\{-1,0,1\}^n$ preserves $E$.

\begin{theorem}\label{th_sharpMOD}
There is an absolute constant $c_3$ such that if $0<q_-<p_-<1/e$, $0<p_+<q_+<1/e$,
$E\subset \{-1,0,1\}^n$
is symmetric and increasing, and
$\Pr_{p_-,p_+}^n(E)>\eta$, then $\Pr_{q_-,q_+}^n(E)>1-\eta$ whenever
\begin{equation}\label{e_sharpMOD}
 \min\{q_+-p_+,p_--q_-\} \ge c_3\log(1/\eta)\pmax\log(1/\pmax)/\log n,
\end{equation}
where $\pmax=\max\{q_+,p_-\}$.
\end{theorem}

\begin{proof}
The case $p_-=q_-=0$ of this result, a result about the weighted discrete cube,
is exactly Theorem 3.2 of Friedgut and Kalai~\cite{FK}; the proof
in~\cite{FK} extends to Theorem~\ref{th_sharpMOD} {\em mutatis mutandis}.
Indeed,
the key bound $w(f)\le cp\log(1/p)$ in the proof of Theorem 3.1 in~\cite{FK} holds in our case also,
with a different constant. The corresponding formula for our three-element space is
\[
 w(f)\le cp_+\log(1/p_+) + cp_-\log(1/p_-),
\]
and the right hand side is at most $2c\pmax\log(1/\pmax)$.
\end{proof}

Finally, we shall need two observations concerning $k$-dependent percolation in $\Z^2$.
By a {\em bond (site) percolation measure on $\Z^2$} we shall mean a probability measure
on the space of assignments of a {\em state}, i.e., {\em open} or {\em closed}, to each
edge (vertex) of $\Z^2$, with the usual $\sigma$-field of measurable events.
For bond (site) percolation on $\Z^2$, the {\em open cluster $C_0$ containing the origin}
is the set of all vertices of $\Z^2$ that may be reached from the origin by an {\em open path},
i.e., a path all of whose edges (vertices) are open. In the case of site percolation,
$C_0=\emptyset$ if the state of the origin is closed.

A bond (site) percolation measure on a graph $G$ is {\em $k$-dependent} if
for every pair $S$, $T$ of sets of edges (vertices) of $G$ at graph
distance at least $k$, the states (being open or closed) of the edges (vertices) in $S$ are independent of
the states of the edges (vertices) in $T$. For bond percolation, when $k=1$ the separation
condition is exactly that no edge of $S$ shares a vertex with an edge of $T$.

Liggett, Schonmann and Stacey~\cite{LSS}
proved in a more general context that, for any $k$ and any $p_1<1$, there is a $p_2<1$ such that
any $k$-dependent probability measure in which each edge is
open with probability at least $p_2$ stochastically dominates the product measure in which edges are open
with probability $p_1$.
An immediate consequence of this result is the following.

\begin{lemma}\label{l_1dep}
There is a constant $p_0<1$ such that for any $1$-dependent bond percolation measure on $\Z^2$
satisfying the
additional condition that each edge is open with probability at least $p_0$, the probability
that the origin is in an infinite open cluster is positive.
\end{lemma}

The best value of $p_0$ in this lemma that is  currently known is due to Balister, Bollob\'as and Walters~\cite{BBW},
who showed that one may take $p_0=0.8639$.
For us, the value of $p_0$ is not important.
As stated, Lemma~\ref{l_1dep} is essentially trivial from first principles.

The second observation is a variant of a converse of Lemma~\ref{l_1dep}, giving exponential decay
rather than percolation.
This time, it is easier to work with site percolation.
Recall that in this context the open cluster $C_0$ of the origin
is the set of vertices of $\Z^2$ joined to the origin by a path in $\Z^2$ every one
of whose vertices is open. Although the lemma below is also essentially trivial, we give a proof.

\begin{lemma}\label{l_kdepneg}
Let $k$ be a fixed positive integer, and let $\Prtwo$ be a $k$-dependent site percolation measure on $\Z^2$
in which
every vertex $v\in \Z^2$ is open with probability at most $p$. There is a constant $p_1=p_1(k)>0$ such that
for every $p\le p_1$ there is a $c(p,k)>0$ for which
\[
 \Prtwo( |C_0| \ge n)\le \exp(-c(p,k)n)
\]
for all $n\ge 1$.
\end{lemma}
\begin{proof}
If $|C_0|\ge n$, then the subgraph of $\Z^2$ induced by the open vertices
contains a tree $T$ with $n$ vertices, one of which is the origin.
It is well known and easy to check that the number of such trees in $\Z^2$ grows exponentially, and is at most $(4e)^n$.
Fix any such tree $T$. Then there is a subset $S$ of at least $n/(2k^2-2k+1)$ vertices of $T$ such that any $a,b\in S$
are at distance at least $k$; indeed, one can find such a set by a greedy algorithm:
whenever a vertex $a$ is chosen, the number of other vertices it rules out is at most the number
of other vertices of $\Z^2$ within graph distance $k-1$ of $a$, namely $4\binom{k}{2}=2k^2-2k$.
The vertices of $S$ are open independently, so the probability that every vertex of $T$ is open is at
most $p^{|S|}$.
Hence,
\[
  \Prtwo( |C_0| \ge n)\le (4e)^n p^{n/(2k^2-2k+1)}.
\]
Provided $p$ is small enough that $r=4e p^{1/(2k^2-2k+1)}<1$, the conclusion follows,
taking $c(p,k)=-\log r$.
\end{proof}

\section{Basic preliminaries}\label{sec_prelim}

Before we get down to the real work,
let us eliminate some degenerate cases.
Recall that $Z$ is a Poisson process of intensity $1$ on $\RR^2$.
It is easy to check that with probability $1$ every
Voronoi cell $V_Z(z)$ is bounded.
The probability that $Z$ contains $4$ points lying on a circle is zero. Thus, with probability $1$,
no $x\in \RR^2$ lies in more than three Voronoi cells. Finally, any bounded region contains finitely
many points of $Z$. It follows that, with probability $1$, every Voronoi cell $V_Z(z)$ is a closed
convex $k$-gon for some $k$, and that when two Voronoi cells meet, they share an edge.
We shall assume throughout that these non-degeneracy conditions hold {\em always}.

Here is a simple observation about Voronoi tilings, which we state as a lemma for ease of reference.
\begin{lemma}\label{l_closest}
Suppose that there is a point $x\in \RR^2$ such that the closest two points of $Z$ to $x$ are $z_1$ and $z_2$.
Then the Voronoi cells of $z_1$ and $z_2$ meet.
\end{lemma}

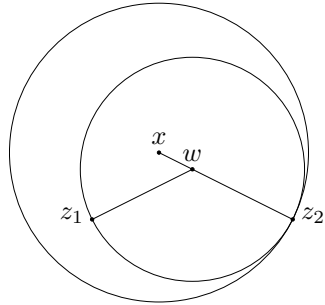
\begin{figure}[htb]
\centering
\input{2closest.pstex_t}
\caption{A point $x\in \RR^2$ and the two closest points $z_1$, $z_2$ of $Z$.
There are no points of $Z$ inside the inner circle shown.}\label{fig_2closest}
\end{figure}

\begin{proof}
Writing $\dist{u}{v}$ for the Euclidean distance between two points $u$, $v\in \RR^2$,
we may assume that $\dist{x}{z_1}\le\dist{x}{z_2}$.
But then there is a (unique) point $w$ on the segment $xz_2$ with
$\dist{w}{z_1}=\dist{w}{z_2}$, and for this point $w$ we have $\dist{w}{z_i}\ge \dist{w}{z_2}$
for every $z_i\in Z$; see Figure~\ref{fig_2closest}.
Hence $V(z_1)$ and $V(z_2)$ do meet -- they both contain $w$.
\end{proof}

Throughout this paper,
there will be a parameter $s$ controlling the {\em scale} of the region we are considering.
We say that an event $E=E(s)$ holds {\em with high probability}, or \whp, if it holds with probability
tending to $1$ as $s\to\infty$, with all other parameters (for example, $p$) fixed.
All $o(.)$ notation will refer to the same limit, so an event
holds \whp\ if and only if it holds with probability $1-o(1)$.
Most of our results concern the case $s\to\infty$. In the statement and proof
of such results we may assume that $s$ is larger than some fixed constant; often
we shall do so without comment.
Note that in a region of area $s^2$, we expect the largest Voronoi cell to have diameter $\Theta(\sqrt{\log s})$.

For a rectangle $R\subset \RR^2$ and a real number $r>0$, let us define the {\em $r$-neighbourhood}
$R[r]$ of $R$ in the usual way, i.e., as the set of all points within distance $r$ of some point of $R$.
Given $\rho>1$ and a $\rho s$ by $s$ rectangle $R_s$, $s>1$,
let $\Es(R_s)$ be the event that for every $x\in R_s[r]$ there is 
some point $z\in Z$ at distance $\dist{x}{z}<r$, 
where $r=2\sqrt{\log s}$.
We collect some simple properties of the event $\Es(R_s)$
in the following easy lemma. 

\begin{lemma}\label{l_small}
Let $\rho\ge 1$ be constant.
Let $Z$ be a Poisson process of intensity $1$ on $\RR^2$, let $R_s\subset \RR^2$ be a $\rho s$ by $s$ rectangle,
and set $r=2\sqrt{\log s}$.
Then $\Es(R_s)$ holds \whp. Also, $\Es(R_s)$ depends only on the restriction of $Z$ to $R_s[2r]$,
and if $\Es(R_s)$ holds, then the colour of every point of $R_s$ is determined by the restriction
of $(Z,\col)$ to $R_s[2r]$.
\end{lemma}

\begin{proof}
We may cover $R_s[r]$ with $O(s^2/\log s)=o(s^2)$ (half-open)
squares $S_i$ of side length $r/\sqrt{2}$ so that the $S_i$
are disjoint, and their union lies in $R_s[2r]$.
Each $S_i$ has area $r^2/2=2\log s$, so the number of points of $Z$ in $S_i$ has a Poisson distribution
with mean $2\log s$. Hence, the probability that a particular $S_i$ contains no points of $Z$ is exactly
$s^{-2}$, and \whp\ every $S_i$ contains at least one point of $Z$.
It follows that $\Es(R_s)$ holds \whp: given any $x\in R_s[r]$, any point $z\in Z$ in the same $S_i$ as $x$
satisfies $\dist{x}{z}<r$.

The remaining two claims are immediate from the definition of $\Es(R_s)$.
\end{proof}

\subsection{Increasing events are positively correlated}

Recall that we are working with $\Pr_p=\Pr_p^{\RR^2}$, the probability measure associated to
the coloured Poisson process $(Z,\col)$ where each point is coloured black independently with probability
$p$. We may think of this process as assigning a state from $\{-1,0,1\}$ to every point of $\RR^2$:
the black points of $Z$ have state $1$, the white points of $Z$ have state $-1$,
and every point of $\RR^2\setminus Z$ has
state $0$.
An event $E$ is {\em black-increasing}, or simply {\em increasing},
if $f$ is increasing in terms of the states of the points of $\RR^2$.
In other words, $E$ is preserved under the addition of black points to $Z$ and under
the deletion of white points from $Z$. The following result follows easily from Harris' Lemma.

\begin{lemma}\label{l_cor}
Let $E_1$ and $E_2$ be increasing events. Then for any $0\le p\le 1$ we have
\[
 \Pr_p(E_1\cap E_2) \ge \Pr_p(E_1)\Pr_p(E_2).
\]
\end{lemma}

We shall apply Lemma \ref{l_cor} to events such as `there is a black path across a certain
rectangle', to be defined below.

\subsection{Crossing a rectangle: two definitions}

Given a rectangle $R=[x_1,x_2]\times [y_1,y_2]$, $x_1<x_2$, $y_1<y_2$,
let $H'(R)=H'_b(R)$ be the event that there is a path $P'=C_1C_2\ldots C_n$ of black cells such that
$C_i\cap C_{i+1}\cap R$ is non-empty for $i=1,2,\ldots,n-1$, with
$C_1$ meeting the left-hand side of $R$ and $C_n$ meeting the right-hand side.
We call such a path a {\em path of black cells across $R$}. In the notation, the $H$
stands for horizontal, and the $b$ for black.

Recall that a point $x$ of $\RR^2$ is black if it lies in a black cell,
i.e., if (at least one of) the point(s) $z\in Z$ closest to $x$ is black.
We say that a subset $A\subset \RR^2$ is {\em black} if every point of $A$ is black.

Let $H(R)=H_b(R)$ be the event that there is a piecewise-linear black path $P$
inside $R$ from some point on the left-hand side to some point on the right-hand side.
We call such a path $P$ a {\em black path across $R$},
and say that $R$ may be {\em crossed (horizontally)}
by a black path if $H(R)$ holds.
Note that $H(R)$ is black-increasing:
the blackness of any point of $\RR^2$, or of any set, is black-increasing.

It is easy to see that $H(R)$ and $H'(R)$ coincide (at least up to probability
zero degenerate cases, which we have ruled out).
While the definition given for $H'(R)$ is perhaps easier to visualize, that given for $H(R)$
is easier to work with, and this is what we shall use most of the time.
We shall write $V(R)=V_b(R)$ for the event that there is a {\em black path crossing $R$ vertically},
defined analogously. Also, we shall write $H_w(R)$ and $V_w(R)$ for the corresponding
events with black replaced by white.

The basic starting point for the study of critical points for random
Voronoi percolation in $\RR^2$ is the following well-known fact.

\begin{lemma}\label{l_half}
If $p=1/2$ and $S$ is any square $[x,x+a]\times [y,y+a]$, $a>0$, then the probability of the event $H(S)$
that there is a black path across $S$ is {\em exactly} $1/2$.
\end{lemma}

\begin{figure}[htb]
\[\epsfig{file=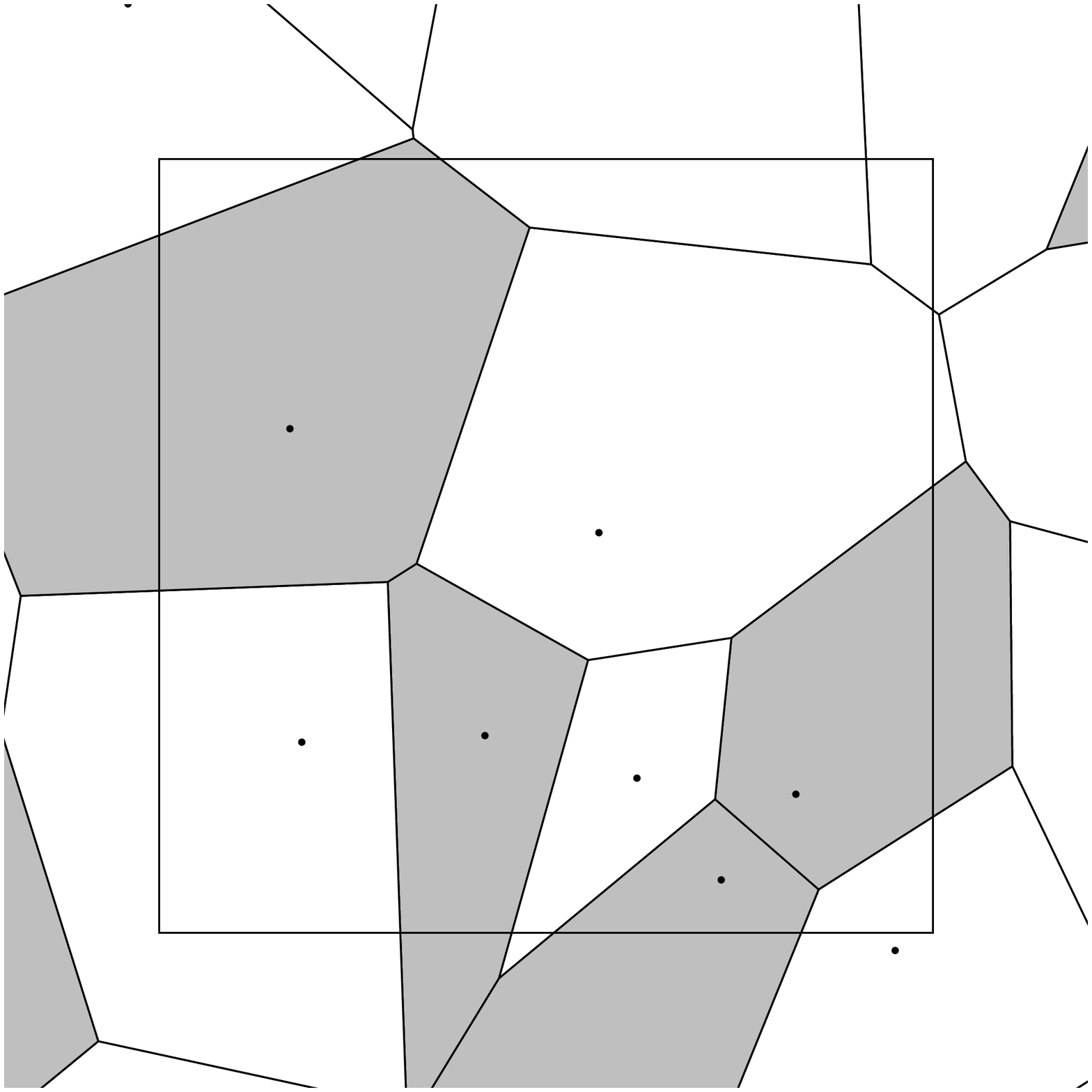,height=1.5in}\]
\caption{An example in which $H_b(S)$ does not hold and $V_w(S)$ does.}\label{fig_badcr}
\end{figure}

This result is the exact analogue of the well-known result for $\Z^2$ based on duality.
It is an immediate consequence of the fact that, with probability $1$, exactly one of the
events $H(S)=H_b(S)$ and $V_w(S)$ holds; see Figure~\ref{fig_badcr}. By symmetry, these events
have the same probability when $p=1/2$.
Random Voronoi percolation differs from bond percolation in $\Z^2$ in the following important respect:
events such as $H(R)$ for different rectangles $R$ are {\em not} independent,
as $H(R)$ depends on the points of $Z$ outside but near $R$, as well as those inside.

\section{Crossing a rectangle at $p=1/2$}\label{sec_RSW}

One of our main intermediate results is Theorem~\ref{th_RSW} below,
which will play a role
corresponding to the Russo-Seymour-Welsh (RSW) Theorem for bond percolation in $\Z^2$.
The RSW Theorem, proved independently by Russo~\cite{Russo} and Seymour and Welsh~\cite{SW} (see also~\cite{Grimmett}),
states that,
if the probability of a horizontal open crossing of an $s$ by $s$ square is at least $c>0$,
then for $\rho>1$ the probability
of a horizontal (i.e., long) open crossing of a $\rho s$ by $s$ rectangle is at least $f(\rho,c)>0$. The key point
is that $f(\rho,c)$ (which is, in fact, polynomial in $c$ when $\rho$ is fixed)
does not depend on $s$. Together with an analogue
of Lemma~\ref{l_half}, the RSW Theorem shows that in $p=1/2$ bond percolation on $\Z^2$ the probability of a horizontal open crossing
of a $10s$ by $s$ rectangle, say, is bounded away from zero as $s\to\infty$. It is very easy to deduce that percolation
does not occur at $p=1/2$.

The result we shall prove here does not directly correspond to the RSW Theorem: starting from the fact that
a square of any given size can be crossed by a black path with probability bounded away from zero (namely,
probability $1/2$ if $p=1/2$), we deduce
that {\em certain} (not necessarily {\em all})
very large $\rho s$ by $s$ rectangles can be crossed with probability bounded away from zero.
The proof is much longer than that of the RSW Theorem: none of the existing proofs of the latter
seems to transfer to this context.
This observation was made by Zvavitch~\cite{Zv}, who pointed out that the simplest
proof of the RSW Theorem, based on finding a `lowest' black left-right crossing of a square,
and noting that this is independent
of edges above the crossing, breaks down for the Voronoi tiling: the fact that various Voronoi cells form a path
depends on the {\em absence} of points of $Z$ even above the path.
Zvavitch notes that the proof of Alexander's result~\cite{Alex_RSW}
for a different Poisson percolation model does not
carry over to the Voronoi setting either.
In proving
that $p_H\ge 1/2$, Zvavitch uses a totally different method, not going via the RSW Theorem, but rather
adapting a proof of Harris' Theorem due to Burton and Keane~\cite{BurtonKeane} and Zhang
(see~\cite[pages 198 and 289]{Grimmett}),
which does carry over to the Voronoi case in a simple way.

In any case, the question of whether a direct analogue of the RSW Theorem holds for random Voronoi percolation is open.
The following weaker result will suffice for our purposes.
In the results below and for the rest of the paper, we shall make use of the function
\[
 f_p(\rho,s) = \Pr_p\big(H([0,\rho s]\times [0,s])\big).
\]
In other words, $f_p(\rho,s)$ is the $\Pr_p$-probability that
a $\rho s$ by $s$ rectangle has a horizontal black crossing.

\begin{theorem}\label{th_RSW}
Let $0<p<1$ and $\rho>1$ be fixed. If $\liminf_{s\to\infty} f_p(1,s)>0$, then $\limsup_{s\to\infty} f_p(\rho,s)>0$.
\end{theorem}

Before we prove this result, let us note that it has the following immediate consequence.

\begin{corollary}\label{c_RSWcon}
Let $\rho>1$ be fixed.
There is a constant $c_0=c_0(\rho)>0$ such that for every $s_0$ there is an $s>s_0$ with
$\Pr_{1/2}(H(R_s))\ge c_0$, where $R_s$ is a $\rho s$ by $s$ rectangle.
\hfill\endproof
\end{corollary}

\begin{proof}
Lemma~\ref{l_half} states exactly that $f_{1/2}(1,s)=1/2$ for all $s>0$. From Theorem~\ref{th_RSW}
it follows that $\limsup f_{1/2}(\rho,s)>0$, which is precisely the conclusion of the corollary.
\end{proof}

In principle, we could give an explicit value for $c_0(\rho)$, but it would be rather small for the
values of $\rho$ we shall consider, for example, $\rho=3$.
It is likely that $f_{1/2}(\rho,s)$ tends to a limit as $s\to\infty$ with $\rho$ fixed,
but as far as we are aware this has not
been proved. Indeed, this statement would be a weak form of a very special case of the conformal invariance
conjecture of~\cite{Langlands_confinvar}
mentioned in the introduction. However, even for bond percolation in
$\Z^2$, which has been studied much more extensively than random Voronoi percolation, the corresponding
question is still open.

\begin{proof}[Proof of Theorem~\ref{th_RSW}]
Assume, for a contradiction, that Theorem~\ref{th_RSW} does not hold, and
fix a particular value of $p$ for which the result fails. Then
there is a constant $c_1>0$ such that
\begin{equation}\label{f1}
 f_p(1,s)\ge c_1
\end{equation}
for all large enough $s$, but, for some fixed $\rho>1$,
\begin{equation}\label{fpto0}
 f_p(\rho,s)\to 0.
\end{equation}
Here, as usual, the limit is taken as $s\to\infty$
with all other parameters fixed. We shall make several claims
during the proof; all of these will be conditional on our assumptions,
which we shall show to be inconsistent.

Note that $f_p(\rho,s)$ is decreasing in $\rho$, since the corresponding events for different $\rho$
are nested. We claim that
for any $a,b\ge 1$ we have 
\begin{equation}\label{e3}
 f_p(a+b-1,s)\ge f_p(a,s)f_p(b,s)f_p(1,s).
\end{equation}
To see this, consider the events $E_1$, $E_2$ that there are
(piecewise-linear) black paths $P_1$, $P_2$ across $[0,as]\times [0,s]$ and $[(a-1)s,(a+b-1)s]\times [0,s]$,
respectively, and the event $E_3$ that there is a piecewise-linear black path $P_3$
\begin{figure}[htb]
 \[\epsfig{file=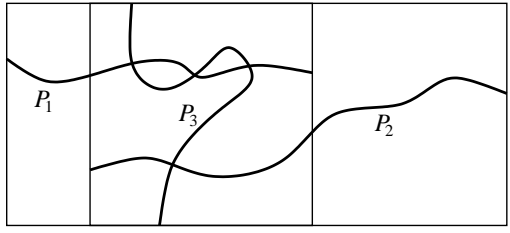,height=1.2in}\]
 \caption{Three paths whose union contains a horizontal crossing of the outer
rectangle.}
 \label{fig_p123}
\end{figure}
crossing the square $S=[(a-1)s,as]\times [0,s]$ from top to bottom; see
Figure~\ref{fig_p123}. The $E_i$ are increasing events
so, by Lemma~\ref{l_cor},
\[
 \Pr_p(E_1\cap E_2\cap E_3) \ge \Pr_p(E_1)\Pr_p(E_2)\Pr_p(E_3) = f_p(a,s)f_p(b,s)\Pr_p(E_3).
\]
But, of course, $\Pr_p(E_3)=\Pr_p(V(S))=\Pr_p(H(S))=f_p(1,s)$, using rotational
symmetry of the model.
Finally, if the $E_i$ all hold then $P_3$ meets both $P_1$ and $P_2$, so the union of the $P_i$ contains
a path across $[0,(a+b-1)s]\times [0,s]$ which is black. This proves \eqref{e3}.

Combining this with our assumptions \eqref{f1} and \eqref{fpto0},
for any $\epsilon>0$ we have 
\begin{equation}\label{cto}
 f_p(1+\epsilon,s)\to 0
\end{equation}
as $s\to\infty$; otherwise, as $f_p(1,s)$ is bounded away from zero by \eqref{f1}, using \eqref{e3}
$\lceil (\rho-1)/\epsilon\rceil$ times we have $f_p(\rho',s)\not\to 0$ for some $\rho'\ge \rho$, contradicting
\eqref{fpto0}.

Using \eqref{cto} we shall be able to deduce increasingly implausible properties of black paths
crossing certain rectangles. Throughout the argument below we work within the strip $T_s=[0,s]\times \RR$.
\setcounter{claim_base}{12}
\setcounter{claim}{0}
The first step is to show that if we try to cross the infinite strip $T_s$, then we almost
always stay nearly within height $s/2$ of the point we start from.

\begin{claim}\label{c_nottoobigNEW}
Let $\epsilon>0$ be fixed, and let $L$ be the line-segment $\{0\}\times [-\epsilon s,\epsilon s]$.
Assuming that \eqref{f1} and \eqref{fpto0} hold,
the probability that there is a black path $P$ in $T_s$ starting from $L$ and going outside
$S'=[0,s]\times [-(1/2+2\epsilon) s, (1/2+2\epsilon)s]$ tends to zero as $s\to\infty$.
\end{claim}
\begin{proof}
Considering the point at which the path $P$ first leaves $S'$,
by symmetry it suffices to show that the event $E$ that there is a black path $P_1$ lying entirely
within $S'$ and connecting some point of the segment $L$ to some point on the top side of $S'$
has probability tending to zero.

Let $E_1$ be the event that there is such a path $P_1$ lying within 
the rectangle $R=[0,s]\times [-s/2,s/2+2\epsilon s]$.
Note that if $E$ holds
and $E_1$ does not, then there is some black path within $T_s$ joining a point at height $-s/2$ to a point
at height $s/2+2\epsilon s$.
The shortest such path is a black path crossing a rectangle of width $s$ and height $(1+2\epsilon)s$
from top to bottom. Since such a path
exists with probability $f_p(1+2\epsilon,s)$, which tends to zero by \eqref{cto}, we have
$\Pr_p(E\setminus E_1)\le f_p(1+2\epsilon,s) \to 0$. Consequently, it suffices to show that $\Pr_p(E_1)\to 0$.

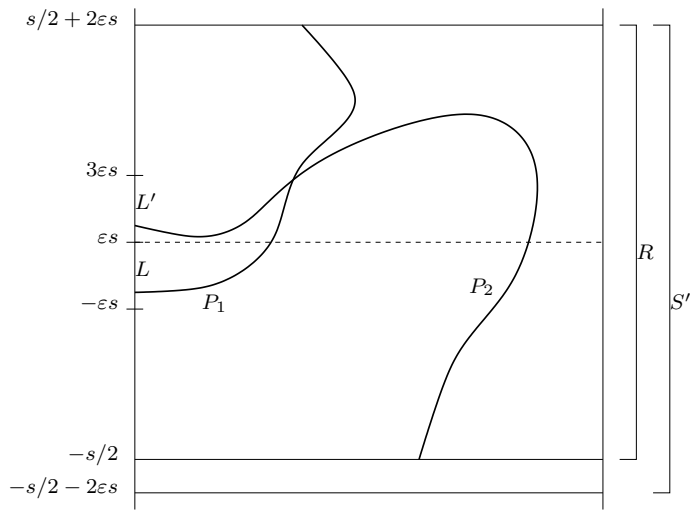
\begin{figure}[htb]
\centering
\input{claimfig1.pstex_t}
\caption{Paths $P_1$ and $P_2$ guaranteeing the events $E_1$ and $E_2$}\label{fig_nottoobig}
\end{figure}

Reflecting vertically in the line $y=\epsilon s$, a symmetry axis of $R$,
let $L'=\{0\}\times [\epsilon s,3\epsilon s]$
be the image of the segment $L$, and let
$E_2$ be the event that there is a black path $P_2$ within $R$
from $L'$ to some point with height $-s/2$; see Figure~\ref{fig_nottoobig}. Then
$\Pr_p(E_2)=\Pr_p(E_1)$ by symmetry.  Now $E_1$ and $E_2$ are increasing
events, so both hold with probability at least $\Pr_p(E_1)^2$.  But if
both hold then the corresponding paths $P_1$ and $P_2$ cross, as $L$ lies entirely below $L'$. Thus $P_1\cup P_2$
contains a black path crossing $R$ from top to
bottom. The probability that such a path exists is
$f_p(1+2\epsilon,s)$, which tends to zero by \eqref{cto}. Consequently
$\Pr_p(E_1)\to 0$ and hence $\Pr_p(E)\to 0$, proving the claim.
\end{proof}

Given a path $P\subset \RR^2$, let $\ymax(P)$ and $\ymin(P)$ be the supremum and infimum of the $y$-coordinates
of points on $P$. In other words, (as $P$ is closed) $\ymax$ and $\ymin$ are the maximum
and minimum `heights' attained by $P$.
Recall that an event holds {\em with high probability}, or \whp, if it holds with probability $1-o(1)$
as $s\to\infty$ with any other parameters (e.g., $p$, $\epsilon$) fixed. Recall also that when we say that a path $P$
crosses a rectangle $R$ horizontally, we mean that $P$ lies entirely within $R$, starts
on the left-hand side of $R$, and ends on the right-hand side of $R$.
We shall write $y_0(P)$ and $y_1(P)$ for the $y$-coordinates of the left and right endpoints of such a path $P$.

\begin{claim}\label{ymm}
Let $\epsilon>0$ and $C>0$ be fixed, and let $R_s$ be an $s$ by $2Cs$ rectangle.
Assuming that \eqref{f1} and \eqref{fpto0} hold, \whp\ every black path $P$ crossing $R_s$ horizontally satisfies
\begin{equation}\label{4heights}
 \big| |\ymaxmin(P)-y_i(P)| - s/2 \big| \le \epsilon s,
\end{equation}
for $\ymaxmin(P)=\ymax(P),\ymin(P)$ and $i=0,1$. In particular,
\begin{equation}\label{straight}
 |y_0(P)-y_1(P)|\le 2\epsilon s.
\end{equation}
\end{claim}
In other words, the maximum and minimum heights attained by $P$ are almost exactly $s/2$ larger and
smaller respectively than its starting and ending heights. Note that the claim becomes stronger
as $C$ increases; for $C<1/2$ it is trivial, as \eqref{cto} implies that in this case
\whp\ there is no black path $P$ crossing $R_s$ horizontally.

\begin{proof}
Without loss of generality we may take $R_s=[0,s]\times [-Cs,Cs]$.
It suffices to prove that \whp\ any black path $P$ crossing $R_s$ horizontally satisfies
\begin{equation}\label{ymax}
 y_0+s/2-\epsilon s \le \ymax \le y_0+s/2+\epsilon s.
\end{equation}
Here, and in what follows, we suppress the dependence of $\ymax(P)$ on $P$ when the path being
considered is clear from context, and use other self-explanatory abbreviations.
Inequality \eqref{ymax} gives one
of the four cases included in \eqref{4heights}. The others follow by reflecting horizontally
and vertically. The bound \eqref{straight} then follows by applying \eqref{4heights}
twice, with $\ymaxmin(P)=\ymax(P)$, say, and $y=0,1$.

Let us cover the left-hand side of $R_s$ by $\lceil 2C/(\epsilon/2)\rceil=O(1)$ line-segments
$L_i$ of length $\epsilon s/2$. Fixing $i$, let us call a black path $P$ {\em eligible} if $P$ 
crosses $R_s$ horizontally, starting from a point $(0,y_0)\in L_i$. As any path $P$ crossing $R_s$ horizontally
must contain a point $(0,y_0)\in L_j$ for some $j$,
it suffices to show that
\whp\ \eqref{ymax} holds for every eligible path $P$.

Let $(0,y)$ be the midpoint of $L_i$.
Applying Claim~\ref{c_nottoobigNEW} with $\epsilon/4$ in place of $\epsilon$, and using translational
invariance of the model, \whp\ every eligible path $P$ satisfies
\begin{equation}\label{ym1}
 \ymax\le y+s/2+\epsilon s/2
\end{equation}
and
\begin{equation}\label{ym2}
 \ymin\ge y-s/2-\epsilon s/2.
\end{equation}
As $|y_0-y|\le \epsilon s/4$ from the definition of eligibility, \eqref{ym1} implies the upper bound in \eqref{ymax}.
It remains to prove the lower bound. Suppose that it is not true that the lower bound
in \eqref{ymax} holds \whp. Then
with probability bounded away from $0$ there is an eligible path with
\begin{equation}\label{ym3}
 \ymax \le y_0+s/2-\epsilon s \le y+s/2-3\epsilon s/4.
\end{equation}
As \whp\ every eligible path satisfies \eqref{ym2},
we see that with probability bounded away from zero there is some black path $P$
crossing $R_s$ horizontally for which both \eqref{ym2} and \eqref{ym3} hold. As such a path crosses a fixed
rectangle of width $s$ and height $(1-\epsilon/4)s$ horizontally, this conclusion
contradicts \eqref{cto}, completing the proof of \eqref{ymax} and hence of the claim.
\end{proof}

The conclusions of Claims~\ref{c_nottoobigNEW} and~\ref{ymm} may seem fairly plausible. Recall, however, that we are aiming
for a contradiction. In the next two claims, we show that our assumptions imply properties
of paths crossing certain rectangles that are far from plausible.
 
\begin{claim}\label{cl_split2}
Let $C>0$ be fixed, and let $R=R_s$ be the $s$ by $2Cs$ rectangle $[0,s]\times [-Cs,Cs]$. For $i=0,1$, set
$R_i=[is/100,(i+99)s/100]\times [-Cs,Cs]$.
Assuming that \eqref{f1} and \eqref{fpto0} hold,
\whp\ every black path $P$ crossing $R$ horizontally contains disjoint black paths $P_0$ and $P_1$ such that
$P_i$ crosses $R_i$ horizontally.
\end{claim}

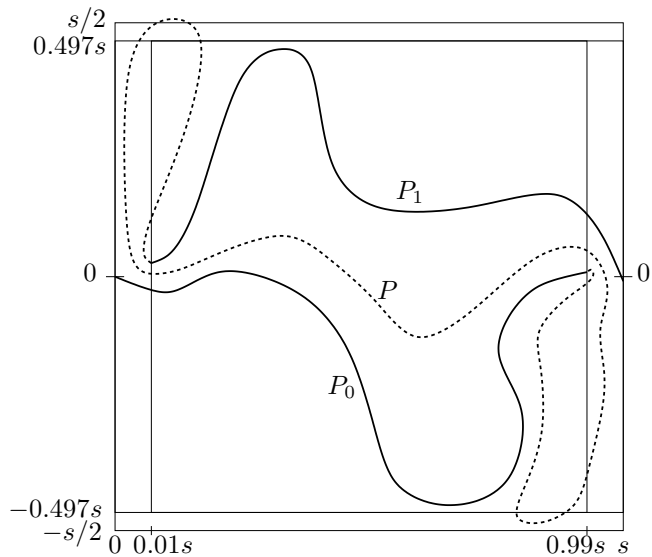
\begin{figure}[htb]
\centering
\input{newsplit2.pstex_t}
\caption{Two sub-paths $P_0$, $P_1$ (solid lines) of a path $P$ crossing horizontally from $x=0$
to $x=s$. The numbers outside the rectangle are $x$- or $y$-coordinates; to avoid clutter,
$y$-coordinates are shown relative to $y_0(P)$.}
\label{fig_split2}
\end{figure}

\begin{proof}
Let $E$ be the event that every black path $P$ crossing $R$ horizontally satisfies
\begin{equation}\label{4hp}
 \big| |\ymaxmin(P)-y_j(P)| - s/2 \big| \le s/1000
\end{equation}
for $j=0,1$.
Note that $\Pr_p(E)=1-o(1)$ by Claim~\ref{ymm}, applied with $\epsilon=1/1000$ to $R$.
Similarly, for $i=0,1$, let $E_i$ be the event that every black path $P_i$ crossing $R_i$
horizontally satisfies
\begin{equation}\label{4hp'}
 \big| |\ymaxmin(P_i)-y_j(P_i)| - 0.495s \big| \le s/1000
\end{equation}
for $j=0,1$. By Claim~\ref{ymm} applied to $R_i$ with $0.99s$ in place of $s$, $C/0.99$ in place of $C$,
and $\epsilon=1/1000$, each event $E_i$ holds \whp.

We shall assume, as we may, that $E\cap E_0\cap E_1$ holds.
Let $P$ be any black path crossing $R$ horizontally, starting at
$(0,y_0)$ and ending at $(s,y_1)$, say.
Consider the first time that $P$ crosses the line $x=0.99s$.
Let $P_0$ be the initial segment of $P$ from $(0,y_0)$ up to this point.
Similarly, let $P_1$ be the final segment of $P$, defined by going backwards along $P$ from
$(s,y_1)$ until the line $x=0.01s$ is first reached; see Figure~\ref{fig_split2}.
Then each $P_i$ crosses $R_i$ horizontally.
Applying \eqref{4hp} twice, as in the proof of \eqref{straight}, we have $|y_0-y_1|\le 2s/1000$.
Hence, from \eqref{4hp'},
each of $P_0$ and $P_1$ remains between the heights $y_0-0.497s$ and $y_0+0.497s$. It follows that
$P_0$ and $P_1$ are disjoint: if they meet, their union contains a path $P'$ crossing
$R$ from left to right with $\ymax(P')-\ymin(P')\le 0.994s$. Such a path $P'$
cannot satisfy \eqref{4hp}, contradicting our assumption that $E$ holds.
\end{proof}

Note that the properties of the path $P$ and its sub-paths $P_0$, $P_1$ given by
\eqref{4hp} and \eqref{4hp'} are even less plausible than Figure~\ref{fig_split2} might suggest:
$P_0$ and $P_1$ both
have maximum and minimum heights very close to $y_0\pm 0.495s$. Taking $P$ as a shortest
path, the rest of $P$
must somehow thread its way between $P_0$ and $P_1$ to join one to the other, without crossing either.

\begin{claim}\label{cl_split16}
Let $C>0$ be fixed, and let $R=R_s$ be the $s$ by $2Cs$ rectangle $[0,s]\times [-Cs,Cs]$. For $0\le j\le 4$, set
$R_j=[js/100,(j+96)s/100]\times [-Cs,Cs]$.
Assuming that \eqref{f1} and \eqref{fpto0} hold,
\whp\ every black path $P$ crossing $R$ horizontally contains $16$ disjoint black paths $P_i$, $1\le i\le 16$,
where each $P_i$ crosses some $R_j$ horizontally.
\end{claim}

\begin{proof}
Arguing as in the proof of Claim~\ref{cl_split2},
\whp\ any path $P_0$ obtained as in that proof has disjoint initial
and final segments $P_{00}$ and $P_{01}$ crossing rectangles of width $0.98s$, and so on.
Alternatively, we could apply a slightly modified Claim~\ref{cl_split2} (with $0.99$ replaced by $0.98/0.99$, etc.)
multiple times.
\end{proof}

At this point it seems that it should be easy to reach a contradiction
and hence complete the proof of Theorem~\ref{th_RSW}.
Note that all paths $Q$ we consider are piecewise linear. Hence the geometric
length $|Q|$ is well-defined. Let $R=R_s=[0,s]\times [-s,s]$.
We have shown that \whp\ any black path $P$ crossing $R$ horizontally
contains 16 disjoint paths $P_i$ crossing strips of width $0.96s$, with
$|P|\ge \sum_{i=1}^{16} |P_i|$.

Digressing from the formal proof for a moment,
suppose for simplicity that $p=1/2$ (which is the important case).
Under our assumption that the conclusion of Theorem~\ref{th_RSW}
does not hold, from \eqref{cto} we have
$f_{1/2}(1+\epsilon,s)\to 0$ for any $\epsilon>0$. Recalling that $V_w(R)$ is the event
that there is a vertical crossing of $R$ by a white path, it follows that
$\Pr_{1/2}(V_w(R))=o(1)$, so $\Pr_{1/2}(H(R))=1-o(1)$. Hence, \whp\ $R$ has a horizontal black crossing;
let $P$ be a shortest such crossing.
Considering the sub-paths $P_i$ of $P$ given by Claim~\ref{cl_split16}, one might hope that the expected length of $P$
is at least $16+o(1)$ times the expected length of the shortest crossing of a slightly smaller rectangle,
which would lead to a contradiction. Unfortunately, we cannot work with expected lengths,
as we have no control over the path lengths in the $o(1)$ probability exceptional cases.
In general, there is a problem with these cases, as their probabilities tend to accumulate. 
To deal with this, we shall subdivide the paths $P_i$ once more, this time
in such a way as to introduce independence.

Returning to the formal proof of Theorem~\ref{th_RSW},
given any rectangle $R$, let $\short(R)$ be the minimum length
of a (piecewise-linear) black path $P$ crossing $R$ horizontally.
We take $\short(R)=\infty$ if there is no such path. Of course, we could define
$\short(R)$ as an infimum, but it is easy to see that this infimum is attained.

Let $R_s$ be an $s$ by $2s$ rectangle, where $s>1$.
Setting $r=2\sqrt{\log s}$, let $\Es(R_s)$ be the event considered in Lemma~\ref{l_small}.
We define a modification $\shortm(R_s)$ of the random variable $\short(R_s)$ as follows:
if $\Es(R_s)$ does not hold, set $\shortm(R_s)=0$; otherwise, set $\shortm(R_s)=\short(R_s)$.
By Lemma~\ref{l_small}, $\Es(R_s)$ holds \whp, so \whp
\begin{equation}\label{same}
 \shortm(R_s)=\short(R_s).
\end{equation}
As any path across $R_s$ has length at least $s$, it follows that the inequality
\begin{equation}\label{start}
 \shortm(R_s)\ge s
\end{equation}
holds \whp.
As $R_s$ has width $s$ and height $2s>s$, with probability at least $f_p(1,s)$
there is some black path $P$ crossing $R$ horizontally. 
From \eqref{f1}
there is a constant $c_1>0$ such that $f_p(1,s)\ge c_1$ for all (large enough) $s$.
Using \eqref{same}, it follows that if $s$ is large enough, then
\begin{equation}\label{finite}
 \Pr_p\big(\shortm(R_s)<\infty\big)\ge c_1/2
\end{equation}
holds for any $s$ by $2s$ rectangle $R_s$.

Let $0<\etatwo<\min\{c_1/3,10^{-4}\}$ be fixed, and
define a function $\fun(s)$ by
\[
 \fun(s) = \sup\{x: \Pr_p\big(\shortm(R_s) < x\big) \le \etatwo\},
\]
where $R_s$ is any $s$ by $2s$ rectangle.
This makes sense as the distribution of $\shortm(R_s)$ depends on $s$ only and not on the location of the rectangle.
From \eqref{finite}, if $s$ is large enough, which we shall assume from now on, then
$\fun(s)<\infty$,
as $c_1/2 > \etatwo$.
It follows that the supremum in the definition of $\fun(s)$ is attained, so 
\begin{equation}\label{att}
 \Pr_p\big(\shortm(R_s) < \fun(s)\big)\le \etatwo.
\end{equation}

The key property of $\shortm(R_s)$ is described in the following claim.

\begin{claim}\label{cl_indep}
Let $R_1$ and $R_2$ be two $s$ by $2s$ rectangles separated by a distance of
at least $s/100$. If $s$ is large enough, then the random variables
$\shortm(R_1)$ and $\shortm(R_2)$ are independent.
\end{claim}

\begin{proof}
Let us write $r=r(s)$ for $2\sqrt{\log s}$ as before. We take $s$ large enough
that $4r(s)\le s/100$.
By Lemma~\ref{l_small}, for $R=R_1$ or $R_2$,
whether $\Es(R)$ holds depends only on the intersection of
the Poisson process $Z$ with the $2r$-neighbourhood $R[2r]$ of $R$.
Also, given that $\Es(R)$ does hold, the colour of every point of $R$ is determined by $Z\cap R[2r]$ and the 
colours of these points.
Hence $\shortm(R)$ depends only on $Z\cap R[2r]$ and the black/white-colouring of these points.
As $R_1[2r]$ and $R_2[2r]$ are disjoint, the claim follows.
\end{proof}

After all this work the proof of Theorem~\ref{th_RSW} is not far away. However, in order to complete
it, we need one more technical assertion.

\begin{claim}\label{cl_sqr}
Let $R_s$ be a fixed $0.96s$ by $2s$ rectangle.
Assuming that \eqref{f1} and \eqref{fpto0} hold,
if $s$ is large enough then
\[
 \Pr_p\big(\short(R_s) < \fun(0.47s)\big) \le 200\etatwo^2.
\]
In other words, with probability at least $1-200\etatwo^2$, every black path $P$ crossing $R$
horizontally has length $|P|\ge \fun(0.47s)$.
\end{claim}

\begin{proof}
We take $R_s=[0,0.96s]\times [-s,s]$.
Let us cover the left-hand side of $R_s$ by $100$ line-segments $L_i$ of length
$0.02s$. Fixing $i$, let us say that a black path $P$ crossing $R$ horizontally
is {\em eligible} if $P$ starts from a point of $L_i$.
Let $B$ ($=B_i$) be the event that there is an eligible path $P$ with $|P|<t(0.47s)$.
As every black $P$ crossing $R$ starts from a point of some $L_j$, it suffices to prove
that $\Pr_p(B)\le 2\etatwo^2$.

Let $(0,y)$ be the midpoint of $L_i$. Consider an eligible path $P$,
and let $(0,y_0)$, $(0.96s,y_1)$ be its endpoints.
Let $P_0=P_0(P)$ be the initial segment of $P$ from $(0,y_0)$, stopping the first time $P$
reaches the line $x=0.47s$, a little less than half way across $R_s$.
Similarly, let $P_1=P_1(P)$ be the final segment defined by going
backwards along $P$ from $(0.96s,y_1)$ to the line $x=0.49s$.
Let $R_0=[0,0.47s]\times [y-0.47s,y+0.47s]$ and
$R_1=[0.49s,0.96s]\times [y-0.47s,y+0.47s]$ be rectangles of width $0.47s$ and height twice this;
see Figure~\ref{fig_resplit}.

\begin{figure}[htb]
\centering
\input{resplit.pstex_t}
\caption{The upper portion of $R_s$, showing a path $P$ and its initial and final segments $P_i$ (drawn solid).}
\label{fig_resplit}
\end{figure}

From our previous claims, \whp\ every eligible path $P$ has the following properties.
Firstly, $|y_0-y|\le 0.01s$ (by definition of eligibility). Secondly, $|y_1-y_0|\le 0.01s$,
by Claim~\ref{ymm} applied to the rectangle $R_s$, with $0.96s$ in place of $s$.
Thirdly, $P_0$ is contained in $R_0$, by Claim~\ref{ymm} applied to the rectangle
$[0,0.47s]\times [-s,s]$ with $0.47s$ in place of $s$, using $|y_0-y|\le 0.01s$.
Fourthly, $P_1$ is contained in $R_1$, using $|y_1-y|\le 0.02s$ and Claim~\ref{ymm}.

It follows that \whp\ every eligible path $P$ has length at least $\short(R_0)+\short(R_1)$.
As the $R_i$ are $s'=0.47s$ by $2s'$ rectangles, using \eqref{same}
it follows that \whp\ every eligible
path satisfies $|P|\ge \shortm(R_0)+\shortm(R_1)$.
But $\shortm(R_0)$ and $\shortm(R_1)$ are independent by Claim~\ref{cl_indep}, applied
with $0.47s$ in place of $s$. Hence,
using \eqref{att},
\[
 \Pr_p\big( \shortm(R_0)+\shortm(R_1) < \fun(0.47s) \big) \le \etatwo^2.
\]
Thus, the probability that some eligible path $P$ has $|P|<\fun(0.47s)$ is at most
$\etatwo^2+o(1)\le 2\etatwo^2$.
\end{proof}

We can now complete the proof of Theorem~\ref{th_RSW}.
Indeed, the results above imply that, if $s$ is large enough, then
\begin{equation}\label{grow}
 \fun(s)\ge 16\fun(0.47s).
\end{equation}
To see this, consider the rectangle $R=R_s=[0,s]\times [-s,s]$.
Consider the $5$ sub-rectangles $R_j=[js/100,(j+96)s/100]\times [-s,s]$, $0\le j\le 4$,
of $R_s$. Note that each $R_j$ has width $0.96s$ and height $2s$.
By Claim~\ref{cl_sqr}, with probability at least $1-1000\etatwo^2$ we
have $\short(R_j)\ge \fun(0.47s)$ for each $j$.
But by Claim~\ref{cl_split16}, \whp\ every black path $P$ crossing $R$ horizontally has length
at least $16\min_j\short(R_j)$.
Hence, with probability at least $1-1001\etatwo^2$ we have $\short(R)\ge 16\fun(0.47s)$.
From \eqref{same} we have $\shortm(R)=\short(R)$ \whp, so it follows that if $s$ is large enough then
\[
 \Pr_p\big(\shortm(R)< 16\fun(0.47s)\big)\le 1002\etatwo^2 \le \etatwo.
\]
But then $\fun(s)\ge 16\fun(0.47s)$ by the definition of $\fun(s)$, proving \eqref{grow}.

Taking logarithms, \eqref{grow} states that when $s$ is large enough,
if $\log s$ increases by $\log(1/0.47)$,
then $\log \fun(s)$ increases by at least $\log(16) > 3\log(1/0.47)$. Using as
`initial condition' the fact that $\fun(s)\ge s$ for large $s$, implied by \eqref{start},
it follows that $\fun(s)$ grows at least as fast as a constant times $s^{\log 16/\log(1/0.47)}$,
and in particular that $\fun(s)\ge s^3$ for
$s$ large enough. This conclusion is absurd because, considering an $s$ by $2s$ rectangle
$R_s$, from \eqref{same}, \eqref{finite} and \eqref{att},
for large $s$ we have
\begin{equation}\label{nl1}
 \Pr_p\big( \fun(s) \le \short(R_s) < \infty\big) \ge c_1/2-\etatwo-o(1) \ge c_1/7>0.
\end{equation}
In other words, if $s$ is large enough, 
with probability bounded away from zero the shortest black path $P$ crossing $R_s$
horizontally exists and  has length at least $\fun(s)\ge s^3$.
Let $r=2\sqrt{\log s}$, as before, and let $N$ be the number of points
of $Z$ in $R_s[r]$. Then $N$ has a Poisson distribution with mean $\area(R_s[r])<(s+2r)(2s+2r)< 3s^2$,
so $N\le 4s^2$ holds \whp. Suppose that $\Es(R_s)$ holds, as it does \whp\ by Lemma \ref{l_small},
and that $N\le 4s^2$.
Then every point of $R_s$ is within distance $r$ of the centre of any Voronoi cell it
lies in, so there are at most $N\le 4s^2$ cells meeting $R_s$, and the diameter of the intersection
of each cell with $R_s$ is at most $2r$.
Any shortest black crossing $P$ meets each cell at most once, and then
in a line segment, and so under our assumptions has length at most $8rs^2<s^3$.
Hence, \whp\ either $\short(R_s)=\infty$ or $\short(R_s)<s^3$, contradicting \eqref{nl1} when $\fun(s)\ge s^3$.
This contradiction completes the proof of Theorem~\ref{th_RSW}.
\end{proof}

Let us note for possible future reference that, while our proof of Theorem~\ref{th_RSW} is rather indirect,
it is also rather general. The only properties of the measure $\Pr_p$ and notion of black crossing
needed for the proof are the following. Firstly, crossings can be defined by the `blackness' of geometric paths
(so, for example, horizontal and vertical crossings of the same rectangle meet, enabling the
combination of certain crossings to form longer crossings). Secondly, the existence of certain
crossings must be `increasing' events, so that positive correlation of such events holds. Thirdly,
the crossing probabilities must be invariant under certain symmetries of $\RR^2$. Note that invariance under
the symmetries of $\Z^2$ certainly suffices, as we need only consider rectangles with integer coordinates.
Fourthly, some kind of asymptotic independence is
needed: certainly independence of regions separated by more than a certain distance (which does
not hold for random Voronoi percolation) is more than enough.

From the remarks above, the proof of Theorem~\ref{th_RSW} applies, for example, to bond or
site percolation in the square lattice; of course, the existing proofs
of the stronger RSW Theorem are much simpler in these contexts. Our proof also applies
to the Poisson model considered by Alexander~\cite{Alex_RSW}, whose stronger result
for this simpler model has a rather long proof, and to the random discrete Voronoi setting
considered in the introduction.

In what follows, all we shall use from this section is Corollary~\ref{c_RSWcon}. The claims used in the proof,
which in any case assumed a false assumption, will not be needed later.

\section{The torus}\label{sec_torus}

{From} now on we shall work mostly in the torus $\TT(s)$, $s>0$,
obtained by identifying opposite sides of the square $[0,s]^2\subset \RR^2$. To reduce unnecessary
clutter, most of the time we shall
suppress the dependence on $s$.
Throughout this paper, we are interested only in the $s\to\infty$ limit: in the end, we shall show that
for any $s_0$, certain events hold with certain probabilities in some $\TT(s)$, $s\ge s_0$.
Thus we may assume that $s$ is larger than any given constant; often we shall do so without comment,
assuming inequalities such as $\log s\ge 2\sqrt{\log s}$ without comment.

We consider the (metric induced by) the Euclidean metric on $\TT=\TT(s)$. We may define a Poisson process $Z=Z^\TT$ on $\TT$
with intensity $1$, and the corresponding random Voronoi tiling, which we may colour as before:
each point $z$ of $Z^\TT$ is black with probability $p$, independently of all other points, and white
otherwise. Every point $x$ in the Voronoi cell associated to $z$ inherits the colour of $z$. (So some
boundary points are both black and white.) We write
$\Pr_p^\TT=\Pr_p^{\TT(s)}$ for the corresponding probability measure.

Of course, we may couple $Z^\TT$ and $Z^{\RR^2}$ so that they agree on $S=[0,s]^2\subset \RR^2$. We shall
see below that \whp\ every
disc of radius $\sqrt{\log s}$ contained in $S$ contains some point of $Z^\TT$. It follows that \whp\ the Voronoi
cells associated with $Z^\TT$ meeting $S'=[\log s,s-\log s]^2$ and the cells associated with $Z^{\RR^2}$
meeting $S'$ coincide. In particular, we shall be able to apply Corollary~\ref{c_RSWcon} to a rectangle
in $\TT$ as long as the longer side has length at most $s-\log s$.

We note two basic facts about the maximum and minimum number of points in discs of certain sizes in the torus.
From now on we write $Z$ for $Z^\TT$.

\begin{lemma}\label{l_disc}
With probability $1-o(1)$
every disc of radius $\sqrt{\log s}$ in $\TT(s)$ contains at least
one point of $Z$. Furthermore, for every constant $a>0$ there is a $\rho=\rho(a)>0$, and for every
constant $\rho>0$ there is an $a=a(\rho)>0$,  
such that
with probability $1-o(s^{-10})$ no disc of radius $\rho\sqrt{\log s}$ in $\TT(s)$ contains more than $a\log s$
points of $Z$.
\end{lemma}

We omit the standard proof based on the first-moment method.

\section{Approximation preliminaries}\label{sec_approx}

We are going to approximate the continuous Poisson process on $\TT$ using a fine grid. There is a
problem with the order of limits, as our application of Theorem~\ref{th_sharpMOD} will restrict how
fine the grid can be as a function of $s$, the scale of the torus. This gives rise to `defects' or
potential ambiguities, i.e., places where the discrete approximation of the location
of the points does not tell us which Voronoi cells actually meet.
If we could take a fine enough grid (in terms of $s$), we could say that with high probability
there are no defects in $\TT(s)$; unfortunately we cannot do this.
While the density of defects will be low, some will occur, and we must argue
that they do not make much difference. More precisely, we shall compare their
effect to the effect of changing $p$ slightly, and show that the latter dominates
in some precise sense.
This is not very surprising, but seems to be a little fiddly to prove rigorously.

An analogous difficulty was encountered by Benjamini and Schramm in~\cite{BS},
and also required considerable work to overcome. Essentially,
the result of~\cite{BS} is that at any $p$, a very small density of `defects' (a bounded number in two dimensions) has
negligible effect on the probability that a large region may be crossed.
(See the introduction for further details.) Here, the density of defects is much higher
(an arbitrarily small negative power of the area of the region). On the other hand, the fact that we can vary
$p$ gives us much more elbow room.

Let us say that a point $x\in \TT$ is {\em $\delta$-robustly black} if
the closest black point of $Z$ is at least a distance
$\delta$ closer to $x$ than the closest white point.
A piecewise-linear path $P$ in $\TT$ is {\em $\delta$-robustly black}
if every point $x\in P$ is.
When we come to make our discrete approximation, $\delta$-robustly black paths will be very useful;
such a path clearly remains black if we adjust the position of each point of $Z$ by a small amount
(less than $\delta/2$).
In the following result, for $p_1<p_2$ we couple the two probability measures (on black/white-coloured Poisson
processes) $\Pr_{p_1}^{\TT(s)}$ and $\Pr_{p_2}^{\TT(s)}$. In other words, we construct
simultaneously Poisson processes $Z_1$ and $Z_2$ of intensity $1$ on $\TT(s)$, with associated colourings
$\col_1$ and $\col_2$, so that given $Z_i$ each point of $Z_i$ is coloured black with probability $p_i$
independently of all other points of $Z_i$, and white otherwise.

\begin{theorem}\label{th_couplerobust}
Let $0<p_1<p_2<1$ and $\epsilon>0$ be given. Let $\delta=\delta(s)$ be any function
with $\delta(s)\le s^{-\epsilon}$.
We may construct in the same probability space
coloured Poisson processes $(Z_1,\col_1)$ and $(Z_2,\col_2)$ on $\TT(s)$ as described above,
in such a way that the following global event holds \whp\ as
$s\to\infty$: for every piecewise-linear path $P_1$ which is black with respect
to $(Z_1,\col_1)$ there is a piecewise-linear path $P_2$ which is $4\delta$-robustly black with respect
to $(Z_2,\col_2)$, such that every point of $P_2$ is within distance $(\log s)^2$ of some point
of $P_1$ and vice-versa.
\end{theorem}

The only aim of this section is to prove the result above. This will involve overcoming several bothersome technicalities.
Nothing used in the (rather long) proof of Theorem~\ref{th_couplerobust} will be used in other parts of the paper.

\begin{figure}[htb]
\centering
\input{baddef.pstex_t}
\caption{A $\delta$-bad quadruple $\{z_1,z_2,z_3,z_4\}\subset Z$. There are no points of $Z$ in the inner
circle.}\label{fig_baddef}
\end{figure}
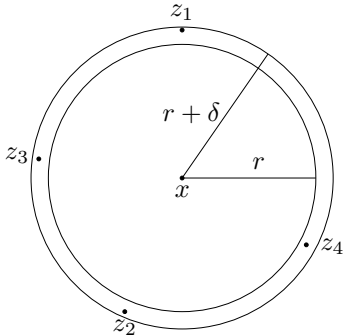

\subsection{Good points are useful}

Let us say that a quadruple of points $\{z_1,z_2,z_3,z_4\}$ of the Poisson process $Z$ on $\TT=\TT(s)$
is {\em $\delta$-bad} if there is a point $x$ of $\TT$ and a distance $r<2\sqrt{\log s}$ such that
(a) each $z_i$ is within distance $r+\delta$ of $x$, and (b) no point of $Z$ is within distance
$r$ of $x$; see Figure~\ref{fig_baddef}.
 In other words, the four points are close to a circle which is empty
of points of $Z$. Let us say that a pair of points $z_1$, $z_2$ of $Z$ is {\em $\delta$-close} if
they are within distance $\delta$.
Finally, a single point $z\in Z$ is {\em $\delta$-bad} if it is in a $\delta$-bad quadruple
or a $\delta$-close pair. Note that if $\delta<\delta'$, then every $\delta$-bad point
is also $\delta'$-bad. A point $z\in Z$ is {\em $\delta$-good} if it is not $\delta$-bad.
The reason for considering $\delta$-good points is that paths of cells associated to 
$\delta$-good points give rise to piecewise-linear paths that are robust. We
write $\dist{x}{y}$ for the Euclidean distance between two points $x$ and $y$ of $\TT$.

\begin{lemma}\label{l_rejoinbad}
Let $s\ge e^{100}$ and $\delta<1/\log s$ be given, and
suppose that no disc of radius $2\sqrt{\log s}$ in $\TT$ is free of points of $Z$.
Let $z_1$, $z_2\in Z$ be points whose Voronoi cells share an edge, and let $M$ be the midpoint of that edge.
If some $z\in Z\setminus\{z_1,z_2\}$ is at distance less than $\dist{z_1}{M}+\delta^6$ from $M$, then
either $\{z_1,z_2\}$ is a $\delta$-close pair, or
there is a $\delta$-bad quadruple $\{z_1,z_2,z,z'\}$, $z'\in Z$.
\end{lemma}

The proof of this deterministic statement is given in Section~\ref{a_l_rejoinbad} of the appendix.
Lemma~\ref{l_rejoinbad} has the following immediate corollary.

\begin{corollary}\label{c_goodsolid}
Suppose that no disc of radius $2\sqrt{\log s}$ in $\TT$ is free of points of $Z$.
Let $z_1,z_2$ be $\delta$-good black points of $Z$ whose Voronoi cells share an edge, and let $M$ be the
midpoint of that edge. Then the piecewise-linear path $z_1Mz_2$ is $\delta^6$-robustly black.
\end{corollary}
\begin{proof}
By assumption, $z_1$ and $z_2$ are not $\delta$-close, and
there is no $\delta$-bad quadruple containing $z_1$ and $z_2$, so 
by Lemma~\ref{l_rejoinbad} every point $z\ne z_1,z_2$
of $Z$ is at least $\delta^6$ further away from $M$ than $z_1$ is. Hence $M$ is $\delta^6$-robustly black,
and so is the line-segment $z_1M$: as we move towards $z_1$ along this line segment at rate 1,
the distance from $z_1$ decreases at rate $1$, while the distances from all other points decrease
at most this fast. Similarly $z_2M$ is $\delta^6$-robustly black.
\end{proof}

We have shown that a `good' black path gives rise to a robustly-black
path. It remains to show that we can find good black paths. This will require a slight increase in the
probability $p$ with which we colour points of $Z$ black.

\subsection{Clusters of bad points}

To avoid bad points, we have to show that they are not very common. In fact, we shall prove (and
need) a much stronger statement, namely that `large clusters' of bad points essentially
never occur. (See below for a precise statement.)
Given $Z$,
a {\em $\delta$-bad component} is a component in the hypergraph on the $\delta$-bad points of $Z$ with edges given by
the $\delta$-close pairs and $\delta$-bad quadruples. In other words, we partition the bad points 
into minimal components such that every bad quadruple or close pair involves points from a single component.

\begin{lemma}\label{l_bad}
Let $\epsilon>0$ and $\eta>0$ be fixed, let $Z$ be a Poisson process of intensity $1$ on $\TT(s)$,
and suppose that $\delta=\delta(s)\le s^{-\epsilon}$.
Then \whp\ no $\delta$-bad component contains more than $\eta\log s$ vertices.
\end{lemma}

The proof will be a little long; it makes use of the purely deterministic statement below,
which will be proved in Section~\ref{a_l_deter} of the appendix.

\begin{lemma}\label{l_deter}
  Let $0<\eta<1/100$ be fixed.
  If $t$ is large enough, then it is impossible to
  arrange $M\ge \eta^{-6}t$ points $P_1,\ldots,P_M$ and $M$
  associated circles $C_1,\ldots,C_M$ in $\RR^2$ so that no two $P_i$
  are within distance $\eta$ of each other, the distance of every
  $P_i$ from the origin is between $\eta t$ and $t$, each $C_i$ has
  radius at most $t$, no $C_i$ contains any $P_j$, and $C_i$ passes
  within distance $\eta^3$ of both $P_i$ and the origin.
\end{lemma}

\begin{proof}[Proof of Lemma~\ref{l_bad}]
Throughout {\em bad} and {\em close} will mean $\delta$-bad and $\delta$-close respectively.
We shall write $f=O^*(g)$ to mean
that there is an absolute constant $C$ such that $f=O((\log s)^Cg)$ as $s\to\infty$.
Suppose the lemma fails for some particular values $\epsilon$, $\eta>0$, which we may assume to be
smaller than $1/10$.

Any bad quadruple or close pair lies in a disc of radius $2\sqrt{\log s}+\delta<\log s/2$, say.
(As noted earlier, throughout we use without comment the assumption that $s$ is larger than some
constant, depending on $\epsilon$ and $\eta$.)
Hence, if $z$ is a bad point in a bad component $C$ of more than
$\eta\log s$ vertices, the intersection of this component with a square $S$ of side $2(\log s)^2$ centred
on $z$
must contain at least $\eta\log s$ points of $C$: one cannot get further than distance $(\log s)^2$ from $z$
in fewer
than $\log s$ steps, each step passing from one point of a bad quadruple/close pair to another.

Fix a square $S$ of side $3(\log s)^2$. As we may cover $\TT$ with $O(s^2)$ such squares $S_i$, in such
a way that any square of side $2(\log s)^2$ lies within one $S_i$, it suffices
to prove
the following statement.
\begin{equation}\label{bs_disp}
 \Pr(\hbox{$S$ contains more than $\eta\log s$ bad points})=o(s^{-2}).
\end{equation}

The number $N$ of points of $Z$ in $S$ has a Poisson distribution with mean $\area(S)=9(\log s)^4$.
Let $B_1$ be the event that $N>10(\log s)^4$. Then
\begin{equation}\label{B1}
 \Pr(B_1)=\Pr\big(N>10(\log s)^4\big) =\exp\big({-}\Theta((\log s)^4)\big)= o(s^{-2}).
\end{equation}

Our first step will be to count close pairs in $S$.
In fact, we shall count $\delta'$-close pairs, where $\delta'=\delta^{1/10}$, so $\delta'$ is
much larger than $\delta$. Any upper bound on the number of $\delta'$-close pairs applies to $\delta$-close pairs as well.
Let $B_2$ be the event that $S$ contains more than $900/\epsilon^2$ pairs of $\delta'$-close points of $Z$.
Our first aim is to show that $\Pr(B_2)=o(s^{-2})$.
In doing this we shall condition on $N$; in the light of~\eqref{B1}, we may and shall assume that $N\le 10(\log s)^4$.

Given $N$, we may realize the positions of the $N$ points of $Z\cap S$
as a random sequence $z_1,\ldots,z_N$, where the $z_i$ are independent
and each is uniformly distributed in $S$. Consider placing the $z_i$ one by one: let $\CL_i$ be the
event that $z_i$ is $\delta'$-close to one of $\{z_1,\ldots,z_{i-1}\}$. 
As $i\le N\le 10(\log s)^4=O^*(1)$,
there are $O^*(1)$ previously placed points $z_i$ could be close to.
As the area $\delta'$-close to a given point is at most $\pi{\delta'}^2$,
whatever the positions of $z_1,\ldots,z_{i-1}$, the conditional probability of $\CL_i$, given
$N$ and the positions of $z_1,\ldots,z_{i-1}$,
is $O^*({\delta'}^2)$. As $N=O^*(1)$ by assumption, it follows that having placed some 
points $z_1,\ldots,z_j$,
the probability that there is an $i>j$ for which $\CL_i$ holds
is $O^*(N{\delta'}^2)=O^*({\delta'}^2)\le \delta'$.
Thus the probability that more than $30/\epsilon$ points $z_i$ are close to earlier points is 
at most ${\delta'}^{30/\epsilon}=\delta^{3/\epsilon}\le s^{-3}=o(s^{-2})$.
A point $z_i$ may be close to many earlier points, but is very
unlikely to be close to 
more than $30/\epsilon$, say -- the probability that $S$ contains any disc of radius $\delta'$ containing
$30/\epsilon$ points of $Z$ is $o(s^{-2})$.
It follows that with probability $1-o(s^{-2})$ there are at most $900/\epsilon^2$
pairs of $\delta'$-close points in $S$. In other words, $\Pr(B_2)=o(s^{-2})$.

We should like to count points in bad quadruples in a similar way. The problem is that it is not in general
true that, given three points, the probability that there is a fourth point forming a bad quadruple with them
is $O(\delta)$. If the first three points are very close together, the probability may in fact be rather large.
Some quadruples are easy to deal with, however; we start with these.

Let us say that a $\delta$-bad quadruple is {\em separated} if no two of its points are $\delta'$-close.
Let $x_1$, $x_2$, $x_3$ be any three points of $S$, and let $x$ be chosen from $S$ uniformly at random.
It will be straightforward to check that the probability that the set $\{x_1,x_2,x_3,x\}$ forms a
separated bad quadruple is $o(\delta')$. More formally, let us say that a quadruple $\{x_1,x_2,x_3,x_4\}$
is {\em weakly $\delta$-bad} if there is an $x\in \TT$ and an $r<2\sqrt{\log s}$ such that
$r<\dist{x}{x_i}\le r+\delta$ for every $i$. This is the same as the definition of $\delta$-badness,
except that we have omitted the condition that no other point of $Z$ is within distance $r$ of $x$.
In particular, $\delta$-badness implies weak $\delta$-badness.
\setcounter{claim_base}{18}
\setcounter{claim}{0}
\begin{claim}\label{cl_badsep}
Suppose that $x_1$, $x_2$, and $x_3$ are points of $S$, no two within distance
$\delta'$ of each other. Let $B$ be the set of $x\in S$ for which the quadruple $\{x_1,x_2,x_3,x\}$ is
weakly $\delta$-bad. Then $\area(B)=o(\delta')$.
\end{claim}

The straightforward proof of this purely deterministic statement is given in Section~\ref{a_cl_badsep}
of the appendix.

Let $B_3$ be the event that there are at least $M=200/\epsilon$ points in separated bad quadruples in $S$.
Recall that $\epsilon<1/10$ is constant, so $M$ is a constant with $M\ge 2000$. We claim
that $\Pr(B_3)=o(s^{-2})$. We condition on $N$, the number of points of $Z$ in $S$.
We shall assume, as we may, that $B_1$ does not hold, so
$N\le 10(\log s)^4$. Given such an $N$, construct $Z\cap S$ as $\{z_1,\ldots,z_N\}$ with the $z_i$
independent and uniformly distributed, as above.
If $B_3$ holds, then we can find subsets $\emptyset=Y_0,Y_1,\ldots,Y_r$ of $[N]=\{1,2,\ldots,N\}$
with $Y_i\subset Y_{i+1}$,
$1\le |Y_{i+1}\setminus Y_i|\le 4$ and $M\le|Y_r|\le M+3$, where $Y_{i+1}$ is the union of $Y_i$ with a quadruple
$\{a,b,c,d\}$ such that $\{z_a,z_b,z_c,z_d\}$ is a separated bad quadruple. (We can continue the sequence
from $Y_i$ to $Y_{i+1}$ using the indices $\{a,b,c,d\}$ of any separated
bad quadruple such that $\{a,b,c,d\}\not\subset Y_i$.
By assumption such a quadruple must exist unless $|Y_i|\ge M$, when we stop.)

It follows that if $B_3$ holds, then
we can find a subset $I$ of $[N]$ of size at most $M+3$, an order on $I$, and a subset $J$ of 
$I$ of size $r\ge M/4$, so that if we add the points $\{z_i:i\in I\}$ in this particular order,
then for each $j\in J$ the point $z_j$ 
completes a separated bad quadruple with three points of $\{z_i:i\in I\}$ already placed.
Fixing $I$, the order, $J$, $j\in J$, and the three earlier points, the probability that
$z_j$ completes the specified separated bad quadruple
is $o(\delta')=o(\delta^{1/10})$, using uniform distribution of $z_j$ and Claim~\ref{cl_badsep}.
Hence, the probability that a given $z_j$, $j\in J$, completes a separated bad quadruple is
$o(\binom{|I|}{3}\delta^{1/10})=o(\delta^{1/10})$, and the probability
that every $z_j$, $j\in J$, does so is $o(\delta^{M/40})$.
Crudely, the number of choices for $I$ is at most $4N^{M+3}\le (\log s)^{5M}$. The number of
choices for the order is $|I|!\le |I|^{|I|}\le M^{2M}$, while the number of choices for $J$ is at most $2^{2M}$.
Hence,
\[
 \Pr(B_3\mid B_1^c) \le \left(4(\log s)^5M^2\delta^{1/40}\right)^M \le \delta^{M/50} =o(s^{-2}).
\]
Since $\Pr(B_1)=o(s^{-2})$, it follows that $\Pr(B_3)=o(s^{-2})$ as claimed.

We have now shown that there is a constant $C$ (depending on the parameters
$\epsilon$ and $\eta$ appearing in the statement of Lemma~\ref{l_bad}) such that with
probability $1-o(s^{-2})$ the square $S$ contains (a) at most $C$
pairs that are $\delta'$-close, and (b) at most $200/\epsilon$ points in separated $\delta$-bad
quadruples. To complete the proof of \eqref{bs_disp} it suffices to show that with high
enough probability there are at most $\eta\log s/2$ points $z\in Z\cap S$ with the property
that $z$ is not in a $\delta'$-close pair, but is in a $\delta$-bad quadruple of points of $Z\cap S$
containing a $\delta'$-close pair $\{x_1,x_2\}$.

Let $B_4$ be the event
that there are more than $\eta\log s/2$ such points $z$.
If $B_4$ holds but $B_2$ does not, then there is a $\delta'$-close pair $\{x_1,x_2\}\subset Z\cap S$
and a set $D$ of points of $Z\cap S$ so that every $z\in D$ is in a $\delta$-bad quadruple with $\{x_1,x_2\}$,
where $D$ contains at least
$\eta'\log s$ points, $\eta'=\eta'/(2C)$.
From the definition of $\delta$-badness,
each point of $D$ is within distance $4\sqrt{\log s}+2\delta\le
5\sqrt{\log s}$ of $x_1$. By Lemma~\ref{l_disc}, there is a constant $\eta_2>0$ such that the
event $B_5$ that some disc in $\TT$ of radius $\eta_2\sqrt{\log s}$ contains more than $\eta'\log s/2$
points of $Z$ has probability $\Pr(B_5)=o(s^{-2})$. Assuming that $B_5$ does not hold, we see that $D$
contains a set $D'$ of at least $\eta'\log s/2$ points whose distance from $x_1$ is between $\eta_2\sqrt{\log s}$
and $5\sqrt{\log s}$.

We shall next show that for some constant $\eta_3$, with very high probability
not many points in $D'$ are $\eta_3$-close to other points in $D'$.
The argument is as at the beginning of the lemma for $\delta'$-close points,
but we must consider a total area of the right order of magnitude.

Fix any disc $X$ of radius $6\sqrt{\log s}$. From Lemma~\ref{l_disc},
there is a constant $C'$ such
that the event $B_1'(X)$ that $X$ contains more than $C'\log s$ points of $Z$ has probability $o(s^{-4})$.
Let $\eta_3>0$ be constant. Arguing as at the start of the proof of this lemma,
as we place points one by one in $X$, assuming that $B_1'(X)$ does not hold, at each stage the probability that
the new point $z_i$ is within distance $\eta_3$ of an earlier point is $O(\eta_3^2)$:
the area $\eta_3$-close to some earlier point is at most $\pi\eta_3^2C'\log s$, while $z_i$ is chosen
uniformly from $X$, a domain of area $\area(X)=36\pi\log s$.
Hence the number $Y(X)$ of points $z_i$
$\eta_3$-close to an earlier point
is dominated by a binomial distribution $\Bi(n,p)$ with parameters $n=C'\log s$ and $p=O(\eta_3^2)$.
Choosing $\eta_3$ small enough, and using the weak bound $\Pr(\Bi(n,p)\ge t)\le (enp/t)^t$,
the probability that $Y(X)$ exceeds $\eta'\log s/4$ is $o(s^{-4})$.

Let $B_6$ be the event that some disc $X'$ of radius $5\sqrt{\log s}$ meeting $S$ contains 
at least $\eta'\log s/4$ points $z_i\in Z$ each
closer than $\eta_3$ to some earlier point of $Z$ in $X'$ (taking a random order on the points of $Z\cap X'$
as above). Then 
(placing $O(s^2)$ discs of radius $6\sqrt{\log s}$ so that any disc of radius $5\sqrt{\log s}$
is contained in one), we
have $\Pr(B_6)=o(s^{-2})$. Assuming that $B_6$ does not hold, and deleting from $D'$ each point $z_i$ closer
than $\eta_3$ to an earlier point, we find a set $D''\subset D'$ in which any two points are at least $\eta_3$ apart,
with $|D''|\ge \eta'\log s/4$.

Let us write $D''=\{P_1,\ldots,P_M\}$ where $M\ge \eta'\log s/4$, and set $t=5\sqrt{\log s}$.
Then any two $P_i$ are at least $\eta_3$ apart, and the distance from any $P_i$ to $x_1$
is between $\eta_2\sqrt{\log s}=\eta_2t/5$ and $t$. Furthermore, as each $P_i$ is in a $\delta$-bad quadruple
with $x_1$ and $x_2$,
for each $P_i$ there is a circle $C_i$ with radius at most $2\sqrt{\log s}<t$
passing within distance $\delta=o(1)$ of $P_i$ and $x_1$,
and containing no points of $Z$, and hence no $P_j$.
For $s$ large enough, the existence of points $P_i$ and circles $C_i$ with the above properties
contradicts the deterministic Lemma~\ref{l_deter},
applied with $t=5\sqrt{\log s}$ and $\eta=\min\{\eta_2/5,\eta_3,1/100\}=\Theta(1)$.

In summary, if $S$ contains many bad points, then
one of the `bad' events $B_i$, $1\le i\le 6$, holds. Each has probability $o(s^{-2})$, so with probability
$1-o(s^{-2})$, $S$ contains fewer than $\eta\log s$ bad points, proving \eqref{bs_disp}.
As noted at the start of the proof,
this suffices to prove Lemma~\ref{l_bad}.
\end{proof}

\begin{figure}[htb]
\centering
\input{badeg.pstex_t}
\caption{Three points $\{x_1,x_2,x_3\}$
of $Z$ are very close to $x$. If one places one more point $z_i$ of $Z$ in each small circle
in an arbitrary way, and $C$ is otherwise free of points of $Z$, then each $z_i$ forms a bad quadruple
with $\{x_1,x_2,x_3\}$; the circle $C_1$ illustrates this for $z_1$.}
\label{fig_badeg}
\end{figure}
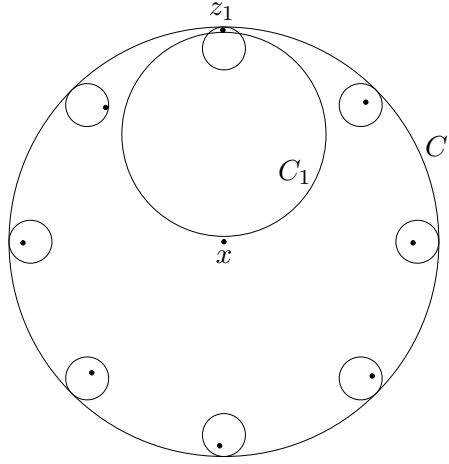

The last part of the proof above may seem unnecessarily complicated, but the situation
is a little delicate. One might expect that
there is some argument that the size of all bad components will be at most a constant, as (locally) each
individual `badness' has probability $n^{-\Theta(1)}$, but this is not
true. In fact, with high probability, $\TT(s)$ will contain a bad component with $\Theta(\log s/\log\log s)$
vertices: consider three points very close together, and $k$ points spaced almost evenly around
a circle of radius $1$, say. Each of these $k$ points can be placed anywhere in a small disc of radius
$\Theta(k^{-2})$ so that all $k$ points form bad quadruples with the three central points;
this is illustrated for $k=8$ in Figure~\ref{fig_badeg}.
Such a configuration will occur somewhere if $k=c\log s/\log \log s$ for a small enough constant $c$.
Thus the bound $o(\log s)$ given by Lemma~\ref{l_bad} is not far from best possible.

\subsection{Bad points can be avoided}
In the next lemma, for $p_1<p_2$ we couple the probability measures (on black/white-coloured Poisson
processes) $\Pr_{p_1}^{\TT(s)}$ and $\Pr_{p_2}^{\TT(s)}$. In other words, we shall construct
simultaneously Poisson processes $Z_1$ and $Z_2$ of intensity $1$
on $\TT(s)$, with associated colourings
$\col_1$, $\col_2$, so that given $Z_i$ each point of $Z_i$ is coloured black with probability $p_i$
independently of all other points of $Z_i$, and white otherwise.

Let us write $\Egl$ for the global event that the following three conditions hold:

 (a) there is a bijection $\phi$ from the points of $Z_1$ to those
of $Z_2$ such that $\phi(z)$ is black whenever $z$ is,

 (b) each point $z$ of $Z_1$ is within distance $1$ of the corresponding point
$\phi(z)$ of $Z_2$, and

 (c) if $z$, $z'$ are black points of $Z_1$ whose Voronoi cells share an edge,
 then there is a sequence $\phi(z)=z_1,z_2,\ldots,z_t=\phi(z')$, $t=O(\log s)$,
 of black, $\delta$-good points of $Z_2$, such that the Voronoi cells of
 $z_i$ and $z_{i+1}$ with respect to $Z_2$ share an edge.

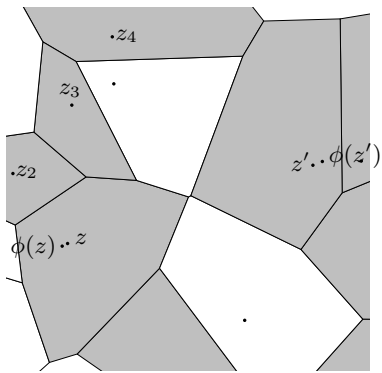
\begin{figure}[htb]
\centering
\input{couple.pstex_t}
\caption{Part of the Voronoi tiling with respect to $Z_2$, together with two points $z$, $z'$ of $Z_1$
whose Voronoi cells with respect to $Z_1$ meet; these cells are not drawn.}
\label{fig_couple}
\end{figure}

The event $\Egl$ is illustrated in Figure~\ref{fig_couple}: part of the Voronoi tiling with respect
to $Z_2$ is shown, together with the points $z$, $z'$ of $Z_1$. For simplicity the figure is drawn
with
$\phi(z'')=z''$ for all $z''\in Z_1\setminus\{z,z'\}$, so $z$ and $z'$ are the only points
that move when we pass from $Z_1$ to $Z_2$. The Voronoi cells of $z$ and
$z'$ with respect to $Z_1$, which are not shown, meet. The points $z_i$, with $z_1=\phi(z)$
and $z_5=\phi(z')$, exhibit the path required for condition (c).

As usual, in the following statement $o(1)$ notation refers to $s\to\infty$ with all
other parameters (here $p_1$, $p_2$ and $\epsilon$) fixed. The statement is not very appealing, but will
allow us to deduce Theorem~\ref{th_couplerobust}, at which point the lemma can be forgotten.

\begin{lemma}\label{l_coupling}
Let $0<p_1<p_2<1$ and $\epsilon>0$ be given. Let $\delta=\delta(s)$ be any function
with $\delta(s)\le s^{-\epsilon}$.
We may construct in the same probability space
coloured Poisson processes $(Z_1,\col_1)$ and $(Z_2,\col_2)$ on $\TT(s)$ as described above,
in such a way that the global event $\Egl$ defined above holds \whp.
\end{lemma}

\begin{proof}
An outline of the idea is as follows: we first decide roughly where the points $z_i$ of $Z$ are.
Then in most cases
we can already see that a particular $z_i$ is $\delta$-good. There will be some clusters of points
that may be $\delta$-bad; however, these will have size $o(\log s)$. Different clusters will behave
independently when we decide the precise location of the $z_i$, and each will contain some bad point with
only small probability. We can allow for this by using the difference in the two colouring probabilities.

We now turn to the details.
As usual, throughout the proof we assume that $s$ is larger than some very large constant. We may assume
that $\delta=s^{-\epsilon}$.

Pick $\delta_1$ with $\delta_1\sim \delta^{1/5}=s^{-\epsilon/5}$ so that $s/\delta_1$ is an integer,
and partition $\TT(s)$ into $K=(s/\delta_1)^2$ squares $S_k$ of side $\delta_1$ in the natural way.
Set $\delta_2=\delta_1^{1/10}$. Note that $\delta << \delta_1 << \delta_2 \sim s^{-\epsilon/50}$.

Let $N$ be a Poisson random variable with mean $s^2$. For $i=1,2$ we shall realize $Z_i$ as 
$\{z_{i,1},\ldots,z_{i,N}\}$, where the $z_{i,j}$ are uniformly distributed on $\TT$, and
for each $i$, the $z_{i,j}$, $1\le j\le N$, are independent.
For $1\le j\le N$, choose $s_j$ independently and uniformly at random from $[K]=\{1,2,\ldots,K\}$. We shall place
both $z_{1,j}$ and $z_{2,j}$ into the square $S_{s_j}$. As long as the locations of each $z_{i,j}$
within $S_{s_j}$ are chosen uniformly at random, with independence as $j$ varies and $i$ is fixed,
we shall obtain coupled Poisson processes $Z_1$, $Z_2$ so that there is a bijection
$\phi$ satisfying condition (b) of $\Egl$.

Having chosen the $s_j$, we say that a pair $\{a,b\}\subset [N]$ is {\em potentially close}
if it is possible, given that $z_{2,j}\in S_{s_j}$ for all $j$, 
that $z_{2,a}$ and $z_{2,b}$ are $\delta_2$-close, i.e., if some two points of the squares
$S_{s_a}$ and $S_{s_b}$ are within distance $\delta_2$.
Note that in this case, wherever $z_{2,a}$ and $z_{2,b}$ are placed in their squares,
they are certainly $(\delta_2+2\sqrt{2}\delta_1)$-close and hence $(2\delta_2)$-close.

A quadruple $\{a,b,c,d\}\subset [N]$ is {\em potentially bad} if it is possible, given
that $z_{2,j}\in S_{s_j}$ for all $j$, 
that $Q=\{z_{2,a},z_{2,b},z_{2,c},z_{2,d}\}$ forms a $\delta_2$-bad quadruple. This time it follows
that wherever the points $z_{2,j}$ are within their squares, $Q$ is a $(2\delta_2)$-bad quadruple.
Indeed, by definition there is a placement of the $z_{2,j}$ so that $Q$ is a $\delta_2$-bad quadruple.
Let $x\in \TT$ and $r$
be the point and radius witnessing this, so no $z_{2,j}$ is within distance $r$ of $x$ and
all four points of $Q$ are within distance $r+\delta_2$. Moving each $z_{2,j}$ arbitrarily
within $S_{s_j}$, each $z_{2,j}$ moves by at most a distance $\sqrt{2}\delta_1$. Reducing $r$
by $\sqrt{2}\delta_1$ (stopping at $0$ if $r$ is smaller than this to start with!), we find a smaller
$r'$ so that none of the relocated $z_{2,j}$ is within distance $r'$ of $x$, but all four
relocated points of $Q$ are within distance $r+\delta_2+\sqrt{2}\delta_1\le r'+\delta_2+2\sqrt{2}\delta_1
\le r'+2\delta_2$ of $x$.

In analogy with previous definitions we say that a point $z_{2,a}$, or the
corresponding index $a\in [N]$,
is {\em potentially bad}
if it is in a potentially close pair or potentially bad quadruple, and define {\em potentially bad components}
as components of the hypergraph induced on the set of potentially bad points by the potentially close
pairs/bad quadruples. Note that these definitions all depend only on $N$ and the $s_j$, i.e.,
on which squares each $z_{2,j}$ lies in. As we place $z_{1,j}$ and $z_{2,j}$ into the same square $S_{s_j}$,
we could have written $z_{1,j}$ instead of $z_{2,j}$ throughout all definitions of potential badness.

We claim that for any fixed $\eta>0$, with probability $1-o(1)$ no potentially bad component
contains more than $\eta\log s$ points. This follows from Lemma~\ref{l_bad},
applied with $2\delta_2\le s^{-\epsilon/51}$ in place of $\delta$: we just realize the exact
locations of the $z_{2,j}$, and apply Lemma~\ref{l_bad} to see that \whp\ there is no component
of more than $\eta\log s$ $(2\delta_2)$-bad points. But every potentially bad pair/quadruple
gives rise to a $(2\delta_2)$-bad pair/quadruple,
so the largest $(2\delta_2)$-bad component is certainly at least as large
as the largest potentially bad component.

As outlined above, we construct the coupling as follows. First, we decide $N$, and the $s_j$.
Later we shall position $z_{1,j}$ and $z_{2,j}$ in $S_{s_j}$ appropriately. Also, we shall
define the colourings so that whenever $\col_1(z_{1,j})$ is black, so is $\col_2(z_{2,j})$.

Let $B_1$ be the `bad' event that there is some disc $X$ of radius
$2\sqrt{\log s}$ in $\TT$ so that every square $S_k$ contained entirely in $X$ is free of points
of $Z_1$ (and hence of $Z_2$).
Then $\Pr(B_1)=o(1)$ by Lemma~\ref{l_disc}: if $B_1$ holds, then,
when we place the points $z_{1,j}$ randomly in the corresponding $S_{s_j}$, there is always a smaller
disc within $X$, of radius at least $2\sqrt{\log s}-\sqrt{2}\delta_1\ge \sqrt{\log s}$, containing
no points of $Z_1$, a Poisson process of intensity $1$ on $\TT$. The existence of such a disc
has probability $o(1)$ by Lemma~\ref{l_disc}, so $\Pr(B_1)=o(1)$.
Note that if $B_1$ does not hold, then every disc of radius $2\sqrt{\log s}$ will contain points of
both $Z_1$ and $Z_2$, however we complete the coupling.
Let $B_2$ be the event that there is a potentially bad component of more than $\eta\log s$ points,
where $\eta>0$ is a very small constant to be chosen later.
As shown above, $\Pr(B_2)=o(1)$. Note that $B_1$ and $B_2$ depend only on $N$ and the $s_j$.

\label{page_B1}
If $B_1$ or $B_2$ holds, we abandon the construction and complete the coupling any old way; $\Egl$ will
not hold in this case, but this does not matter, as $\Pr(B_1\cup B_2)=o(1)$. Suppose from now on
that $N$ and the $s_j$ are fixed and that $B_1\cup B_2$ does not hold.

Let $\{a,b\}\subset [N]$ be a pair not in the same potentially bad component. Thus $a$ and $b$
are not potentially close, and are not two of the four points in a potentially bad quadruple.
We claim that, no matter how we position the $z_{i,j}$ in their squares $S_{s_j}$,
if the Voronoi cells of $z_{1,a}$ and $z_{1,b}$ with respect to $Z_1$ are adjacent,
then the cells of $z_{2,a}$ and $z_{2,b}$ with respect to $Z_2$ will also be adjacent.
Indeed, consider any placements $Z_1$, $Z_2$ consistent with our choice of squares.
By definition of potential badness,
$z_{1,a}$ and $z_{1,b}$ are not $\delta_2$-close, and do not belong to a $\delta_2$-bad quadruple.
Now consider the midpoint $M$ of the common edge of the relevant cells of $Z_1$.
By Lemma~\ref{l_rejoinbad}, no point $z_{1,c}$ of $Z_1\setminus\{z_{1,a},z_{1,b}\}$ is within distance
$\dist{z_1}{M}+\delta_2^6$ of $M$. Moving from $Z_1$ to $Z_2$, all points move a distance of at most $\sqrt{2}\delta_1<\delta_2^6/2$,
so $z_{2,a}$ and $z_{2,b}$ are the two points of $Z_2$ closest to $M$. It follows by Lemma~\ref{l_closest}
that the Voronoi cells of $z_{2,a}$ and $z_{2,b}$ with respect to $Z_2$ meet.

We are now ready to construct the coupling; the observation above suggests that we consider each potentially bad component
separately. Let us start with the trivial components: suppose that $j\in [N]$ is not potentially bad.
Then we couple the $z_{i,j}$ and their colourings in the following natural way: take $z_{1,j}=z_{2,j}$
uniformly random on $S_{s_j}$, colour $z_{2,j}$ black with probability $p_2$ independently of everything
else, and colour $z_{1,j}$ black only if $z_{2,j}$ is black, and then only with conditional probability $p_1/p_2$.

Next, let $C$ be a non-trivial potentially bad component. We start by coupling the $z_{i,j}$, $j\in C$,
in the natural way as above; we shall have to adjust the coupling slightly.
Consider the event $B(C)$ that there is a pair $\{a,b\}\subset C$ so that $z_{2,a}$ and $z_{2,b}$ are $\delta$-close,
or a quadruple $\{a,b,c,d\}\subset C$ so that
$\{z_{2,a},z_{2,b},z_{2,c},z_{2,d}\}$ forms a $\delta$-bad quadruple. Note that here we consider
$\delta$-badness, a much stronger condition than $\delta_2$-badness.

We claim that $\Pr(B(C)) \le s^{-\epsilon/6}$. Indeed, considering quadruples first,
there are at most $|C|^4=O^*(1)$ quadruples
$\{a,b,c,d\}$ to consider, using the assumption that $B_2$ does not hold. For each, the joint distribution
of $(z_{2,a},z_{2,b},z_{2,c},z_{2,d})$ is that these points are chosen independently and uniformly at random from
certain squares $S_{s_a}$, etc. Recall that each square $S_i$ has side-length $\delta_1$, and
that $\delta_1\sim \delta^{1/5}$. In particular, $\delta_1$ is much larger than $\delta$.
It is easy to check that no matter how the squares are located, the probability that the quadruple is
$\delta$-bad is $O^*(\delta/\delta_1^4)=O^*(\delta_1)$. Indeed, one argument goes as follows.
If the quadruple is $\delta$-bad, then
there is a point $x$ and a radius $r$ such that $r$ and both coordinates of $x$ are multiples of $\delta$,
so that all four points are in the annulus centred at $x$ with inner radius $r$ and outer radius $r+5\delta$.
(We find $x$ and $r$ close to the $x'$ and $r'$ witnessing $\delta$-badness.) For any such annulus,
its intersection with each $S_j$ has area $O(\delta\delta_1)$, so the probability that all four points
lie in one given annulus is $O((\delta\delta_1/\delta_1^2)^4)=O(\delta^4/\delta_1^4)$. But the number of choices for the annulus
is $O^*(\delta^{-3})$; the range of the coordinates of $x$ and of $r$ is at most $O(\sqrt{\log s})$,
as the definition of $\delta$-badness limits $r$ to $O(\sqrt{\log s})$.
Hence any particular quadruple is $\delta$-bad with probability $O^*(\delta/\delta_1^4)=O^*(\delta_1)$,
and the expected number of $\delta$-bad quadruples is $O^*(\delta_1)$. The argument for pairs is much
simpler: each of the $O^*(1)$ pairs is $\delta$-close with probability at most $O(\delta\delta_1/\delta_1^2)
=O(\delta_1^4)$.
Thus $\Pr(B(C)) = O^*(\delta_1)\le s^{-\epsilon/6}$, as claimed.

Let $G(C)$ be the event that every point in the component is white with respect to $(Z_1,\col_1)$ and
black with respect to $(Z_2,\col_2)$. From the way the coupling is defined, $\Pr(G(C))=(p_2-p_1)^{|C|}$.
As $\epsilon$, $p_1$ and $p_2$ are constants, while $|C|\le \eta\log s$,
we have $\Pr(G(C))\ge 2s^{-\epsilon/6}$ provided we chose $\eta$ small enough.
Hence $\Pr(G(C)\setminus B(C))\ge \Pr(B(C))$.

Let $G'(C)\subset G(C)\setminus B(C)$ be an event with probability exactly $\Pr(B(C))$. Our final coupling is defined by
`crossing over' from the natural coupling using $B(C)$ and $G'(C)$: deleting the portion of the probability space
in which $B(C)$ or $G'(C)$ holds, we add instead the distribution of $(Z_1,\col_1)$ on $B(C)$ coupled in an arbitrary way
with that of $(Z_2,\col_2)$ on $G'(C)$, and vice versa.
[More formally, let the natural coupling described above be given by random variables 
$\tZ_i(\omega)$, $\tcol_i(\omega)$, $i=1,2$, $\omega\in \Omega$, on a probability space $(\Omega,\Pr)$.
Then $B(C)$ and $G'(C)$ are subsets of $\Omega$ with the same ($\Pr$-)measure. Let $f$
be an arbitrary measure-preserving bijection from $B(C)\cup G'(C)$ to itself, mapping $B(C)$
into $G'(C)$ and vice versa. Such a map exists if $\Omega$ is defined suitably.
We may define the crossed-over coupling $(Z_i,\col_i)$, $i=1,2$, on the same space $(\Omega,\Pr)$
by setting $(Z_1,\col_1)$ equal to $(\tZ_1,\tcol_1)$, and setting
\[
 Z_2(\omega)= \left\{\begin{array}{ll}
   \tZ_2(\omega) & \omega\notin B(C)\cup G'(C),\\
   \tZ_2(f(\omega)) & \omega\in B(C)\cup G'(C), \end{array}\right.
\]
and defining $\col_2$ similarly.]
 
In the final coupling, the partition of the probability space
into $B(C)$, $G'(C)$, and the rest has the following properties:

\label{cases}
1. When $\omega\notin B(C)\cup G'(C)$, then
$z_{1,j}=z_{2,j}$ for all $j\in C$, there is no $\delta$-close pair/bad quadruple in $C$, and if $\col(z_{1,j})$ is black, so is
$\col(z_{2,j})$.

2. When $\omega\in B(C)$, then
$Z_1$ contains a $\delta$-close pair or $\delta$-bad quadruple in $C$, while
[as $f(\omega)\in G'(C)$] $Z_2$ does not, and $\col_2(z_{2,j})$ is black for every $j\in C$.

3. When $\omega\in G'(C)$, then [as $f(\omega)\in B(C)$]
$Z_2$ contains a $\delta$-close pair or $\delta$-bad quadruple in $C$ while $Z_1$ does not, and $\col_1(z_{1,j})$ is white for every $j\in C$.

We extend the coupling from individual components $C$ to the whole of $Z_1$, $Z_2$ by taking the product of the relevant
measures, i.e., treating different components independently.

Note that we have indeed coupled the correct distributions for the $(Z_i,\col_i)$: the `natural' coupling clearly does
this, and the marginal distributions are preserved by the crossing over. From the way we constructed the coupling,
placing points always in the same squares $S_j$, conditions (a) and (b) of the event $\Egl$ hold: note that in each
case above $z_{2,j}$ is black whenever $z_{1,j}$ is.
It remains to check condition (c).
For this we shall use the following claim.
\setcounter{claim_base}{20}
\setcounter{claim}{0}

\begin{claim}\label{cl_pathincpt}
Consider the coupling defined above. Suppose that $B_1$ does not hold, that
$a,b\in C$ for some potentially bad component $C$, and that $z_{1,a}$ and $z_{1,b}$ are black
points whose Voronoi cells with respect to $Z_1$ meet.
Then there is a piecewise-linear path $P$ joining $z_{2,a}$ to $z_{2,b}$ so that,
with respect to $Z_2$,
every point of $P$ lies in the Voronoi cell of some $z_{2,c}$, $c\in C$, where $z_{2,c}$ is black and $\delta$-good.
\end{claim}

The straightforward proof of this deterministic statement is given in Section~\ref{a_cl_pathincpt} of the appendix.

Claim~\ref{cl_pathincpt}
contains everything we need to finish the proof of Lemma~\ref{l_coupling}: it remains only to establish
the existence and properties of the path (of cells) required by condition (c) in the definition of $\Egl$.
Note first that from the properties of the three cases 1,2,3 above, whenever a point is black in $Z_1$ it is black in $Z_2$
and $\delta$-good in $Z_2$. (Recall that if $j\in C$ is in a $\delta$-close pair/bad quadruple in $Z_2$,
 then this pair/quadruple
is contained in $C$, as $C$ is a component in the hypergraph of potentially close pairs/bad quadruples.)

Suppose that $z_{1,a}$ and $z_{1,b}$ are black points of $Z_1$ whose Voronoi cells
meet. Then $z_{2,a}$ and $z_{2,b}$ are black and $\delta$-good. 
We have already shown above that if $a$ and $b$ do not lie in the same potentially bad
component, then the Voronoi cells of $z_{2,a}$ and $z_{2,b}$ must also meet, and (c) holds. The claim
gives (c) in the case
that $a$ and $b$ do lie in the same potentially bad component $C$,
using $|C|=O(\log s)$ and noting that if there is an $x$-$y$ path in the union of certain Voronoi cells,
then there is such a path visiting each of the cells at most once.
\end{proof}

We now have all the tools we need to prove Theorem~\ref{th_couplerobust}.

\begin{proof}[Proof of Theorem~\ref{th_couplerobust}]
We use the coupling guaranteed by Lemma~\ref{l_coupling}, applied with $(4\delta)^{1/6}$ in place of $\delta$.
We shall assume that every disc of radius $\sqrt{\log s}$ contains points of both $Z_1$ and $Z_2$; this holds
\whp\ by Lemma~\ref{l_disc}.
Given $P_1$, there is a sequence of (not necessarily distinct) black points $z_{1,i}$, $1\le i\le t$, in $Z_1$ so that
the Voronoi cells of $z_{1,i}$ and $z_{1,i+1}$ share an edge for each $i$, and so that $P_1$
meets the Voronoi cells corresponding to each $z_{1,i}$ and lies in the union of these cells.
From Lemma~\ref{l_coupling}, \whp\ for each $i$ there is a sequence $S_i$
of $(4\delta)^{1/6}$-good black points $z^{(i)}_{2,j}$
of $Z_2$ so that for each $j$ the Voronoi cells of $z^{(i)}_{2,j}$ and $z^{(i)}_{2,j+1}$
with respect to $Z_2$ share an edge, where $S_i$ starts with the unique point $\phi(z_{1,i})\in Z_2$ corresponding
to $z_{1,i}$, ends with $\phi(z_{1,i+1})$, and has length $O(\log s)$.
Piecing together the sequences $S_i$, we obtain a single sequence $z_{2,j}$, $1\le j\le t'$.
By Corollary~\ref{c_goodsolid}, the piecewise-linear path $P_2$ formed by the $z_{2,j}$ and the midpoints of the
edges where the cells of consecutive $z_{2,j}$ meet is $4\delta$-robustly black with respect to $Z_2$.
Finally, as all Voronoi cells have diameter $O(\sqrt{\log s})$, each $S_i$ has length at most $O(\log s)$,
and $z^{(i)}_{2,1}=\phi(z_{i,1})$ is within distance $1$ of $z_{1,i}$, the paths
$P_1$ and $P_2$ are geometrically `close' to each other as claimed.
\end{proof}

The entire purpose of the present section was to develop a tool (Theorem~\ref{th_couplerobust})
that will enable us to deal with problems arising when we discretize the continuous Poisson process $Z$.
As noted in the introduction, these problems would not arise in, for example, the random discrete
Voronoi setting.

\section{Crossing a rectangle when $p>1/2$}\label{sec_above}

Using Corollary~\ref{c_RSWcon}, Theorem~\ref{th_sharpMOD} and Theorem~\ref{th_couplerobust},
we shall prove that for any $p>1/2$,
some large enough long thin rectangle can be crossed the long way
by a black path with probability at least 99\%.
Once we have this, it will be easy to prove our main results, Theorems~\ref{th_doesperc} and~\ref{th_decay}.
Recall that $f_p(\rho,s)$ is the probability that $H(R)$ holds, i.e., that there is a (piecewise-linear)
black path crossing $R$ from left to right,
when $R$ is a $\rho s$ by $s$ rectangle, in the probability measure $\Pr_p=\Pr_p^{\RR^2}$.

\begin{theorem}\label{th_above}
Let $p>1/2$, $c_1<1$, $\rho>1$ and $s_1$ be given. There is an $s>s_1$ such that $f_p(\rho,s)>c_1$.
\end{theorem}

\begin{proof}
It suffices to prove the result for any fixed $\rho>1$, for example, $\rho=3$. The general result follows
by using positive correlation (Lemma~\ref{l_cor}) as in the proof of equation \eqref{e3}.

Throughout we fix $p>1/2$, $c_1<1$, $\rho=3$ and $s_1$. 
All implicit constants in our notation ($O(.)$, etc) will depend on $p$ and $c_1$ only.
We shall choose constants $\epsilon$ and $s_0$ as follows. We first fix $\epsilon=\epsilon(p,c_1)>0$
small enough that a certain condition we shall encounter in the course of the proof holds; this condition does
not involve $s_0$ (or $s$). Then we choose an $s_0$ depending on all parameters chosen so far,
large enough that the statements `... provided $s$ is sufficiently large' in what follows hold for all
$s\ge s_0$. (In principle, we could give an explicit value.) We assume throughout that $s_0\ge 6s_1$.

By Corollary~\ref{c_RSWcon}, applied with $\rho=10$, there is an absolute constant $c_0>0$
such that for some $s\ge s_0$ we have with probability at least $c_0$ a horizontal (i.e., long)
black crossing of $R_1=[0,10s/13]\times [0,s/13]$ in the $p=1/2$ random Voronoi tiling in $\RR^2$.
Let us fix such an $s$ throughout the proof.

Let us consider the fixed rectangle $R_1$ as a rectangle in the torus $\TT=\TT(s)$.
Let $E_1$ be the event that $R_1$ has a black crossing $P_1$ in the Voronoi tiling defined on the torus.
As noted in Section~\ref{sec_torus}, since $R_1$ does not come close to `wrapping round' the torus,
if $s$ is large enough we have $\Pr_{1/2}^\TT(E_1)=\Pr_{1/2}^{\RR^2}(H(R_1)) + o(1) \ge c_0/2$.

Let $\delta$ be chosen so that $s^{-\epsilon} \le \delta \le (1+o(1)) s^{-\epsilon}$,
and so that $s/\delta$ is an integer.
Set $p'=(1/2+p)/2$, so $1/2<p'<p$.
Consider the coupling given by Theorem~\ref{th_couplerobust}, applied with $p_1=1/2$, $p_2=p'$, and
$\epsilon/2$ in place of $\epsilon$, noting that $\delta\le s^{-\epsilon/2}$ if $s$ is large enough.
Suppose that $E_1$ holds with respect to the process $(Z_1,\col_1)$; as noted above,
this event has probability at least $c_0/2$. Suppose also that the global event described in Theorem~\ref{th_couplerobust}
holds, as it does with probability $1-o(1)$.
Then with respect to $(Z_2,\col_2)$ there is a $4\delta$-robustly black path $P_2$ geometrically close
to $P_1$. In particular, if $s$ is large enough, $P_2$ remains entirely within the
rectangle $R_1'$ with the same centre as $R_1$ but with corresponding sides
$2(\log s)^2<s/1000$ longer. Furthermore, $P_2$ contains points within distance $s/1000$ of the left and right
sides of $R_1$. Hence, a sub-path $P_2'$ of $P_2$ crosses $R_2$ from left to right, 
where $R_2$ is the $3s/4$ by $s/12$ rectangle with the same centre as $R_1$; see Figure~\ref{fig_RP12}.

\begin{figure}[htb]
\centering
\input{RP12.pstex_t}
\caption{The rectangles $R_1$ and $R_2$, together with a black
path $P_1$ (shown solid) crossing $R_1$ horizontally, and the nearby path $P_2$ (shown dotted),
part of which crosses $R_2$ horizontally.}
\label{fig_RP12}
\end{figure}
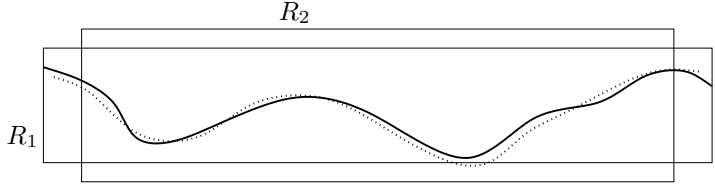

To summarize, let $R_2$ be the fixed $3s/4$ by $s/12$ rectangle in $\TT$ described above, and let 
$E_2$ be the event that there is a $4\delta$-robustly black path $P$ crossing $R_2$ horizontally.
As the assumptions made above hold simultaneously with probability at least $c_0/2-o(1)$,
and these assumptions imply that $E_2$ holds for $(Z_2,\col_2)$,
we have
\[
 \Pr_{p'}^\TT(E_2)\ge c_0/3,
\]
if $s$ is large enough.

Our main aim is to apply Theorem~\ref{th_sharpMOD} to an appropriate symmetric event in a suitable discrete
probability space. Let $E_3$ be the event that
there is a $4\delta$-robustly black path crossing {\em some} $3s/4$ by $s/12$ rectangle in $\TT$
horizontally.
Clearly $E_3$ is `symmetric', and
\[
 \Pr_{p'}^\TT (E_3) \ge \Pr_{p'}^\TT (E_2) \ge c_0/3.
\]

The next step is to define a corresponding symmetric event in a discrete probability space.
Divide $\TT=\TT(s)$ up into $(s/\delta)^2$ small squares $S_i$ of side-length $\delta$ in the natural way.
The {\em crude state} of each $S_i$ will be {\em bad} if $S_i$ contains at least one white point of $Z$ (and perhaps
one or more black points), {\em neutral} if $S_i$ contains no points of $Z$, and {\em good} if $S_i$ contains
at least one black point of $Z$ and no white points. If we generate $Z$ as a Poisson process with intensity $1$
and colour each point independently black with probability $p''$ bounded away from $0$ and $1$
and white with probability $1-p''$,
then writing $\gamma=\delta^2=o(1)$ for the area of $S_i$, each $S_i$ is bad/neutral/good with respective probabilities
\begin{eqnarray}
 \pbad &=& 1-\exp\big({-}\gamma(1-p'')\big) \sim \gamma(1-p''), \nonumber \\
 \pneut &=& \exp(-\gamma),\label{pbng} \\
 \pgood &=& \exp\big({-}\gamma(1-p'')\big)\big(1-\exp(-\gamma p'')\big) \sim \gamma p'',\nonumber
\end{eqnarray}
independently of all other $S_j$.
The {\em crude state} $CS$ of $\TT$ is given by the crude states of the $S_i$.

Define $\Etdc$ to be the event that the crude state $CS$ of $\TT$ is such that $E_3$ is possible given $CS$.
As $\Pr(E_3)=\E(\Pr(E_3\mid CS)) \le \Pr(\Etdc)$, we have
\begin{equation}\label{lb}
 \Pr_{p'}^\TT (\Etdc)\ge c_0/3
\end{equation}
if $s$ is large enough.

If we write $n=(s/\delta)^2$ for the number of squares $S_i$, and take $p_-=\pbad$ and $p_+=\pgood$,
then $\Etdc$ may be regarded as an event $\Etd$ in the finite probability space
$W_{p_-,p_+}^n$ defined in Section~\ref{sec_ingred}. Here the $i$th coordinate in
the state space $\{-1,0,1\}^n$ of $W_{\pbad,\pgood}^n$
is $-1$, $0$ or $1$ if the crude state of $S_i$ is bad, neutral or good respectively.

Now $E_3$ is `symmetric' in the sense that it is preserved under translations of the torus. It follows
that $\Etdc$ is preserved under horizontal and vertical translations of the torus by integer multiples of $\delta$.
Such translations can
map any $S_i$ into any other $S_j$. Thus $\Etd$ is symmetric with respect to the coordinates
$[n]=\{1,2,\ldots,n\}$ of the product space $W_{\pbad,\pgood}^n$:
there is a permutation group acting transitively on $[n]$ whose
induced action on $\{-1,0,1\}^n$, the state space of $W_{\pbad,\pgood}^n$, preserves $\Etd$.
This is precisely the symmetry condition required
for Theorem~\ref{th_sharpMOD}.

We claim that $\Etd$ is increasing in $\{-1,0,1\}^n$.
In other words,
we claim that $\Etdc$ is preserved if we change the crude
state of one $S_i$ from bad to neutral
or from neutral to good. To see this,
suppose that we are given $CS$ (i.e., the crude states of all the $S_j$) such that $\Etdc$ holds. Fix one $S_i$ and
suppose that this $S_i$ is neutral. By the definition
of $\Etdc$, the event $E_3$ is possible given this overall crude state.
Changing the state of $S_i$ to good clearly preserves the possibility of $E_3$, as it corresponds
to adding one or more black points, and $E_3$ is black-increasing. Suppose now that the state of $S_i$ is bad,
and pick a particular realization of $(Z,\col)$ consistent with $CS$ for which $E_3$ holds.
In $Z$, there is (at least one) white point $w$ in $S_i$; there may also be one or more black points $b$.
To show that $E_3$ is still possible after changing $S_i$ to neutral, it suffices to check that $E_3$ still
holds after deleting all points of $Z$ in $S_i$. As $E_3$ is black-increasing, the only possible problem is that deleting
the black points destroys $E_3$. But this cannot happen, because of the point $w$, and the fact that
a black point $b$ within distance $2\delta$ of some white point cannot help $E_3$ to hold -- such 
a black point cannot be the closest point of $Z$ to a $4\delta$-robustly black point of $\TT$.

We are ready to apply Theorem~\ref{th_sharpMOD}.
Let $\pbadn0$ and $\pgoodn0$ be the values defined by \eqref{pbng} with $p''=p'$,
and let $\pbadn1$ and $\pgoodn1$ be the values defined by \eqref{pbng} with $p''=p$.
From the correspondence between $W_{\pbad,\pgood}^n$ and the crude state of the torus,
we have
\begin{equation}\label{equiv0}
 \Pr_{\pbadn0,\pgoodn0}^n(\Etd) = \Pr_{p'}^\TT(\Etdc)
\end{equation}
and
\begin{equation}\label{equiv1}
 \Pr_{\pbadn1,\pgoodn1}^n(\Etd) = \Pr_{p}^\TT(\Etdc).
\end{equation}
From \eqref{lb} and \eqref{equiv0} it follows that
\begin{equation}\label{Wlb}
 \Pr_{\pbadn0,\pgoodn0}^n(\Etd) \ge c_0/3.
\end{equation}

Now
\begin{equation}\label{pbadn0}
 \pbadn0 \sim (1-p')\gamma \hbox{\quad and\quad} \pgoodn0 \sim p'\gamma,
\end{equation}
while
\begin{equation}\label{pbadn1}
  \pbadn1\sim (1-p)\gamma \hbox{\quad and\quad} \pgoodn1\sim p\gamma.
\end{equation}
We shall apply Theorem~\ref{th_sharpMOD}
with $n=(s/\delta)^2$, and with $p_-=\pbadn0$, $p_+=\pgoodn0$, $q_-=\pbadn1$ and $q_+=\pgoodn1$.
Recall that $p'$ and $p$ are constants with $p'<p$, and note that
\begin{equation}\label{psim}
 q_+-p_+,\,p_--q_- \sim (p-p')\gamma.
\end{equation}
Let $c_2<1$ be any absolute constant, 
and set
\[
 \eta = \min\{c_0/4,1-c_2\},
\]
an absolute constant.
From~\eqref{pbadn0} and~\eqref{pbadn1}, all four of $p_-,p_+,q_-,q_+$ are at most $\gamma$ if $s$ is large enough.
Hence the quantity $\pmax$ appearing in \eqref{e_sharpMOD} is at most $\gamma$, and the
right hand side of \eqref{e_sharpMOD} is at most
\[
 \Delta = c_3\log(1/\eta) \gamma\log(1/\gamma)/\log n,
\]
where $c_3$ is an absolute constant.
As $n=(s/\gamma)^2\ge s^2$ and $\gamma\ge s^{-2\epsilon}$, we have
$\log(1/\gamma)/\log n\le \epsilon$. It follows that
$\Delta\le C\epsilon\gamma$ for some absolute constant $C$. By our choice of $\epsilon$ we may ensure
that $C\epsilon<(p-p')/2$, so $\Delta<(p-p')\gamma/2$.
Hence, \eqref{psim} implies that if $s$ is large enough, then $q_+-p_+$ and $p_--q_-$ are both at least
$\Delta$. Thus, by \eqref{Wlb} and Theorem~\ref{th_sharpMOD},
\[
 \Pr_{\pbadn1,\pgoodn1}^n(\Etd) \ge c_2.
\]
Hence, from \eqref{equiv1}, we have $\Pr_p^{\TT}(\Etdc) \ge c_2$.

In summary, we have shown that for any absolute constant $c_2<1$
there is some $s\ge s_0\ge 6s_1$ such that  $\Pr_{p}^{\TT} (\Etdc) \ge c_2$.

Consider any crude state $CS_0$ such that $\Etdc$ holds,
and let $(Z_1,\col_1)$ be a realization of $(Z,\col)$ consistent with $CS=CS_0$ and such that $E_3$ holds,
i.e., such that there is a $4\delta$-robustly black path $P$ crossing some $3s/4$ by $s/12$ rectangle in $\TT$ horizontally.
By the definition of $\Etdc$, such a $(Z_1,\col_1)$ exists.
Let $(Z_2,\col_2)$ be any other realization of $Z$ consistent with $CS=CS_0$.
In $Z_1$, for every point $x$ of $P$ there is a black point $b$ at some distance $r$ and no white point
within distance $r+4\delta$. In particular, the square $S_i$ containing $b$ contains no white points
and is thus good.
As $Z_2$ corresponds to the same crude state, there is (at least one) black
point of $Z_2$ in $S_i$; this point is within distance $r+2\delta$ of $x$. Also,
there is no white point $w\in Z_2$ within
distance $r+2\delta$ of $x$: otherwise, the square $S_j$ containing $w$ is bad, and there is a white point $w'$
of $Z_1$ in this square, and hence within distance $r+4\delta$ of $x$, contradicting our assumption.
Hence every point $x$ of $P$ is black with respect to $Z_2$.

Let $E_4$ be the event that there is a black path crossing {\em some} $3s/4$ by $s/12$ rectangle in $\TT$
horizontally.
We have shown that whenever $CS$ is such that $\Etdc$ holds, the conditional probability, given $CS$, that $E_4$ holds 
is $1$. Hence,
\[
 \Pr_p^{\TT}(E_4) \ge \Pr_p^{\TT}(\Etdc) \ge c_2.
\]

We can cover $\TT$ with a fixed number $M$ of $s/2$ by $s/6$ rectangles $R_i$ so that whenever $E_4$ holds, there
is a black path crossing some $R_i$ horizontally. For example, we may take $M=48$ and use rectangles whose $x$-coordinates are all
multiples of $s/4$ and whose $y$-coordinates are all multiples of $s/12$.
Recall that $f_{p}(3,s/6)$ is the probability in the random Voronoi percolation on $\RR^2$ that an $s/2$ by $s/6$
rectangle has a horizontal black crossing.
As before the corresponding probability in the torus is within $o(1)$. The events $H_i$ that $R_i$ has
a horizontal black crossing are black-increasing. Hence, their complements $H_i^c$ are white-increasing.
Applying Lemma~\ref{l_cor}, which holds in the torus, with black and white exchanged,
each $H_i^c$ is positively correlated
with $\cap_{j<i} H_j^c$, so
\[
 \Pr_p^{\TT}\left(\bigcap_{i=1}^{M} H_i^c\right) \ge \prod_{i=1}^M\Pr_p^{\TT}(H_i^c) = \Pr_p^{\TT}(H_1^c)^M = (1-f_p(3,s/6)-o(1))^M.
\]
If no $H_i$ holds then $E_4$ fails, so the probability of the intersection above is at most $1-c_2$,
and
\[
 1-f_p(3,s/6) \le (1-c_2)^{1/M} +o(1).
\]
(This is sometimes known as the `square-root' trick.)
Choosing $c_2$ sufficiently close to $1$, we see that for some
$s\ge s_0\ge 6s_1$ we have
\[
 f_p(3,s/6) > c_1.
\]
This completes the proof of Theorem~\ref{th_above}.
\end{proof}

\section{Proofs of the main results}\label{sec_proofs}

Now all the technicalities are behind us. To prove our main results, we only need the ingredients
in Section~\ref{sec_ingred}, Theorem~\ref{th_above} and Lemma~\ref{l_cor}.

We first prove Theorem~\ref{th_doesperc}, that for any $p>1/2$ we have $\theta(p)>0$,
i.e., that percolation occurs in the black cells of the random
Voronoi tiling in which cells are coloured black independently with probability $p$. This corresponds
to Kesten's result for bond-percolation in $\Z^2$. It turns out that it is easy to deduce
Theorem~\ref{th_doesperc} from Theorem~\ref{th_above} using the idea of $1$-dependent percolation,
and in particular the simple observation given as Lemma~\ref{l_1dep}.

\begin{proof}[Proof of Theorem~\ref{th_doesperc}]
Let $p>1/2$ be fixed.
Let $p_0<1$ be a value for which Lemma~\ref{l_1dep} holds, and set $\epsilon=(1-p_0)/4$.
For $s>1$, let $R_s$ be a $3s$ by $s$ rectangle.
Let us recall the definition of the event $\Es(R_s)$ considered in Lemma~\ref{l_small}:
writing $R_s[d]$ for the set of all points within Euclidean distance $d$ of $R_s$,
$\Es(R_s)$ is the event that for every $x\in R_s[r]$ there is 
some point $z\in Z$ at distance $\dist{x}{z}<r$, where $r=r(s)=2\sqrt{\log s}$.
By Lemma~\ref{l_small}, there is an $s_1$ such that $\Pr_p(\Es(R_s))\ge 1-\epsilon$ for all $s\ge s_1$.
Also, the event $\Es(R_s)$ depends only on the restriction of $Z$ to $R_s[2r]$.
We may and shall assume that $s_1$ is large enough that $r=r(s)<s/4$ for $s\ge s_1$.

By Theorem~\ref{th_above}, there is an $s\ge s_1$ such that $f_p(3,s)\ge 1-\epsilon$.
From now on we fix such an $s$ and suppress the dependence on $s$ in our notation.
Thus,
for $R$ a $3s$ by $s$ rectangle with the long side horizontal, the probability of the event $H(R)$ that there is a black
path within $R$ crossing $R$ from left to right is at least $1-\epsilon$.
Let $S_1$ and $S_2$ be the two $s$ by $s$ squares making up the left and right thirds of $R$.
As $H(R)$ implies $H(S_1)$, we have $\Pr_p(V(S_i))=\Pr_p(H(S_i))\ge 1-\epsilon$ for $i=1,2$.

Let
\begin{equation}\label{GRdef}
 G(R) = H(R) \cap V(S_1) \cap V(S_2) \cap \Es(R).
\end{equation}
This event is illustrated in Figure~\ref{fig_GR}.

\begin{figure}[htb]
 \[\epsfig{file=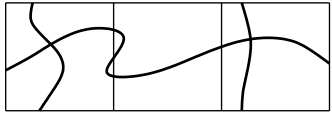,height=0.6in}\]
\caption{A $3s$ by $s$ rectangle $R$ such that $H(R)$, $V(S_1)$ and $V(S_2)$ hold.}
\label{fig_GR}
\end{figure}
\noindent

For $R$ an $s$ by $3s$ rectangle we define $G(R)$ similarly: $G(R)$ is the event that $G(R')$ holds
after rotating $\RR^2$ through 90 degrees mapping $R$ to $R'$.
As each of the four events intersected in \eqref{GRdef} has probability at least $1-\epsilon$, we have
\[
 \Pr_p(G(R))\ge 1-4\epsilon = p_0
\]
for any $3s$ by $s$ or $s$ by $3s$ rectangle $R$.

We define a $1$-dependent bond percolation measure $\Prone$ on $\Z^2$
as follows: for $x,y\in \Z$ the edge $e$ from $(x,y)$ to $(x+1,y)$
is open in $\Prone$ if and only if $G(R_e)$ holds in $\Pr_p$ for the $3s$ by $s$ rectangle $R_e=[2sx,2sx+3s]\times [2sy,2sy+s]$.
Similarly, the edge $e$ from $(x,y)$ to $(x,y+1)$ is open in $\Prone$ if and only if $G(R_e)$ holds in $\Pr_p$ for the $s$
by $3s$ rectangle $R_e=[2sx,2sx+s]\times [2sy,2sy+3s]$.
\begin{figure}[htb]
 \[\epsfig{file=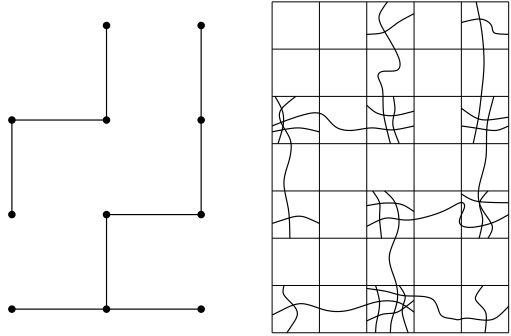,height=1.75in}\]
\caption{A set of open edges in $\Prone$ (left), and corresponding rectangles $R$ drawn with black paths showing that
$G(R)$ holds in $\Pr_p$.}
\label{fig_1indep}
\end{figure}

It is easy to check that $\Prone$ is $1$-dependent. Indeed,
recall that the event $\Es(R)$ depends only on the restriction of $(Z,\col)$
to the $2r$-neighbourhood $R[2r]$ of $R$, and that if $\Es(R)$ holds, then the colours of all points in $R$
are determined by this restriction. It follows that the event $G(R)$ depends only on this restriction.
As $r<s/4$, if $S$ and $T$ are sets of edges of $\Z^2$ at graph distance at least $1$,
then the
regions $D_U=\bigcup_{e\in U} R_e[2r]$, $U=S,T$, are disjoint. Hence the restrictions of $(Z,\col)$ to $D_S$ and $D_T$
are independent, implying 1-dependence of $\Prone$.

By Lemma~\ref{l_1dep} we have percolation in $\Prone$, i.e., with positive $\Prone$-probability
the origin is in an infinite open path in $\Z^2$.
However,  we have defined $G(R)$ in such a way
that a $\Prone$-open path in $\Z^2$ guarantees a corresponding (much longer)
black path in the coloured Voronoi tiling of $\RR^2$,
using only the fact that horizontal and vertical crossings of a square must meet; see Figure~\ref{fig_1indep}.
Hence, with positive $\Pr_p$-probability,
some point of $[0,s]\times [0,s]$ is in an infinite path of black Voronoi cells. It follows that $\theta(p)>0$,
completing the proof of Theorem~\ref{th_doesperc}.
\end{proof}

We next turn to our second main result, Theorem~\ref{th_decay},
giving exponential decay of the size (defined in any reasonable way)
of $C_0$, the black component of the origin, when $p<1/2$. Again the result follows easily from
Theorem~\ref{th_above} using the idea of locally dependent percolation. This time,
we use the observation given as Lemma~\ref{l_kdepneg}.

\begin{proof}[Proof of Theorem~\ref{th_decay}]
The basic idea is to divide $\RR^2$ up into $s$ by $s$ squares, and to show that for each
square $S$ it is very unlikely
that any point of $S$ is connected to any point far outside $S$. Then we shall apply Lemma~\ref{l_kdepneg}.

Let $S$ be an $s$ by $s$ square.
We can surround $S$ with a ring of four overlapping $3s$ by $s$ (or $s$ by $3s$)
rectangles $R_i$, $1\le i\le 4$, as shown in Figure~\ref{fig_wann}.

\begin{figure}[htb]
 \[\epsfig{file=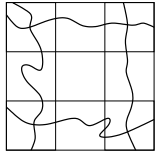,height=0.8in}\]
\caption{A square $S$ surrounded by four rectangles $R_i$, containing four white paths $P_i$
(drawn in black) whose union must contain a cycle surrounding $S$.}
\label{fig_wann}
\end{figure}

\def\Linf{L_\infty}
\def\Sss{{\overline S}}
\def\Sssr{{\overline S}[2r]}
\def\Svssr{{\overline S_v}[2r]}
\def\Swssr{{\overline S_w}[2r]}
\def\Svs{{\overline S_v}}

Let $\Sss=S\cup\bigcup_{i=1}^4 R_i$ be the $\Linf$ $s$-neighbourhood of $S$.
Suppose that each $R_i$ is crossed the long way by a {\em white} path $P_i$. (As before, this means that
there is a piecewise-linear path of white points inside $R_i$ starting and ending on appropriate sides of $R_i$.)
Then the paths $P_i$ meet pairwise, so their union contains a piecewise-linear white cycle $C$ surrounding $S$,
as in Figure~\ref{fig_wann}.
Given the existence of such a cycle $C$, it follows with probability 1 that no point
$x\in S$ is connected to any point $x'$ outside $\Sss$ by a black path $B$.
(Recall that a point may be both black and white, if it is in two Voronoi cells of different colours. But
this is not a problem -- with probability $1$ the Voronoi cells meet three at a vertex. Corresponding to the white
cycle $C$ and to the black path $B$ there is a cycle/path of white/black Voronoi cells, and these cannot cross.)

Let $L(S)$ be the event that some point in $S$ is connected by a black path to a point in the boundary
of $\Sss$. As each $R_i$ is crossed the long way by a white path with probability $f_{1-p}(3,s)$,
and these events are white-increasing and hence positively correlated by Lemma~\ref{l_cor},
we find a white cycle $C$ as above
with probability at least $f_{1-p}(3,s)^4$, and
\begin{equation}\label{LS}
 \Pr_p(L(S)) \le 1 - f_{1-p}(3,s)^4.
\end{equation}

Let $G(S)= L(S) \cup \Es(\Sss)^c$, and set
$\epsilon=p_1(7)/2$, where $p_1(7)$ is the quantity appearing in Lemma~\ref{l_kdepneg}
for $k=7$, i.e., for $7$-dependent site percolation.
Applying Lemma~\ref{l_small} to the $3s$ by $3s$ square $\Sss$,
we have $\Pr_p(\Es(\Sss))\to 1$ as $s \to \infty$, so if $s$ is large enough, then $\Pr_p(\Es(\Sss)^c)<\epsilon$.
Also, by Theorem~\ref{th_above} and \eqref{LS}, there are arbitrarily large $s$ for which
$\Pr_p(L(S)) < \epsilon$.
It follows that there is some $s$ for which
\begin{equation}\label{pgsm}
 \Pr_p(G(S)) < 2\epsilon = p_1(7).
\end{equation}
Let us fix such an $s$, chosen large enough that $r=r(3s)=2\sqrt{\log(3s)}<s/4$.

As before, by Lemma~\ref{l_small}, $\Es(\Sss)$ depends only on the restriction of $(Z,\col)$ to the
(Euclidean) $2r$-neighbourhood $\Sssr$ of $\Sss$.
If $\Es(\Sss)$ does not hold, then $G(S)$ holds. On the other hand, if $\Es(\Sss)$ does hold, then $G(S)$ holds
if and only if $L(S)$ holds, which depends only on the colours of points of $\Sss$. Hence,
$G(S)$ depends only on the restriction of $(Z,\col)$ to $\Sssr$.
To each $v=(a,b)\in \Z^2$ let us assign an $s$ by $s$ square $S_v=[as,as+s]\times [bs,bs+s]$.
We say that $v\in \Z^2$ is {\em open} if $G(S)$ holds. We claim that
this assignment of states to vertices defines a $7$-dependent site percolation measure $\Prone$ on $\Z^2$.
Indeed,
if $v$, $w\in \Z^2$ are at graph distance at least $7$, then either their $x$-coordinates or their $y$-coordinates
differ by at least $4$, so the sets $\Svssr$ and $\Swssr$ are disjoint. For
$A$, $B\subset \Z^2$ at distance at least $7$, which vertices of $A$ are open and which of $B$ 
depend respectively on the restrictions of $(Z,\col)$ to the domains $D_A$, $D_B$, where $D_C=\bigcup_{v\in C}\Svssr$,
so $D_A$ and $D_B$ are disjoint. Hence the states of vertices in $A$ are independent of the sates of the vertices
in $B$, proving $7$-dependence.

Let $C_0'$ be the set of vertices $v\in \Z^2$ in the open cluster of the origin in the site percolation
$\Prone$ on $\Z^2$ we have just defined. From \eqref{pgsm}, each vertex $v$ is open with probability at most
$p_1(7)$,
so by Lemma~\ref{l_kdepneg} there is some $c>0$ such that
\[
 \Pr_p(|C_0'|\ge m) = \Prone(|C_0'| \ge m) \le \exp(-cm)
\]
for all $m\ge 1$.

Now let $C_0$ be the black component of the origin in the random Voronoi percolation, i.e., the maximal connected union
of black Voronoi cells containing the origin (which is empty if the origin is white). Let $N=|C_0^G|$
be the number of Voronoi cells in $C_0$.
We consider two cases. The first is that $C_0$ is contained in a disc $D$ of radius $10s$ centred at the origin.
Note that in this case $\area(C_0)$ and $\diam(C_0)$ are bounded by absolute constants.
It is easy to check that the number of Voronoi cells contained entirely in $D$ decays exponentially (the number
is dominated by the number of points of $Z$ in the disc, a Poisson distribution with a certain fixed expectation).
Hence, there is a $c'>0$ such that
\begin{equation}\label{decaysmall}
 \Pr\big(C_0\subset D\hbox{ and } N\ge n\big) \le \exp(-c'n)
\end{equation}
for all $n$.

On the other hand, if $C_0\not\subset D$, then $C_0$ has geometric diameter at least $10s$,
and for every point $x \in C_0$
there is another point $x'\in C_0$ at distance at least $5s$.
It follows that if $x\in C_0\cap S_v$, then $G(S_v)$ holds; indeed,
$x$ is a point of $S_v$ connected by a black path in $C_0$
to a point $x'$ outside $\Svs$. Stopping the first time the path hits the boundary of $\Svs$,
we see that $L(S_v)$ holds, and hence so does $G(S_v)$.
Let $C_1=\{v\in \Z^2:C_0\cap S_v\ne \emptyset\}$.
We have shown that if $C_0\not\subset D$, then $G(S_v)$ holds for every $v\in C_1$,
i.e., that every $v\in C_1$ is open.
But $C_1$ is with probability $1$ a connected set in $\Z^2$ that, if nonempty, contains the
origin, so $C_1\subset C_0'$, and we have
\begin{equation}\label{c1}
 \Pr(|C_1|\ge m)\le \exp(-cm)
\end{equation}
for all $m$, with the same $c>0$ as above.
As the area of $C_0$ is at most $s^2|C_1|$, this establishes exponential decay of the area of $C_0$.
As $\diam(C_0)\le \sqrt{2}s+s(|C_1|-1)$, say, exponential decay of $\diam(C_0)$ also follows.

Finally, to show exponential decay of the number $N=|C_0^G|$ of Voronoi cells in $C_0$
note that, given that $|C_1|=m$, there are are most $(4e)^m$ possibilities
for $C_1$. As every cell in $C_0$ is contained entirely in the union $U$ of the $m$ squares $S_v$, $v\in C_1$,
the number of cells in $C_0$ is at most the number of points of $Z$ in $U$, which has area $s^2m$.
For any possible $C_1$, the probability that the corresponding $U$ contains more than $100s^2m$ points of $Z$ is certainly
at most $\exp(-10s^2m)\le \exp(-10m)$. (Recall that $s$ is a large constant.)
Multiplying by the number of possibilities for $C_1$,
we see that
\[
 \Pr\big(C_0\not\subset D,\hbox{ }N\ge 100s^2 m \hbox{ and } |C_1|=m\big)\le \exp(-5m).
\]
From \eqref{c1} it follows that
\[
 \Pr\big(C_0\not\subset D\hbox{ and }N\ge 100s^2m\big) \le \exp(-5m)+\exp(-cm).
\]
Combined with \eqref{decaysmall} this establishes exponential decay of $N=|C_0^{G}|$, completing the proof.
\end{proof}

As noted in Section~\ref{sec_thms}, our proofs of Theorems~\ref{th_doesperc} and~\ref{th_decay}
also give a new proof of Theorem~\ref{th_doesnt}, Zvavitch's result that $\theta(1/2)=0$.
As the deduction of Theorem~\ref{th_doesnt} from our intermediate results is rather short, we give it here.

\begin{proof}[Proof of Theorem~\ref{th_doesnt}]
The result follows almost immediately from Corollary~\ref{c_RSWcon}, applied in a similar way to the Russo-Seymour-%
Welsh Theorem for $\Z^2$.
Let $A_s$ be the square annulus centered on the origin, with inner and outer diameters $8s$ and $10s$.
Let $C(A_s)$ be the event that $A_s$ contains a (piecewise-linear) white cycle surrounding the origin.
Arguing as in the proof of Theorem~\ref{th_decay},
from Corollary~\ref{c_RSWcon}, applied with $\rho=10$ and with white and black swapped,
and positive correlation (Lemma~\ref{l_cor}), there is an absolute constant $\epsilon=c_0(10)^4/2>0$
and an unbounded set $\mathcal S$ of values of $s$
for which $\Pr_{1/2}(C(A_s))\ge c_0(10)^4=2\epsilon$.
Let $L(A_s)$ be the event that $\Es(R)$ holds for each of the four rectangles $R$ making up $A_s$. By
Lemma~\ref{l_small}, $\Pr_{1/2}(\Es(R))\to 1$ as $s\to \infty$, so $\Pr_{1/2}(L(A_s))\ge 1-\epsilon$ if $s$
is large enough.
Thus, deleting small values from $\mathcal S$ if necessary, the event
$G(A_s)=C(A_s)\cap L(A_s)$ holds with probability at least $\epsilon$ for every $s\in{\mathcal S}$.
Finally, as in the proof of Theorem~\ref{th_doesperc}, it follows from Lemma~\ref{l_small}
that $G(A_s)$ depends only on the restriction of $(Z,\col)$ to the $4\sqrt{\log(10s)}$ neighbourhood $N_s$
of $A_s$.
If $s_2\ge 2s_1$ and $s_1$ is large enough, then $N_{s_1}$ and $N_{s_2}$ are disjoint.
As $\mathcal S$ is unbounded,
we can pick an infinite sequence $s_i$, $i=1,2,\ldots$, of values from ${\mathcal S}$ with $s_1$ large enough
and $s_{i+1}\ge 2s_i$ for every $i$.
Then the events $G(A_{s_i})$ are independent. As each has probability at least $\epsilon$,
with probability $1$ some $G(A_{s_i})$ holds. But then $C(A_{s_i})$ holds, and the origin is surrounded
by a white cycle, so the black component $C_0$ of the origin is bounded.
This proves that $\theta(1/2)=0$.
\end{proof}

The arguments given in this section were very simple. Indeed, once one thinks of using
$k$-dependent percolation measures, the main results are more or less immediate from
Theorem~\ref{th_above}. Also, the basic idea needed to derive Theorem~\ref{th_above}
from Theorem~\ref{th_RSW}, our analogue of the RSW Theorem, is very simple:
construct a suitable symmetric event and apply a sharp-threshold result such as
Theorem~\ref{th_sharpMOD}. Together, these
ideas give a very simple proof of, for example, Kesten's results on the critical
probability for bond percolation in $\Z^2$; we shall return to this in future work~\cite{ourKesten,ourKesten2}.

\section{Appendix}

In this appendix we give the proofs of several deterministic statements used in earlier sections.
Some of these statements stand alone and will be repeated here; others only make sense in the context in
which they arise.

\subsection{Proof of Lemma~\ref{l_rejoinbad}}\label{a_l_rejoinbad}

We start with a little geometric lemma concerning distances in $\RR^2$. As before,
for two points $A$, $B$ of $\RR^2$ we write $\dist{A}{B}$ for the Euclidean distance between them.

\begin{lemma}\label{l_geom}
Let $C_1$ and $C_2$ be two circles with centres
$O_1$, $O_2$ and radii
$r_1$, $r_2$, with $r_1,r_2\le 1/\delta$, where $0<\delta<1/100$. 
Suppose that $C_1$, $C_2$ meet in two points $P_1$, $P_2$ with $\dist{P_1}{P_2}\ge \delta$,
and that $\dist{O_1}{O_2}\ge \delta$.
Let $C_3$, $C_4$ be the circles with centres $P_1$, $P_2$ and radii $\delta/4$,
and let $M$ be the midpoint of the line-segment $O_1O_2$.
Then every point outside $\cup_{i=1}^4 C_i$ is at distance at least $\dist{P_1}{M}+\delta^6$
from $M$.
\end{lemma}

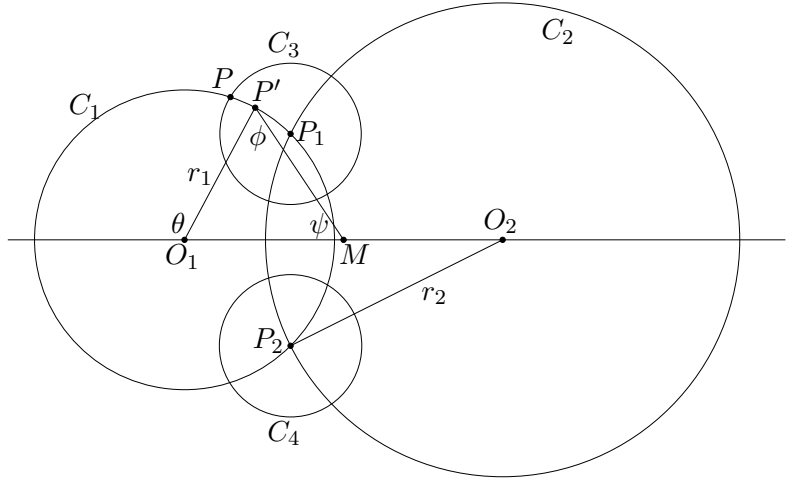
\begin{figure}[htb]
\centering
\input{fig1.pstex_t}
\caption{An illustration of the notation in Lemma~\ref{l_geom}. The circles
$C_3$ and $C_4$ are drawn larger than their true size for clarity.}
\label{figCi}
\end{figure}

\begin{proof}
We may assume that the centres of $C_1$, $C_2$ are $(-a,0)$ and $(a,0)$, $a\ge\delta/2$,
so $M$ is the origin, and that $P_1=(c,d)$ with $d>0$. In fact, as $\dist{P_1}{P_2}\ge \delta$,
we have $d\ge \delta/2$, so $C_3$ lies
entirely in the upper half plane.
Consider a point $P=(x,y)$ not in the interior of any $C_i$ at minimal distance from $M$.
By symmetry we may assume $y\ge 0$.

Note that $M$ lies inside at least one of $C_1$, $C_2$, as otherwise these circles cannot meet
in two points. Hence $P\ne M$, so $P$ must lie on at least one of the $C_i$.
By minimality of $\dist{P}{M}$, either $P$ is the closest point of some $C_i$ to $M$,
or $P$ lies on at least two $C_i$ (otherwise, there is a point very close to $P$ on the same circle
that is closer to $M$).
The first case is easy to exclude: let $Q_i$ be the closest
point of $C_i$ to $M$. Whether or not $M$ lies
inside $C_1$, $Q_1$ is the rightmost point of $C_1$, i.e., the point
$(-a+r_1,0)$. This lies strictly inside $C_2$, as otherwise $C_1$ and $C_2$ do not meet in two points.
So we cannot have $P=Q_1$. Similarly, $P\ne Q_2$.
Now $Q_3$ is a point strictly between $P_1$ and $M$ along the line-segment $P_1M$. (Recall
that $C_3$ lies entirely in the upper half plane, so $M$ lies outside $C_3$.)
As $M$ is inside $C_i$ for $i=1$ or $i=2$, and
$P_1$ is on the boundary of this $C_i$, the point $Q_3$ is inside $C_i$, so we cannot have $P=Q_3$.

It follows that $P$ lies on two of the circles $C_i$. As $C_4$ lies in the lower half plane and, clearly,
$P\ne P_1, P_2$, we have without loss of generality that $P$ lies on $C_1$ and $C_3$. Note that $P$ lies
anticlockwise around the upper-half of $C_1$ from $P_1$; the upper half of $C_1$ on the other
side of $P_1$ lies inside $C_2$. Let us consider moving
a point $P'=(x',y')$ along $C_1$, starting from $P_1$ and moving anticlockwise towards $P$.
As we move, since $y'\ge 0$, the distance $\dist{P'}{M}$ from $M$ increases monotonically, and it suffices
to show that it increases by at least $\delta^6$. To get to $P$ we must move a distance at least
$\delta/4$ along $C_1$; note that $y'\ge d-\delta/4\ge \delta/4$ holds throughout.

Now consider the angles $\phi=O_1P'M$, $\psi=O_1MP'$ and $\theta$ shown in Figure~\ref{figCi}.
Note that if $a=r_1$, i.e., if $M$ is on $C_1$, then $\phi=\psi=\theta/2$.
Hence if $a>r_1$, then $\phi\ge\theta/2$. On the other hand, if $a<r_1$, then $\psi\ge \theta/2$,
so, by the sine rule, $\sin\phi=(a/r_1)\sin\psi \ge (a/r_1)\sin(\theta/2)$.
Using $a\ge\delta/2$, $r_1\le 1/\delta$, it follows that in either case
\[
 \sin\phi\ge \delta^2\sin(\theta/2)/2.
\]
Now $y'\ge \delta/4$ while $r_1\le 1/\delta$, so $\sin(\theta/2)\ge \sin\theta/2\ge\delta^2/8$,
and $\sin\phi\ge \delta^4/16$.

Finally, the rate at which the distance $\dist{P'}{M}$ increases as we move $P'$ at rate $1$ around $C_1$
is exactly $\sin\phi$, so after moving distance $\delta/4$ starting at $P'=P_1$ we have increased $\dist{P'}{M}$
by at least $\delta^5/64\ge \delta^6$. Thus
$\dist{P}{M}\ge \dist{P_1}{M}+\delta^6$, completing the proof.
\end{proof}

It is now easy to deduce Lemma~\ref{l_rejoinbad}. We recall the statement.
\setcounter{theorem}{16}
\addtocounter{theorem}{-1}

\begin{lemma}
Let $s\ge e^{100}$ and $\delta<1/\log s$ be given, and
suppose that no disc of radius $2\sqrt{\log s}$ in $\TT$ is free of points of $Z$.
Let $z_1$, $z_2\in Z$ be points whose Voronoi cells share an edge, and let $M$ be the midpoint of that edge.
If some $z\in Z\setminus\{z_1,z_2\}$ is at distance less than $\dist{z_1}{M}+\delta^6$ from $M$, then
either $\{z_1,z_2\}$ is a $\delta$-close pair, or
there is a $\delta$-bad quadruple $\{z_1,z_2,z,z'\}$, $z'\in Z$.
\end{lemma}

\begin{proof}[Proof of Lemma~\ref{l_rejoinbad}]
Note first that the assumption concerning point-free discs ensures that all Voronoi cells have
diameter much less than $s$. In particular, they are all contractible in the torus, and are hence
convex polygons. Thus the cells of $z_1$ and $z_2$ do meet in a line-segment $O_1O_2$, say.
Furthermore, there are distinct points $x_1$, $x_2\in Z\setminus\{z_1,z_2\}$ so that $O_1$ is a vertex of the cell
of $x_1$ and $O_2$ of the cell of $x_2$. Thus, 
$\dist{O_1}{z_1}=\dist{O_1}{z_2}=\dist{O_1}{x_1}=r_1$, say, and $\dist{O_2}{z_1}=\dist{O_2}{z_2}=\dist{O_2}{x_2}=r_2$, say,
and no point of $Z$ is closer to $O_i$ than distance $r_i$. The notation is shown in Figure~\ref{fig_rejoinbad}.
As $x_1\ne x_2$, we may assume without loss of generality that $z\ne x_1$.

\begin{figure}[htb]
\centering
\input{rejoinbad.pstex_t}
\caption{The region where four Voronoi cells $V_Z(z_1)$, $V_Z(z_2)$, $V_Z(x_1)$ and $V_Z(x_2)$ come together.
The cells associated to $z_1$ and $z_2$ meet in the line segment $O_1O_2$ with midpoint $M$. The thick
lines are parts of the boundaries of the Voronoi cells.}
\label{fig_rejoinbad}
\end{figure}
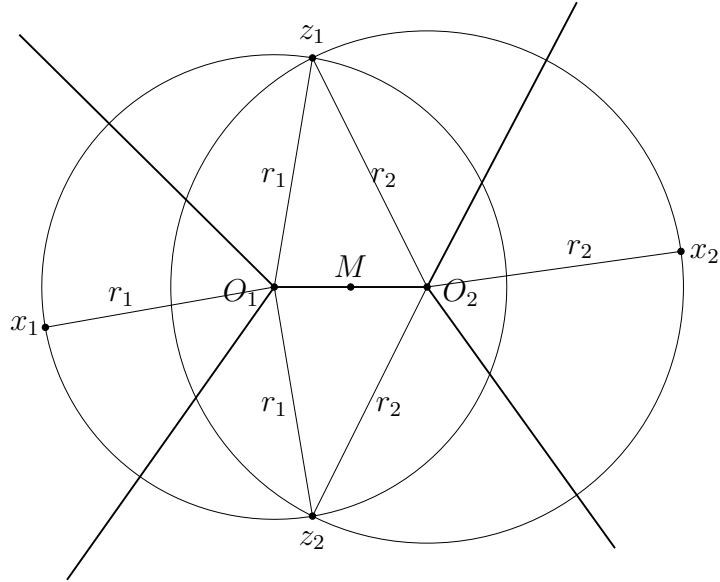

Suppose first that $\dist{O_1}{O_2} < \delta$. Let $r=\dist{z_1}{M}=\dist{z_2}{M}$. Using the triangle inequality twice,
$\left|\dist{x_1}{M}-\dist{z_1}{M}\right|
  \le \left|\dist{x_1}{O_1}-\dist{z_1}{O_1}\right|+2\dist{O_1}{M} =0 +\dist{O_1}{O_2} < \delta$.
By assumption, $\dist{z}{M} < r+\delta^6$.
As $M$ is in the Voronoi cells of both $z_1$ and $z_2$, these are the two closest
points of $Z$ to $M$, so no point of $Z$ is at distance less than $r$ from $M$.
As $\dist{z_1}{M},\dist{z_2}{M},\dist{x_1}{M},\dist{z}{M} < r+\delta$, 
it follows that $\{z_1,z_2,x_1,z\}$ is a $\delta$-bad quadruple.
(The assumption on point free discs ensures that $r<2\sqrt{\log s}$.)
We may thus suppose that $\dist{O_1}{O_2} \ge \delta$.
We may also assume that $\dist{z_1}{z_2}\ge \delta$; otherwise there is nothing to prove.

We shall now apply Lemma~\ref{l_geom}. There is no problem applying this lemma in the torus: by
our assumption on point free discs, $r_1, r_2\le 2\sqrt{\log s}$, so all points we consider are
within distance $10\sqrt{\log s}<s/3$ of $z_1$, say, and we need only consider a region of $\TT(s)$
isomorphic to the corresponding region of $\RR^2$. We apply Lemma~\ref{l_geom} to the circles
$C_i$ with centres $O_i$ and radii $r_i$, which meet at the points $z_1$, $z_2$. The conditions
do apply; we are assuming $\dist{z_1}{z_2}\ge \delta$ and $\dist{O_1}{O_2}\ge \delta$, while, as just noted,
$r_i\le 2\sqrt{\log s}\le 1/\delta$.

As $C_1$ and $C_2$ have no points of $Z$ in their interiors,
we conclude from Lemma~\ref{l_geom} that $z$ must be ($\delta/4$)-close
to $z_1$ or $z_2$. But then $z_1$, $z_2$ and $x_1$ are at distance exactly $r_1$ from $O_1$, while
$z$, which is distinct from these three points, can be at distance at most $r_1+\delta/4$. 
Hence, as $C_1$ is free of points of $Z$, the quadruple $\{z_1,z_2,x_1,z\}$ is $\delta$-bad as required.
\end{proof}

\subsection{Proof of Lemma~\ref{l_deter}}\label{a_l_deter}

\setcounter{theorem}{19}
\addtocounter{theorem}{-1}
\begin{lemma}
  Let $0<\eta<1/100$ be fixed.
  If $t$ is large enough, then it is impossible to
  arrange $M\ge \eta^{-6}t$ points $P_1,\ldots,P_M$ and $M$
  associated circles $C_1,\ldots,C_M$ in $\RR^2$ so that no two $P_i$
  are within distance $\eta$ of each other, the distance of every
  $P_i$ from the origin is between $\eta t$ and $t$, each $C_i$ has
  radius at most $t$, no $C_i$ contains any $P_j$, and $C_i$ passes
  within distance $\eta^3$ of both $P_i$ and the origin.
\end{lemma}

\begin{proof}[Proof of Lemma~\ref{l_deter}]
Consider any arrangement as described in the statement of the lemma. By assumption, all the $P_i$ lie in an annulus of inner
radius $\eta t$ and outer radius $t$ centred on the origin $O$. We may write this annulus
as a union of fewer than $\eta^{-2}$ annuli $A_i$, each of which has outer radius $(1+\eta)$
times its inner radius. It suffices to show that we cannot arrange at least $\eta^2 M$ points as
described in any one of the $A_i$. Suppose we can. Then there are two $P_i$, say $P_1$ and $P_2$, such that
the angle $P_1OP_2$ is at most $\theta_0=2\pi/(\eta^2 M)\le \eta^3/t$. Suppose that $\dist{O}{P_1}
\ge \dist{O}{P_2}$. Let us rotate the coordinates
so that $P_1$ is the point $(0,y)$, $y>0$,
and $P_2$ the point $(x,y')$, $y'<y$. As $\dist{O}{P_2}$ is at most $t$, we see
that $|x|\le t\theta_0\le \eta^3$. Using $\dist{P_1}{P_2}\ge \eta$, it follows that $|y-y'|\ge\eta/2$.

Now it is easy to see that $C_1$ must contain $P_2$: otherwise, translating through a distance
of at most $\eta^3$ so that the image $C_1'$ of $C_1$ passes
exactly through the origin, the image $P_2'$ of $P_2$ lies outside $C_1'$ and satisfies $\dist{P_2}{P_2'}\le \eta^3$.
Also, there is a point $P_1'$ (close to the image of $P_1$) which lies exactly on $C_1'$,
with $\dist{P_1}{P_1'}\le 2\eta^3$. Rotating again, we may write
$P_1'=(0,y_1)$, $P_2'=(x_2,y_2)$ with, crudely, $\eta t/2<y_2<y_1$, $y_1-y_2\ge \eta/3$ and $|x_2|\le 5\eta^3$.

\begin{figure}[htb]
\centering
\input{detfig.pstex_t}
\caption{The notation in the proof of Lemma~\ref{l_deter}.}\label{fig_detfig}
\end{figure}
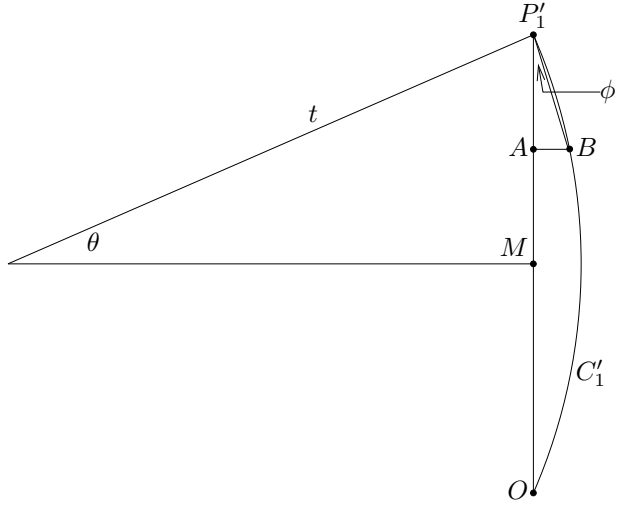

As the circle $C_1'$ passes through the origin and $P_1'=(0,y_1)$, it contains $A=(0,y_2)$ and
certain other points on the line $y=y_2$. The extreme case is when the circle has maximum
radius, $t$, by assumption. This case is shown in Figure~\ref{fig_detfig}. Let $B$ be the unique point $(x_3,y_2)$
on both $C_1'$ and the line $y=y_2$ for which $x_3>0$.
As $\dist{O}{P_1'}=y_1\ge \eta t/2$, writing $M$ for the midpoint of $OP_1'$ we have $\dist{M}{P_1'}\ge \eta t/4$,
so the angle $\theta$ shown is at least $\arcsin(\eta/4)$.
By assumption $\dist{O}{P_1}\le (1+\eta)\dist{O}{P_2}$, as $P_1$ and $P_2$ lie in the same
annulus $A_i$. It follows that, crudely,
$y_2>y_1/2$, so the angle $\phi=BP_1'A$ shown in Figure~\ref{fig_detfig} is at least $\theta/2$, the value
it takes when $y_2=y_1/2$, i.e., when $A=M$.
Thus
\[
 x_3=\dist{A}{B}=(y_1-y_2)\tan \phi \ge (\eta/3)\sin \phi \ge (\eta/3)\sin(\arcsin(\eta/4)/2) \ge \eta^2/24.
\]
As the case $r=t$ was extreme, in all cases
$C_1'$ contains all points $(x,y_2)$ with $|x|\le \eta^2/24$. In particular, $C_1'$ contains $P_2'$.
But then $C_1$ contains $P_2$, contradicting our assumptions.
\end{proof}

\subsection{Proof of Claim~\ref{cl_badsep}}\label{a_cl_badsep}

\setcounter{claim_base}{18}
\setcounter{claim}{0}
Claim~\ref{cl_badsep}, repeated below, is part of the proof of Lemma~\ref{l_bad}. 
The context is as follows: $\epsilon>0$ is fixed, and we have $\delta\le s^{-\epsilon}$ and $\delta'=\delta^{1/10}$.
The square $S$ is a square of side length $\ell=3(\log s)^2=o(1/\delta')$ in the torus $\TT(s)$.
Recall that a quadruple $\{x_1,x_2,x_3,x_4\}$ is {\em weakly $\delta$-bad}
if there is a point $x$ and a radius $r<2\sqrt{\log s}$ such that $r<\dist{x}{x_i}\le r+\delta$ holds for $1\le i\le 4$.
Note that this definition makes no reference to the Poisson process $Z$. Finally, all limits are as $s\to\infty$.

\begin{claim}
Suppose that $x_1$, $x_2$, and $x_3$ are points of $S$, no two within distance
$\delta'$ of each other. Let $B$ be the set of $x\in S$ for which the quadruple $\{x_1,x_2,x_3,x\}$ is
weakly $\delta$-bad. Then $\area(B)=o(\delta')$.
\end{claim}

\begin{proof}[Proof of Claim~\ref{cl_badsep}]
Although in the proof of Lemma~\ref{l_bad} we are working in the torus $\TT(s)$, here we shall work in $\RR^2$,
treating $S$ as a square in $\RR^2$ of side-length $\ell=o(1/\delta')$. We may do this because of the restriction
$r<2\sqrt{\log s}$ in the definition of weak $\delta$-badness, which ensures that no circle we need consider
comes close to `wrapping round' the torus.

To avoid degeneracy we shall assume that $x_1$, $x_2$ and $x_3$ are not collinear. (This is valid as we may shift one
slightly and increase $\delta$ slightly.)

We say that four points of $\RR^2$ are {\em concyclic} if there is a circle passing through them. Suppose
that $x\in B$, i.e., that the quadruple $\{x_1,x_2,x_3,x\}$ is weakly $\delta$-bad.
Then there are points $y_i\in \RR^2$ such that $y_1$, $y_2$, $y_3$ and $x$ are concyclic, with $\dist{x_i}{y_i}\le \delta$.
(Fixing the centre, adjust the radius of the circle witnessing weak $\delta$-badness by at most $\delta$ so that this
circle passes through $x$.)
As the statement of the claim depends only on the set $\{x_1,x_2,x_3\}$, we may renumber so that $y_1$, $y_2$, $y_3$ and $x$ appear
in this order around the circle on which they lie.

Consider the linear transformation $T:\RR^2\to \RR^2$ obtained as follows: first translate $\RR^2$ so that
the point $y_1$ is mapped to $x_1$. Then rotate and rescale around $x_1$ so that (the image of) $y_3$ is mapped to $x_3$.
As $\dist{x_i}{y_i}\le \delta$, while $\dist{x_1}{x_3}\ge \delta'$, the angle of rotation is $O(\delta/\delta')$, and
the scale factor is $1+O(\delta/\delta')$. It follows that $T$, which maps $x_1$ to $y_1$ and $x_3$ to $y_3$, moves
each point inside $S$, or indeed within distance $\ell$ of $S$, through a distance of at most $O(\ell \delta/\delta')=O({\delta'}^8)$.
Hence,
\begin{equation}\label{Ty2x2}
 \dist{T(y_2)}{x_2} \le \dist{T(y_2)}{y_2} + \dist{y_2}{x_2} = O({\delta'}^8) + O(\delta) = O({\delta'}^8).
\end{equation}

\begin{figure}[htb]
\centering
\input{cross.pstex_t}
\caption{The points $x_1,T(y_2),x_3,T(x)$ are concyclic, as are the points $x_1,x_2,x_3,z$.}
\label{fig_cross}
\end{figure}
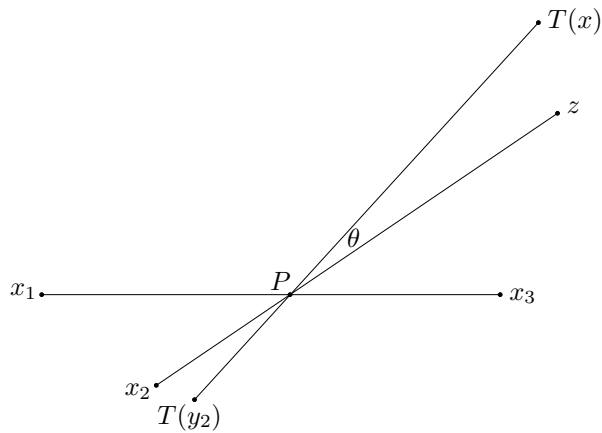

Applying $T$ to $\{y_1,y_2,y_3,x\}$ we see that
$x_1$, $T(y_2)$, $x_3$ and $T(x)$ are concyclic, and appear in this order around the circle on which they lie.
Let $P$ be the intersection of the line-segments $x_1x_3$ and $T(y_2)T(x)$, so $P$ lies in the interior of these
two line-segments; see Figure~\ref{fig_cross}. 
Then
\begin{equation}\label{cdist}
 \dist{x_1}{P}\ \dist{x_3}{P} = \dist{T(y_2)}{P}\ \dist{T(x)}{P}.
\end{equation}
We claim that $\dist{T(y_2)}{P}\ge {\delta'}^3$: otherwise, as $x$ and $T(y_2)$ lie in $S$,
$\dist{T(x)}{P}\le \dist{T(x)}{T(y_2)} \le 2\ell$, so from \eqref{cdist}, $\dist{x_1}{P}\ \dist{x_3}{P}\le 2\ell{\delta'}^3$.
Assuming without loss of generality that $\dist{x_1}{P}\le \dist{x_3}{P}$, as $\dist{x_1}{x_3}\ge \delta'$ by assumption,
we have $\dist{x_3}{P}\ge \delta'/2$. It follows that $\dist{x_1}{P}\le 4\ell{\delta'}^2$. But then using \eqref{Ty2x2},
if $s$ is large enough we have
\[
 \dist{x_1}{x_2} \le \dist{x_1}{P} + \dist{P}{T(y_2)} + \dist{T(y_2)}{x_2}
 \le 4\ell{\delta'}^2 + {\delta'}^3 + O({\delta'}^8) < \delta',
\]
contradicting one of our assumptions.

Let $C$ be the circle through $x_1$, $x_2$ and $x_3$, and let $z$ be the point on $C$ obtained by extending
the line-segment $x_2P$. Let $\theta$ be the angle $T(y_2)Px_2=T(x)Pz$, as shown in Figure~\ref{fig_cross}.
Now from \eqref{Ty2x2} and the fact that $\dist{T(y_2)}{P}\ge {\delta'}^3$, we see that $\theta=O({\delta'}^5)$.
Also, $\dist{x_2}{P} / \dist{T(y_2)}{P}  = 1+O({\delta'}^5)$. As $x_1$, $x_2$, $x_3$ and $z$ are concyclic,
\[
 \dist{z}{P} = \frac{ \dist{x_1}{P}\ \dist{x_3}{P} } { \dist{x_2}{P} }
 = \dist{T(x)}{P} \frac{ \dist{T(y_2)}{P} }{ \dist{x_2}{P} } = \dist{T(x)}{P}(1+O({\delta'}^5)) = \dist{T(x)}{P}+O(\ell {\delta'}^5).
\]
As $\dist{T(x)}{P}=O(\ell)$ and $\theta=O({\delta'}^5)$,
it follows that $\dist{T(x)}{z}=O(\ell {\delta'}^5)$. Thus, as $\dist{T(x)}{x}=O({\delta'}^8)$ from the properties of $T$,
we have $\dist{x}{z}=O(\ell {\delta'}^5)$.

In summary, we have shown that any point $x$ for which $\{x_1,x_2,x_3,x\}$ is weakly $\delta$-bad is
within distance $O(\ell {\delta'}^5)$ of some point $z$ on the circle through $\{x_1,x_2,x_3\}$. Thus, such an $x$
lies in a certain fixed annulus whose inner and outer radii differ by $O(\ell {\delta'}^5)$.
Any such annulus meets $S$ in a total area of $O(\ell^2 {\delta'}^5)$, completing the proof.
\end{proof}

\subsection{Proof of Claim~\ref{cl_pathincpt}}\label{a_cl_pathincpt}

Although purely deterministic, Claim~\ref{cl_pathincpt} is part of the proof of Lemma~\ref{l_coupling},
and the proof below should be read in this context.
We refer the reader to page~\pageref{cl_pathincpt} for the statement of Claim~\ref{cl_pathincpt},
and to page~\pageref{cases} for the three cases describing the properties
of the coupling constructed for a particular potentially bad component $C$.

\begin{proof}[Proof of Claim~\ref{cl_pathincpt}]
In the context of the proof of Lemma~\ref{l_coupling}, recall that we have constructed the coupling
of $(Z_1,\col_1)$ and $(Z_2,\col_2)$ in such a way that one of the
three cases 1,2,3 defined on page~\pageref{cases} holds. Recall also that we are assuming $B_1$ does not hold.
As $a,b\in C$ and $z_{1,a}$ and $z_{1,b}$ are black, case 3 does not hold, and we must consider only cases
1 and 2. Unfortunately, within each of these cases
we shall need to divide further into the two cases (i) that $z_{1,a}$ and $z_{1,b}$
are not $\delta_2$-close, and (ii) that they are. Let us write $M$ for the midpoint of the common edge of the Voronoi
cells of $z_{1,a}$ and $z_{1,b}$ (with respect to $Z_1$).

Suppose first that case 1(i) holds.
Then as $z_{1,a}$ and $z_{1,b}$ are not $\delta_2$-close, for any $c\notin C$ we have by Lemma~\ref{l_rejoinbad}
that $\dist{z_{1,c}}{M}\ge \dist{z_{1,a}}{M} + \delta_2^6$; otherwise $z_{1,c}$ would be in a
$\delta_2$-bad quadruple with $z_{1,a}$ and $z_{1,b}$. (The assumption of Lemma~\ref{l_rejoinbad} on 
point free discs is guaranteed by our present assumption
that $B_1$ does not hold; see page~\pageref{page_B1}.)
Hence, by definition of potential badness,
$c$ would be in a potentially bad quadruple with $a$
and $b$, contradicting $a\in C$, $c\notin C$.
Passing from $Z_1$ to $Z_2$, as we are in case 1 only points outside
$C$ may move (and then by at most $\sqrt{2}\delta_1<\delta_2^6$),
so $z_{2,a}$ and $z_{2,b}$ are the closest points of $Z_2$ to $M$, and by Lemma~\ref{l_closest}
their Voronoi cells meet (indeed, meet at $M$).
As $z_{2,a}$ and $z_{2,b}$ are black and $\delta$-good with respect to $Z_2$, this proves the claim in this case.

Now suppose that case 2(i) holds.
Here all points of $C$ are black and $\delta$-good with respect to $Z_2$, so we need only show that there
is a path joining the cells of $z_{2,a}$ and $z_{2,b}$ staying within cells associated to the component.
We claim that the path $P=z_{2,a}Mz_{2,b}$ will do: as above, all points $z_{1,c}$, $c\notin C$,
are at least $\delta_2^6\ge 2\sqrt{2}\delta_1$ further from $M$ than $z_{1,a}$ and $z_{1,b}$ are.
Passing to $Z_2$, it follows that $z_{2,a}$ and $z_{2,b}$ are both closer to $M$ than any point
$z_{2,c}$, $c\notin C$. Hence, as any point $x\in P$ can be reached from $M$ by moving in a straight line
towards $z_{2,a}$ or $z_{2,b}$, for any $x\in P$ one of $z_{2,a}$, $z_{2,b}$ is closer to $x$
than any $z_{2,c}$, $c\notin C$. Thus the closest point of $Z_2$ to $x$ is some $z_{2,d}$, $d\in C$,
as required.

The case 2(ii) is similar but simpler: now $z_{2,a}$ and $z_{2,b}$ are at distance at most
$\delta_2+2\sqrt{2}\delta_1$, and we can use the path $P=z_{2,a}z_{2,b}$: letting $M'$ be the midpoint
of this line-segment, it suffices to show that no $z_{2,c}$, $c\notin C$, is closer to $M'$ than
$z_{2,a}$, $z_{2,b}$ are. But if some $z_{2,c}$ were, it would be $\delta_2$-close to $z_{2,a}$ or $z_{2,b}$,
contradicting $a,b\in C$, $c\notin C$.

Finally, we have the case 1(ii). As in 1(i), it suffices to show that the Voronoi cells of $z_{2,a}$ and $z_{2,b}$
meet. They certainly do unless some $z_{2,c}$ is closer to $M'$, the midpoint of $z_{1,a}z_{1,b}=z_{2,a}z_{2,b}$,
than $z_{2,a}$ is, so we may suppose that there is such a $c$.
As in case 2(ii) above, for any such $c$ we must have $c\in C$. 

\begin{figure}[htb]
\centering
\input{pathincpt.pstex_t}
\caption{The points $z_i$ represent $z_{1,i}=z_{2,i}$.}
\label{fig_pathincpt}
\end{figure}
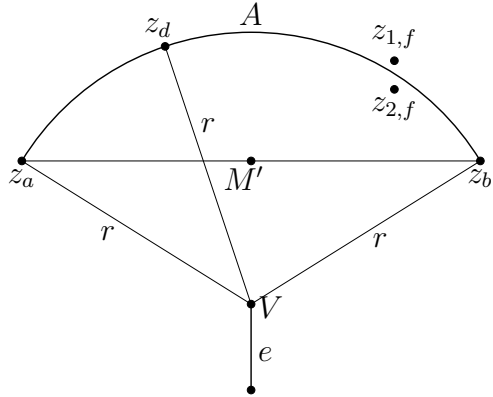

Consider the common edge $e$ of the Voronoi cells of $z_{1,a}$ and $z_{1,b}$ with respect to $Z_1$.
The location of $z_{1,c}=z_{2,c}$ ensures that $M'$ is not a point of this edge, so the edge lies
entirely to one side of $M'$. Let $V$ be the endpoint of $e$ closest to $M'$, and let $z_{1,d}$ be the third
point of $Z_1$ at the same distance $r$ from $V$ as $z_{1,a}$ and $z_{1,b}$. Then (as $V$ is the closer end
of $e$ to $M'$), $z_{1,d}$ lies on an arc-segment $A$ joining $z_{1,a}$ and $z_{1,b}$ every point of which
is within distance $\dist{z_{1,a}}{M'}$ of $M'$; see Figure~\ref{fig_pathincpt}.
Hence $d\in C$, as above. Passing to $Z_2$, the Voronoi
cells of $z_{2,a}=z_{1,a}$ and $z_{2,b}=z_{1,b}$ with respect to $Z_2$ meet at $V$ unless for some $f$
we have $\dist{z_{2,f}}{V}<\dist{z_{2,a}}{V}=\dist{z_{2,b}}{V}=r$.
Since only points outside $C$ move, and then only by a distance
at most $\sqrt{2}\delta_1$, we have $f\notin C$, so $a,b,d$ and $f$ are all distinct.
Furthermore, $\dist{z_{2,f}}{V}\ge r-\sqrt{2}\delta_1$ for any such $f$.
The quadruple $\{z_{2,a},z_{2,b},z_{2,d},z_{2,f}\}$ is thus $\delta_2$-bad, as witnessed by the annulus
with centre $V$, inner radius $r'=r-\sqrt{2}\delta_1$, and outer radius $r'+\delta_2>r$, in which
all four points lie. Note that the condition $r'<2\sqrt{\log s}$ required for $\delta_2$-badness holds
as our assumption that $B_1$ holds guarantees that $r<2\sqrt{\log s}$.
As $\{z_{2,a},z_{2,b},z_{2,d},z_{2,f}\}$ is $\delta_2$-bad, the quadruple $\{a,b,d,f\}$ is potentially bad,
contradicting $a\in C$, $f\notin C$ and completing the proof.
\end{proof}

\end{document}

%% file: 2closest.pstex_t
\begin{picture}(0,0)%
\includegraphics{2closest.pstex}%
\end{picture}%
\setlength{\unitlength}{1381sp}%
\begingroup\makeatletter\ifx\SetFigFont\undefined%
\gdef\SetFigFont#1#2#3#4#5{%
  \reset@font\fontsize{#1}{#2pt}%
  \fontfamily{#3}\fontseries{#4}\fontshape{#5}%
  \selectfont}%
\fi\endgroup%
\begin{picture}(5382,5382)(2410,-6952)
\put(5701,-4411){\makebox(0,0)[b]{\smash{{\SetFigFont{10}{12.0}{\rmdefault}{\mddefault}{\updefault}{\color[rgb]{0,0,0}$w$}%
}}}}
\put(5101,-4111){\makebox(0,0)[b]{\smash{{\SetFigFont{10}{12.0}{\rmdefault}{\mddefault}{\updefault}{\color[rgb]{0,0,0}$x$}%
}}}}
\put(3751,-5461){\makebox(0,0)[rb]{\smash{{\SetFigFont{10}{12.0}{\rmdefault}{\mddefault}{\updefault}{\color[rgb]{0,0,0}$z_1$}%
}}}}
\put(7651,-5461){\makebox(0,0)[lb]{\smash{{\SetFigFont{10}{12.0}{\rmdefault}{\mddefault}{\updefault}{\color[rgb]{0,0,0}$z_2$}%
}}}}
\end{picture}%

%% file: claimfig1.pstex_t
\begin{picture}(0,0)%
\includegraphics{claimfig1.pstex}%
\end{picture}%
\setlength{\unitlength}{1381sp}%
\begingroup\makeatletter\ifx\SetFigFont\undefined%
\gdef\SetFigFont#1#2#3#4#5{%
  \reset@font\fontsize{#1}{#2pt}%
  \fontfamily{#3}\fontseries{#4}\fontshape{#5}%
  \selectfont}%
\fi\endgroup%
\begin{picture}(9930,9073)(2086,-8483)
\put(2101,-4861){\makebox(0,0)[rb]{\smash{{\SetFigFont{8}{9.6}{\familydefault}{\mddefault}{\updefault}{\color[rgb]{0,0,0}$-\varepsilon s$}%
}}}}
\put(2101,-2461){\makebox(0,0)[rb]{\smash{{\SetFigFont{8}{9.6}{\familydefault}{\mddefault}{\updefault}{\color[rgb]{0,0,0}$3\varepsilon s$}%
}}}}
\put(2401,-4261){\makebox(0,0)[lb]{\smash{{\SetFigFont{8}{9.6}{\familydefault}{\mddefault}{\updefault}{\color[rgb]{0,0,0}$L$}%
}}}}
\put(2401,-3061){\makebox(0,0)[lb]{\smash{{\SetFigFont{8}{9.6}{\familydefault}{\mddefault}{\updefault}{\color[rgb]{0,0,0}$L'$}%
}}}}
\put(2101,-3661){\makebox(0,0)[rb]{\smash{{\SetFigFont{8}{9.6}{\familydefault}{\mddefault}{\updefault}{\color[rgb]{0,0,0}$\varepsilon s$}%
}}}}
\put(12001,-4861){\makebox(0,0)[lb]{\smash{{\SetFigFont{8}{9.6}{\familydefault}{\mddefault}{\updefault}{\color[rgb]{0,0,0}$S'$}%
}}}}
\put(2101,-7561){\makebox(0,0)[rb]{\smash{{\SetFigFont{8}{9.6}{\familydefault}{\mddefault}{\updefault}{\color[rgb]{0,0,0}$-s/2$}%
}}}}
\put(2101,-8161){\makebox(0,0)[rb]{\smash{{\SetFigFont{8}{9.6}{\familydefault}{\mddefault}{\updefault}{\color[rgb]{0,0,0}$-s/2-2\varepsilon s$}%
}}}}
\put(11401,-3961){\makebox(0,0)[lb]{\smash{{\SetFigFont{8}{9.6}{\familydefault}{\mddefault}{\updefault}{\color[rgb]{0,0,0}$R$}%
}}}}
\put(2101,239){\makebox(0,0)[rb]{\smash{{\SetFigFont{8}{9.6}{\familydefault}{\mddefault}{\updefault}{\color[rgb]{0,0,0}$s/2+2\varepsilon s$}%
}}}}
\put(8401,-4561){\makebox(0,0)[lb]{\smash{{\SetFigFont{8}{9.6}{\familydefault}{\mddefault}{\updefault}{\color[rgb]{0,0,0}$P_2$}%
}}}}
\put(3601,-4861){\makebox(0,0)[lb]{\smash{{\SetFigFont{8}{9.6}{\familydefault}{\mddefault}{\updefault}{\color[rgb]{0,0,0}$P_1$}%
}}}}
\end{picture}%

%% file: newsplit2.pstex_t
\begin{picture}(0,0)%
\includegraphics{newsplit2.pstex}%
\end{picture}%
\setlength{\unitlength}{1500sp}%
\begingroup\makeatletter\ifx\SetFigFont\undefined%
\gdef\SetFigFont#1#2#3#4#5{%
  \reset@font\fontsize{#1}{#2pt}%
  \fontfamily{#3}\fontseries{#4}\fontshape{#5}%
  \selectfont}%
\fi\endgroup%
\begin{picture}(8955,9183)(886,-9268)
\put(1201,-9136){\makebox(0,0)[b]{\smash{{\SetFigFont{10}{12.0}{\rmdefault}{\mddefault}{\updefault}{\color[rgb]{0,0,0}$0$}%
}}}}
\put(1051,-436){\makebox(0,0)[rb]{\smash{{\SetFigFont{10}{12.0}{\rmdefault}{\mddefault}{\updefault}{\color[rgb]{0,0,0}$s/2$}%
}}}}
\put(901,-4636){\makebox(0,0)[rb]{\smash{{\SetFigFont{10}{12.0}{\rmdefault}{\mddefault}{\updefault}{\color[rgb]{0,0,0}$0$}%
}}}}
\put(9826,-4636){\makebox(0,0)[lb]{\smash{{\SetFigFont{10}{12.0}{\rmdefault}{\mddefault}{\updefault}{\color[rgb]{0,0,0}$0$}%
}}}}
\put(6076,-3286){\makebox(0,0)[b]{\smash{{\SetFigFont{10}{12.0}{\rmdefault}{\mddefault}{\updefault}{\color[rgb]{0,0,0}$P_1$}%
}}}}
\put(4951,-6511){\makebox(0,0)[b]{\smash{{\SetFigFont{10}{12.0}{\rmdefault}{\mddefault}{\updefault}{\color[rgb]{0,0,0}$P_0$}%
}}}}
\put(5701,-4936){\makebox(0,0)[b]{\smash{{\SetFigFont{10}{12.0}{\rmdefault}{\mddefault}{\updefault}{\color[rgb]{0,0,0}$P$}%
}}}}
\put(9601,-9136){\makebox(0,0)[b]{\smash{{\SetFigFont{10}{12.0}{\rmdefault}{\mddefault}{\updefault}{\color[rgb]{0,0,0}$s$}%
}}}}
\put(1501,-9136){\makebox(0,0)[lb]{\smash{{\SetFigFont{10}{12.0}{\rmdefault}{\mddefault}{\updefault}{\color[rgb]{0,0,0}$0.01s$}%
}}}}
\put(9301,-9136){\makebox(0,0)[rb]{\smash{{\SetFigFont{10}{12.0}{\rmdefault}{\mddefault}{\updefault}{\color[rgb]{0,0,0}$0.99s$}%
}}}}
\put(1051,-886){\makebox(0,0)[rb]{\smash{{\SetFigFont{10}{12.0}{\rmdefault}{\mddefault}{\updefault}{\color[rgb]{0,0,0}$0.497s$}%
}}}}
\put(976,-8536){\makebox(0,0)[rb]{\smash{{\SetFigFont{10}{12.0}{\rmdefault}{\mddefault}{\updefault}{\color[rgb]{0,0,0}$-0.497s$}%
}}}}
\put(976,-8911){\makebox(0,0)[rb]{\smash{{\SetFigFont{10}{12.0}{\rmdefault}{\mddefault}{\updefault}{\color[rgb]{0,0,0}$-s/2$}%
}}}}
\end{picture}%

%% file: resplit.pstex_t
\begin{picture}(0,0)%
\includegraphics{resplit.pstex}%
\end{picture}%
\setlength{\unitlength}{1500sp}%
\begingroup\makeatletter\ifx\SetFigFont\undefined%
\gdef\SetFigFont#1#2#3#4#5{%
  \reset@font\fontsize{#1}{#2pt}%
  \fontfamily{#3}\fontseries{#4}\fontshape{#5}%
  \selectfont}%
\fi\endgroup%
\begin{picture}(7980,8034)(2611,-8143)
\put(8476,-3286){\makebox(0,0)[b]{\smash{{\SetFigFont{10}{12.0}{\rmdefault}{\mddefault}{\updefault}{\color[rgb]{0,0,0}$P_1$}%
}}}}
\put(3451,-436){\makebox(0,0)[rb]{\smash{{\SetFigFont{10}{12.0}{\rmdefault}{\mddefault}{\updefault}{\color[rgb]{0,0,0}$s$}%
}}}}
\put(7051,-6811){\makebox(0,0)[lb]{\smash{{\SetFigFont{10}{12.0}{\rmdefault}{\mddefault}{\updefault}{\color[rgb]{0,0,0}$0.49s$}%
}}}}
\put(6901,-6811){\makebox(0,0)[rb]{\smash{{\SetFigFont{10}{12.0}{\rmdefault}{\mddefault}{\updefault}{\color[rgb]{0,0,0}$0.47s$}%
}}}}
\put(3601,-8011){\makebox(0,0)[b]{\smash{{\SetFigFont{10}{12.0}{\rmdefault}{\mddefault}{\updefault}{\color[rgb]{0,0,0}$0$}%
}}}}
\put(10201,-8011){\makebox(0,0)[b]{\smash{{\SetFigFont{10}{12.0}{\rmdefault}{\mddefault}{\updefault}{\color[rgb]{0,0,0}$0.96s$}%
}}}}
\put(5101,-1411){\makebox(0,0)[b]{\smash{{\SetFigFont{10}{12.0}{\rmdefault}{\mddefault}{\updefault}{\color[rgb]{0,0,0}$R_0$}%
}}}}
\put(8701,-1411){\makebox(0,0)[b]{\smash{{\SetFigFont{10}{12.0}{\rmdefault}{\mddefault}{\updefault}{\color[rgb]{0,0,0}$R_1$}%
}}}}
\put(4951,-2686){\makebox(0,0)[b]{\smash{{\SetFigFont{10}{12.0}{\rmdefault}{\mddefault}{\updefault}{\color[rgb]{0,0,0}$P_0$}%
}}}}
\put(5476,-5386){\makebox(0,0)[b]{\smash{{\SetFigFont{10}{12.0}{\rmdefault}{\mddefault}{\updefault}{\color[rgb]{0,0,0}$P$}%
}}}}
\put(3451,-1636){\makebox(0,0)[rb]{\smash{{\SetFigFont{10}{12.0}{\rmdefault}{\mddefault}{\updefault}{\color[rgb]{0,0,0}$y+0.47s$}%
}}}}
\put(3301,-4036){\makebox(0,0)[rb]{\smash{{\SetFigFont{10}{12.0}{\rmdefault}{\mddefault}{\updefault}{\color[rgb]{0,0,0}$y$}%
}}}}
\put(3451,-6436){\makebox(0,0)[rb]{\smash{{\SetFigFont{10}{12.0}{\rmdefault}{\mddefault}{\updefault}{\color[rgb]{0,0,0}$y-0.47s$}%
}}}}
\put(10576,-4036){\makebox(0,0)[lb]{\smash{{\SetFigFont{10}{12.0}{\rmdefault}{\mddefault}{\updefault}{\color[rgb]{0,0,0}$y_1$}%
}}}}
\put(2626,-4486){\makebox(0,0)[rb]{\smash{{\SetFigFont{10}{12.0}{\rmdefault}{\mddefault}{\updefault}{\color[rgb]{0,0,0}$y_0$}%
}}}}
\end{picture}%

%% file: baddef.pstex_t
\begin{picture}(0,0)%
\includegraphics{baddef.pstex}%
\end{picture}%
\setlength{\unitlength}{1579sp}%
\begingroup\makeatletter\ifx\SetFigFont\undefined%
\gdef\SetFigFont#1#2#3#4#5{%
  \reset@font\fontsize{#1}{#2pt}%
  \fontfamily{#3}\fontseries{#4}\fontshape{#5}%
  \selectfont}%
\fi\endgroup%
\begin{picture}(4795,5425)(2386,-5897)
\put(5401,-2461){\makebox(0,0)[rb]{\smash{{\SetFigFont{10}{12.0}{\rmdefault}{\mddefault}{\updefault}{\color[rgb]{0,0,0}$r+\delta$}%
}}}}
\put(6001,-3211){\makebox(0,0)[b]{\smash{{\SetFigFont{10}{12.0}{\rmdefault}{\mddefault}{\updefault}{\color[rgb]{0,0,0}$r$}%
}}}}
\put(2401,-3061){\makebox(0,0)[rb]{\smash{{\SetFigFont{10}{12.0}{\rmdefault}{\mddefault}{\updefault}{\color[rgb]{0,0,0}$z_3$}%
}}}}
\put(3901,-5761){\makebox(0,0)[b]{\smash{{\SetFigFont{10}{12.0}{\rmdefault}{\mddefault}{\updefault}{\color[rgb]{0,0,0}$z_2$}%
}}}}
\put(4801,-811){\makebox(0,0)[b]{\smash{{\SetFigFont{10}{12.0}{\rmdefault}{\mddefault}{\updefault}{\color[rgb]{0,0,0}$z_1$}%
}}}}
\put(6976,-4486){\makebox(0,0)[lb]{\smash{{\SetFigFont{10}{12.0}{\rmdefault}{\mddefault}{\updefault}{\color[rgb]{0,0,0}$z_4$}%
}}}}
\put(4801,-3661){\makebox(0,0)[b]{\smash{{\SetFigFont{10}{12.0}{\rmdefault}{\mddefault}{\updefault}{\color[rgb]{0,0,0}$x$}%
}}}}
\end{picture}%

%% file: badeg.pstex_t
\begin{picture}(0,0)%
\includegraphics{badeg.pstex}%
\end{picture}%
\setlength{\unitlength}{1776sp}%
\begingroup\makeatletter\ifx\SetFigFont\undefined%
\gdef\SetFigFont#1#2#3#4#5{%
  \reset@font\fontsize{#1}{#2pt}%
  \fontfamily{#3}\fontseries{#4}\fontshape{#5}%
  \selectfont}%
\fi\endgroup%
\begin{picture}(6016,6511)(1793,-6368)
\put(4801,-3661){\makebox(0,0)[b]{\smash{{\SetFigFont{11}{13.2}{\rmdefault}{\mddefault}{\updefault}{\color[rgb]{0,0,0}$x$}%
}}}}
\put(4771,-196){\makebox(0,0)[b]{\smash{{\SetFigFont{11}{13.2}{\rmdefault}{\mddefault}{\updefault}{\color[rgb]{0,0,0}$z_1$}%
}}}}
\put(7771,-2161){\makebox(0,0)[b]{\smash{{\SetFigFont{11}{13.2}{\rmdefault}{\mddefault}{\updefault}{\color[rgb]{0,0,0}$C$}%
}}}}
\put(5791,-2491){\makebox(0,0)[b]{\smash{{\SetFigFont{11}{13.2}{\rmdefault}{\mddefault}{\updefault}{\color[rgb]{0,0,0}$C_1$}%
}}}}
\end{picture}%

%% file: couple.pstex_t
\begin{picture}(0,0)%
\includegraphics{couple.pstex}%
\end{picture}%
\setlength{\unitlength}{1500sp}%
\begingroup\makeatletter\ifx\SetFigFont\undefined%
\gdef\SetFigFont#1#2#3#4#5{%
  \reset@font\fontsize{#1}{#2pt}%
  \fontfamily{#3}\fontseries{#4}\fontshape{#5}%
  \selectfont}%
\fi\endgroup%
\begin{picture}(6000,6000)(2401,-7561)
\put(3526,-5461){\makebox(0,0)[lb]{\smash{{\SetFigFont{9}{10.8}{\rmdefault}{\mddefault}{\updefault}{\color[rgb]{0,0,0}$z$}%
}}}}
\put(2551,-4336){\makebox(0,0)[lb]{\smash{{\SetFigFont{9}{10.8}{\rmdefault}{\mddefault}{\updefault}{\color[rgb]{0,0,0}$z_2$}%
}}}}
\put(3451,-2986){\makebox(0,0)[b]{\smash{{\SetFigFont{9}{10.8}{\rmdefault}{\mddefault}{\updefault}{\color[rgb]{0,0,0}$z_3$}%
}}}}
\put(3226,-5611){\makebox(0,0)[rb]{\smash{{\SetFigFont{9}{10.8}{\rmdefault}{\mddefault}{\updefault}{\color[rgb]{0,0,0}$\phi(z)$}%
}}}}
\put(7726,-4111){\makebox(0,0)[lb]{\smash{{\SetFigFont{9}{10.8}{\rmdefault}{\mddefault}{\updefault}{\color[rgb]{0,0,0}$\phi(z')$}%
}}}}
\put(7426,-4261){\makebox(0,0)[rb]{\smash{{\SetFigFont{9}{10.8}{\rmdefault}{\mddefault}{\updefault}{\color[rgb]{0,0,0}$z'$}%
}}}}
\put(4201,-2086){\makebox(0,0)[lb]{\smash{{\SetFigFont{9}{10.8}{\rmdefault}{\mddefault}{\updefault}{\color[rgb]{0,0,0}$z_4$}%
}}}}
\end{picture}%

%% file: RP12.pstex_t
\begin{picture}(0,0)%
\includegraphics{RP12.pstex}%
\end{picture}%
\setlength{\unitlength}{1579sp}%
\begingroup\makeatletter\ifx\SetFigFont\undefined%
\gdef\SetFigFont#1#2#3#4#5{%
  \reset@font\fontsize{#1}{#2pt}%
  \fontfamily{#3}\fontseries{#4}\fontshape{#5}%
  \selectfont}%
\fi\endgroup%
\begin{picture}(10623,2911)(211,-4883)
\put(226,-4261){\makebox(0,0)[rb]{\smash{{\SetFigFont{10}{12.0}{\rmdefault}{\mddefault}{\updefault}{\color[rgb]{0,0,0}$R_1$}%
}}}}
\put(4501,-2311){\makebox(0,0)[rb]{\smash{{\SetFigFont{10}{12.0}{\rmdefault}{\mddefault}{\updefault}{\color[rgb]{0,0,0}$R_2$}%
}}}}
\end{picture}%

%% file: fig1.pstex_t
\begin{picture}(0,0)%
\includegraphics{fig1.pstex}%
\end{picture}%
\setlength{\unitlength}{1460sp}%
\begingroup\makeatletter\ifx\SetFigFont\undefined%
\gdef\SetFigFont#1#2#3#4#5{%
  \reset@font\fontsize{#1}{#2pt}%
  \fontfamily{#3}\fontseries{#4}\fontshape{#5}%
  \selectfont}%
\fi\endgroup%
\begin{picture}(13224,8066)(-11,-7994)
\put(9376,-511){\makebox(0,0)[b]{\smash{\SetFigFont{11}{13.2}{\rmdefault}{\mddefault}{\updefault}{\color[rgb]{0,0,0}$C_2$}%
}}}
\put(7276,-4936){\makebox(0,0)[b]{\smash{\SetFigFont{11}{13.2}{\rmdefault}{\mddefault}{\updefault}{\color[rgb]{0,0,0}$r_2$}%
}}}
\put(4726,-7336){\makebox(0,0)[b]{\smash{\SetFigFont{11}{13.2}{\rmdefault}{\mddefault}{\updefault}{\color[rgb]{0,0,0}$C_4$}%
}}}
\put(2926,-3811){\makebox(0,0)[b]{\smash{\SetFigFont{11}{13.2}{\rmdefault}{\mddefault}{\updefault}{\color[rgb]{0,0,0}$\theta$}%
}}}
\put(5326,-3811){\makebox(0,0)[b]{\smash{\SetFigFont{11}{13.2}{\rmdefault}{\mddefault}{\updefault}{\color[rgb]{0,0,0}$\psi$}%
}}}
\put(5176,-2236){\makebox(0,0)[b]{\smash{\SetFigFont{11}{13.2}{\rmdefault}{\mddefault}{\updefault}{\color[rgb]{0,0,0}$P_1$}%
}}}
\put(8401,-3736){\makebox(0,0)[b]{\smash{\SetFigFont{11}{13.2}{\rmdefault}{\mddefault}{\updefault}{\color[rgb]{0,0,0}$O_2$}%
}}}
\put(4726,-5761){\makebox(0,0)[rb]{\smash{\SetFigFont{11}{13.2}{\rmdefault}{\mddefault}{\updefault}{\color[rgb]{0,0,0}$P_2$}%
}}}
\put(4276,-2311){\makebox(0,0)[b]{\smash{\SetFigFont{11}{13.2}{\rmdefault}{\mddefault}{\updefault}{\color[rgb]{0,0,0}$\phi$}%
}}}
\put(3301,-2911){\makebox(0,0)[b]{\smash{\SetFigFont{11}{13.2}{\rmdefault}{\mddefault}{\updefault}{\color[rgb]{0,0,0}$r_1$}%
}}}
\put(5926,-4336){\makebox(0,0)[b]{\smash{\SetFigFont{11}{13.2}{\rmdefault}{\mddefault}{\updefault}{\color[rgb]{0,0,0}$M$}%
}}}
\put(3001,-4336){\makebox(0,0)[b]{\smash{\SetFigFont{11}{13.2}{\rmdefault}{\mddefault}{\updefault}{\color[rgb]{0,0,0}$O_1$}%
}}}
\put(1351,-1786){\makebox(0,0)[b]{\smash{\SetFigFont{11}{13.2}{\rmdefault}{\mddefault}{\updefault}{\color[rgb]{0,0,0}$C_1$}%
}}}
\put(3676,-1336){\makebox(0,0)[b]{\smash{\SetFigFont{11}{13.2}{\rmdefault}{\mddefault}{\updefault}{\color[rgb]{0,0,0}$P$}%
}}}
\put(4426,-1561){\makebox(0,0)[b]{\smash{\SetFigFont{11}{13.2}{\rmdefault}{\mddefault}{\updefault}{\color[rgb]{0,0,0}$P'$}%
}}}
\put(4726,-736){\makebox(0,0)[b]{\smash{\SetFigFont{11}{13.2}{\rmdefault}{\mddefault}{\updefault}{\color[rgb]{0,0,0}$C_3$}%
}}}
\end{picture}

%% file: rejoinbad.pstex_t
\begin{picture}(0,0)%
\includegraphics{rejoinbad.pstex}%
\end{picture}%
\setlength{\unitlength}{1579sp}%
\begingroup\makeatletter\ifx\SetFigFont\undefined%
\gdef\SetFigFont#1#2#3#4#5{%
  \reset@font\fontsize{#1}{#2pt}%
  \fontfamily{#3}\fontseries{#4}\fontshape{#5}%
  \selectfont}%
\fi\endgroup%
\begin{picture}(10577,9142)(2564,-9195)
\put(7201,-8611){\makebox(0,0)[b]{\smash{{\SetFigFont{12}{14.4}{\familydefault}{\mddefault}{\updefault}{\color[rgb]{0,0,0}$z_2$}%
}}}}
\put(7801,-4411){\makebox(0,0)[b]{\smash{{\SetFigFont{12}{14.4}{\familydefault}{\mddefault}{\updefault}{\color[rgb]{0,0,0}$M$}%
}}}}
\put(8401,-6511){\makebox(0,0)[b]{\smash{{\SetFigFont{12}{14.4}{\familydefault}{\mddefault}{\updefault}{\color[rgb]{0,0,0}$r_2$}%
}}}}
\put(6601,-2911){\makebox(0,0)[b]{\smash{{\SetFigFont{12}{14.4}{\familydefault}{\mddefault}{\updefault}{\color[rgb]{0,0,0}$r_1$}%
}}}}
\put(6601,-6511){\makebox(0,0)[b]{\smash{{\SetFigFont{12}{14.4}{\familydefault}{\mddefault}{\updefault}{\color[rgb]{0,0,0}$r_1$}%
}}}}
\put(7201,-661){\makebox(0,0)[b]{\smash{{\SetFigFont{12}{14.4}{\familydefault}{\mddefault}{\updefault}{\color[rgb]{0,0,0}$z_1$}%
}}}}
\put(8326,-2911){\makebox(0,0)[b]{\smash{{\SetFigFont{12}{14.4}{\familydefault}{\mddefault}{\updefault}{\color[rgb]{0,0,0}$r_2$}%
}}}}
\put(6076,-4786){\makebox(0,0)[b]{\smash{{\SetFigFont{12}{14.4}{\familydefault}{\mddefault}{\updefault}{\color[rgb]{0,0,0}$O_1$}%
}}}}
\put(9526,-4786){\makebox(0,0)[b]{\smash{{\SetFigFont{12}{14.4}{\familydefault}{\mddefault}{\updefault}{\color[rgb]{0,0,0}$O_2$}%
}}}}
\put(13126,-4111){\makebox(0,0)[lb]{\smash{{\SetFigFont{12}{14.4}{\familydefault}{\mddefault}{\updefault}{\color[rgb]{0,0,0}$x_2$}%
}}}}
\put(4201,-4786){\makebox(0,0)[b]{\smash{{\SetFigFont{12}{14.4}{\familydefault}{\mddefault}{\updefault}{\color[rgb]{0,0,0}$r_1$}%
}}}}
\put(11401,-4036){\makebox(0,0)[b]{\smash{{\SetFigFont{12}{14.4}{\familydefault}{\mddefault}{\updefault}{\color[rgb]{0,0,0}$r_2$}%
}}}}
\put(2926,-5236){\makebox(0,0)[rb]{\smash{{\SetFigFont{12}{14.4}{\familydefault}{\mddefault}{\updefault}{\color[rgb]{0,0,0}$x_1$}%
}}}}
\end{picture}%

%% file: detfig.pstex_t
\begin{picture}(0,0)%
\includegraphics{detfig.pstex}%
\end{picture}%
\setlength{\unitlength}{1579sp}%
\begingroup\makeatletter\ifx\SetFigFont\undefined%
\gdef\SetFigFont#1#2#3#4#5{%
  \reset@font\fontsize{#1}{#2pt}%
  \fontfamily{#3}\fontseries{#4}\fontshape{#5}%
  \selectfont}%
\fi\endgroup%
\begin{picture}(9344,7860)(1179,-7636)
\put(10126,-5761){\makebox(0,0)[lb]{\smash{\SetFigFont{10}{12.0}{\rmdefault}{\mddefault}{\updefault}{\color[rgb]{0,0,0}$C_1'$}%
}}}
\put(9376,-2236){\makebox(0,0)[rb]{\smash{\SetFigFont{10}{12.0}{\rmdefault}{\mddefault}{\updefault}{\color[rgb]{0,0,0}$A$}%
}}}
\put(10126,-2236){\makebox(0,0)[lb]{\smash{\SetFigFont{10}{12.0}{\rmdefault}{\mddefault}{\updefault}{\color[rgb]{0,0,0}$B$}%
}}}
\put(9376,-7636){\makebox(0,0)[rb]{\smash{\SetFigFont{10}{12.0}{\rmdefault}{\mddefault}{\updefault}{\color[rgb]{0,0,0}$O$}%
}}}
\put(2551,-3736){\makebox(0,0)[b]{\smash{\SetFigFont{10}{12.0}{\rmdefault}{\mddefault}{\updefault}{\color[rgb]{0,0,0}$\theta$}%
}}}
\put(5926,-1711){\makebox(0,0)[lb]{\smash{\SetFigFont{10}{12.0}{\rmdefault}{\mddefault}{\updefault}{\color[rgb]{0,0,0}$t$}%
}}}
\put(10501,-1336){\makebox(0,0)[lb]{\smash{\SetFigFont{10}{12.0}{\rmdefault}{\mddefault}{\updefault}{\color[rgb]{0,0,0}$\phi$}%
}}}
\put(9376,-3811){\makebox(0,0)[rb]{\smash{\SetFigFont{10}{12.0}{\rmdefault}{\mddefault}{\updefault}{\color[rgb]{0,0,0}$M$}%
}}}
\put(9451,-136){\makebox(0,0)[b]{\smash{\SetFigFont{10}{12.0}{\rmdefault}{\mddefault}{\updefault}{\color[rgb]{0,0,0}$P_1'$}%
}}}
\end{picture}

%% file: cross.pstex_t
\begin{picture}(0,0)%
\includegraphics{cross.pstex}%
\end{picture}%
\setlength{\unitlength}{1579sp}%
\begingroup\makeatletter\ifx\SetFigFont\undefined%
\gdef\SetFigFont#1#2#3#4#5{%
  \reset@font\fontsize{#1}{#2pt}%
  \fontfamily{#3}\fontseries{#4}\fontshape{#5}%
  \selectfont}%
\fi\endgroup%
\begin{picture}(8355,6700)(2011,-6722)
\put(3826,-6136){\makebox(0,0)[rb]{\smash{{\SetFigFont{10}{12.0}{\rmdefault}{\mddefault}{\updefault}{\color[rgb]{0,0,0}$x_2$}%
}}}}
\put(5851,-4486){\makebox(0,0)[b]{\smash{{\SetFigFont{10}{12.0}{\rmdefault}{\mddefault}{\updefault}{\color[rgb]{0,0,0}$P$}%
}}}}
\put(2026,-4561){\makebox(0,0)[rb]{\smash{{\SetFigFont{10}{12.0}{\rmdefault}{\mddefault}{\updefault}{\color[rgb]{0,0,0}$x_1$}%
}}}}
\put(4426,-6586){\makebox(0,0)[b]{\smash{{\SetFigFont{10}{12.0}{\rmdefault}{\mddefault}{\updefault}{\color[rgb]{0,0,0}$T(y_2)$}%
}}}}
\put(9451,-4636){\makebox(0,0)[lb]{\smash{{\SetFigFont{10}{12.0}{\rmdefault}{\mddefault}{\updefault}{\color[rgb]{0,0,0}$x_3$}%
}}}}
\put(10351,-1711){\makebox(0,0)[lb]{\smash{{\SetFigFont{10}{12.0}{\rmdefault}{\mddefault}{\updefault}{\color[rgb]{0,0,0}$z$}%
}}}}
\put(10051,-361){\makebox(0,0)[lb]{\smash{{\SetFigFont{10}{12.0}{\rmdefault}{\mddefault}{\updefault}{\color[rgb]{0,0,0}$T(x)$}%
}}}}
\put(6901,-3811){\makebox(0,0)[lb]{\smash{{\SetFigFont{10}{12.0}{\rmdefault}{\mddefault}{\updefault}{\color[rgb]{0,0,0}$\theta$}%
}}}}
\end{picture}%

%% file: pathincpt.pstex_t
\begin{picture}(0,0)%
\includegraphics{pathincpt.pstex}%
\end{picture}%
\setlength{\unitlength}{2368sp}%
\begingroup\makeatletter\ifx\SetFigFont\undefined%
\gdef\SetFigFont#1#2#3#4#5{%
  \reset@font\fontsize{#1}{#2pt}%
  \fontfamily{#3}\fontseries{#4}\fontshape{#5}%
  \selectfont}%
\fi\endgroup%
\begin{picture}(4896,4115)(2353,-5508)
\put(2401,-3286){\makebox(0,0)[b]{\smash{{\SetFigFont{12}{14.4}{\rmdefault}{\mddefault}{\updefault}{\color[rgb]{0,0,0}$z_a$}%
}}}}
\put(4726,-3361){\makebox(0,0)[b]{\smash{{\SetFigFont{12}{14.4}{\rmdefault}{\mddefault}{\updefault}{\color[rgb]{0,0,0}$M'$}%
}}}}
\put(3826,-1711){\makebox(0,0)[b]{\smash{{\SetFigFont{12}{14.4}{\rmdefault}{\mddefault}{\updefault}{\color[rgb]{0,0,0}$z_d$}%
}}}}
\put(7201,-3286){\makebox(0,0)[b]{\smash{{\SetFigFont{12}{14.4}{\rmdefault}{\mddefault}{\updefault}{\color[rgb]{0,0,0}$z_b$}%
}}}}
\put(4876,-5161){\makebox(0,0)[lb]{\smash{{\SetFigFont{12}{14.4}{\rmdefault}{\mddefault}{\updefault}{\color[rgb]{0,0,0}$e$}%
}}}}
\put(4876,-4711){\makebox(0,0)[lb]{\smash{{\SetFigFont{12}{14.4}{\rmdefault}{\mddefault}{\updefault}{\color[rgb]{0,0,0}$V$}%
}}}}
\put(6301,-1786){\makebox(0,0)[b]{\smash{{\SetFigFont{12}{14.4}{\rmdefault}{\mddefault}{\updefault}{\color[rgb]{0,0,0}$z_{1,f}$}%
}}}}
\put(6301,-2536){\makebox(0,0)[b]{\smash{{\SetFigFont{12}{14.4}{\rmdefault}{\mddefault}{\updefault}{\color[rgb]{0,0,0}$z_{2,f}$}%
}}}}
\put(6151,-3961){\makebox(0,0)[b]{\smash{{\SetFigFont{12}{14.4}{\rmdefault}{\mddefault}{\updefault}{\color[rgb]{0,0,0}$r$}%
}}}}
\put(4351,-2761){\makebox(0,0)[b]{\smash{{\SetFigFont{12}{14.4}{\rmdefault}{\mddefault}{\updefault}{\color[rgb]{0,0,0}$r$}%
}}}}
\put(3301,-3886){\makebox(0,0)[b]{\smash{{\SetFigFont{12}{14.4}{\rmdefault}{\mddefault}{\updefault}{\color[rgb]{0,0,0}$r$}%
}}}}
\put(4801,-1636){\makebox(0,0)[b]{\smash{{\SetFigFont{12}{14.4}{\rmdefault}{\mddefault}{\updefault}{\color[rgb]{0,0,0}$A$}%
}}}}
\end{picture}%